%% @texfile{
%%     filename="amssym.def",
%%     version="2.1",
%%     date="5-APR-1991",
%%     filetype="TeX: option",
%%     copyright="Copyright (C) American Mathematical Society,
%%            all rights reserved.  Copying of this file is
%%            authorized only if either:
%%            (1) you make absolutely no changes to your copy
%%                including name; OR
%%            (2) if you do make changes, you first rename it to some
%%                other name.",
%%     author="American Mathematical Society",
%%     address="American Mathematical Society,
%%            Technical Support Group,
%%            P. O. Box 6248,
%%            Providence, RI 02940,
%%            USA",
%%     telephone="401-455-4080 or (in the USA) 800-321-4AMS",
%%     email="Internet: Tech-Support@Math.AMS.com",
%%     codetable="ISO/ASCII",
%%     checksumtype="line count",
%%     checksum="108",
%%     keywords="amsfonts, tex",
%%     abstract="This file contains definitions that perform the same
%%            functions as similar ones in AMS-TeX, so that the file
%%            AMSSYM.TEX can be used outside of AMS-TeX. Instructions
%%            for using this file and the AMS symbol fonts are
%%            included in the AMSFonts 2.0 User's Guide."
%%     }
%
%%%%%%%%%%%%%%%%%%%%%%%%%%%%%%%%%%%%%%%%%%%%%%%%%%%%%%%%%%%%%%%%%%%%%%%%
\expandafter\ifx\csname amssym.def\endcsname\relax \else\endinput\fi
%
%  Store the catcode of the @ in the csname so that it can be restored later.
\expandafter\edef\csname amssym.def\endcsname{%
       \catcode`\noexpand\@=\the\catcode`\@\space}
%  Set the catcode to 11 for use in private control sequence names.
\catcode`\@=11
%
%  Include all definitions related to the fonts msam, msbm and eufm, so that
%  when this file is used by itself, the results with respect to those fonts
%  are equivalent to what they would have been using AMS-TeX.
%  Most symbols in fonts msam and msbm are defined using \newsymbol;
%  however, a few symbols that replace composites defined in plain must be
%  defined with \mathchardef.

\def\undefine#1{\let#1\undefined}
\def\newsymbol#1#2#3#4#5{\let\next@\relax
 \ifnum#2=\@ne\let\next@\msafam@\else
 \ifnum#2=\tw@\let\next@\msbfam@\fi\fi
 \mathchardef#1="#3\next@#4#5}
\def\mathhexbox@#1#2#3{\relax
 \ifmmode\mathpalette{}{\m@th\mathchar"#1#2#3}%
 \else\leavevmode\hbox{$\m@th\mathchar"#1#2#3$}\fi}
\def\hexnumber@#1{\ifcase#1 0\or 1\or 2\or 3\or 4\or 5\or 6\or 7\or 8\or
 9\or A\or B\or C\or D\or E\or F\fi}

\font\tenmsa=msam10
\font\sevenmsa=msam7
\font\fivemsa=msam5
\newfam\msafam
\textfont\msafam=\tenmsa
\scriptfont\msafam=\sevenmsa
\scriptscriptfont\msafam=\fivemsa
\edef\msafam@{\hexnumber@\msafam}
\mathchardef\dabar@"0\msafam@39
\def\dashrightarrow{\mathrel{\dabar@\dabar@\mathchar"0\msafam@4B}}
\def\dashleftarrow{\mathrel{\mathchar"0\msafam@4C\dabar@\dabar@}}

\def\ulcorner{\delimiter"4\msafam@70\msafam@70 }
\def\urcorner{\delimiter"5\msafam@71\msafam@71 }
\def\llcorner{\delimiter"4\msafam@78\msafam@78 }
\def\lrcorner{\delimiter"5\msafam@79\msafam@79 }
\def\yen{{\mathhexbox@\msafam@55 }}
\def\checkmark{{\mathhexbox@\msafam@58 }}
\def\circledR{{\mathhexbox@\msafam@72 }}
\def\maltese{{\mathhexbox@\msafam@7A }}

\font\tenmsb=msbm10
\font\sevenmsb=msbm7
\font\fivemsb=msbm5
\newfam\msbfam
\textfont\msbfam=\tenmsb
\scriptfont\msbfam=\sevenmsb
\scriptscriptfont\msbfam=\fivemsb
\edef\msbfam@{\hexnumber@\msbfam}
\def\Bbb#1{{\fam\msbfam\relax#1}}
\def\widehat#1{\setbox\z@\hbox{$\m@th#1$}%
 \ifdim\wd\z@>\tw@ em\mathaccent"0\msbfam@5B{#1}%
 \else\mathaccent"0362{#1}\fi}
\def\widetilde#1{\setbox\z@\hbox{$\m@th#1$}%
 \ifdim\wd\z@>\tw@ em\mathaccent"0\msbfam@5D{#1}%
 \else\mathaccent"0365{#1}\fi}
\font\teneufm=eufm10
\font\seveneufm=eufm7
\font\fiveeufm=eufm5
\newfam\eufmfam
\textfont\eufmfam=\teneufm
\scriptfont\eufmfam=\seveneufm
\scriptscriptfont\eufmfam=\fiveeufm

%  Restore the catcode value for @ that was previously saved.
\csname amssym.def\endcsname
\expandafter\ifx\csname pre amssym.tex at\endcsname\relax \else \endinput\fi
%%     Otherwise we store the catcode of the @ in the csname.
\expandafter\chardef\csname pre amssym.tex at\endcsname=\the\catcode`\@
%%     Set the catcode to 11 for use in private control sequence
%%     names.
\catcode`\@=11
%  Most symbols in fonts msam and msbm are defined using \newsymbol.  A few
%  that are delimiters or otherwise require special treatment have already
%  been defined as soon as the fonts were loaded.  Finally, a few symbols
%  that replace composites defined in plain must be undefined first.
\newsymbol\boxdot 1200
\newsymbol\boxplus 1201
\newsymbol\boxtimes 1202
\newsymbol\square 1003
\newsymbol\blacksquare 1004
\newsymbol\centerdot 1205
\newsymbol\lozenge 1006
\newsymbol\blacklozenge 1007
\newsymbol\circlearrowright 1308
\newsymbol\circlearrowleft 1309
\undefine\rightleftharpoons
\newsymbol\rightleftharpoons 130A
\newsymbol\leftrightharpoons 130B
\newsymbol\boxminus 120C
\newsymbol\Vdash 130D
\newsymbol\Vvdash 130E
\newsymbol\vDash 130F
\newsymbol\twoheadrightarrow 1310
\newsymbol\twoheadleftarrow 1311
\newsymbol\leftleftarrows 1312
\newsymbol\rightrightarrows 1313
\newsymbol\upuparrows 1314
\newsymbol\downdownarrows 1315
\newsymbol\upharpoonright 1316
 \let\restriction\upharpoonright
\newsymbol\downharpoonright 1317
\newsymbol\upharpoonleft 1318
\newsymbol\downharpoonleft 1319
\newsymbol\rightarrowtail 131A
\newsymbol\leftarrowtail 131B
\newsymbol\leftrightarrows 131C
\newsymbol\rightleftarrows 131D
\newsymbol\Lsh 131E
\newsymbol\Rsh 131F
\newsymbol\rightsquigarrow 1320
\newsymbol\leftrightsquigarrow 1321
\newsymbol\looparrowleft 1322
\newsymbol\looparrowright 1323
\newsymbol\circeq 1324
\newsymbol\succsim 1325
\newsymbol\gtrsim 1326
\newsymbol\gtrapprox 1327
\newsymbol\multimap 1328
\newsymbol\therefore 1329
\newsymbol\because 132A
\newsymbol\doteqdot 132B
 
\newsymbol\triangleq 132C
\newsymbol\precsim 132D
\newsymbol\lesssim 132E
\newsymbol\lessapprox 132F
\newsymbol\eqslantless 1330
\newsymbol\eqslantgtr 1331
\newsymbol\curlyeqprec 1332
\newsymbol\curlyeqsucc 1333
\newsymbol\preccurlyeq 1334
\newsymbol\leqq 1335
\newsymbol\leqslant 1336
\newsymbol\lessgtr 1337
\newsymbol\backprime 1038
\newsymbol\risingdotseq 133A
\newsymbol\fallingdotseq 133B
\newsymbol\succcurlyeq 133C
\newsymbol\geqq 133D
\newsymbol\geqslant 133E
\newsymbol\gtrless 133F
\newsymbol\sqsubset 1340
\newsymbol\sqsupset 1341
\newsymbol\vartriangleright 1342
\newsymbol\vartriangleleft 1343
\newsymbol\trianglerighteq 1344
\newsymbol\trianglelefteq 1345
\newsymbol\bigstar 1046
\newsymbol\between 1347
\newsymbol\blacktriangledown 1048
\newsymbol\blacktriangleright 1349
\newsymbol\blacktriangleleft 134A
\newsymbol\vartriangle 134D
\newsymbol\blacktriangle 104E
\newsymbol\triangledown 104F
\newsymbol\eqcirc 1350
\newsymbol\lesseqgtr 1351
\newsymbol\gtreqless 1352
\newsymbol\lesseqqgtr 1353
\newsymbol\gtreqqless 1354
\newsymbol\Rrightarrow 1356
\newsymbol\Lleftarrow 1357
\newsymbol\veebar 1259
\newsymbol\barwedge 125A
\newsymbol\doublebarwedge 125B
\undefine\angle
\newsymbol\angle 105C
\newsymbol\measuredangle 105D
\newsymbol\sphericalangle 105E
\newsymbol\varpropto 135F
\newsymbol\smallsmile 1360
\newsymbol\smallfrown 1361
\newsymbol\Subset 1362
\newsymbol\Supset 1363
\newsymbol\Cup 1264
 
\newsymbol\Cap 1265
 
\newsymbol\curlywedge 1266
\newsymbol\curlyvee 1267
\newsymbol\leftthreetimes 1268
\newsymbol\rightthreetimes 1269
\newsymbol\subseteqq 136A
\newsymbol\supseteqq 136B
\newsymbol\bumpeq 136C
\newsymbol\Bumpeq 136D
\newsymbol\lll 136E
 
\newsymbol\ggg 136F
 
\newsymbol\circledS 1073
\newsymbol\pitchfork 1374
\newsymbol\dotplus 1275
\newsymbol\backsim 1376
\newsymbol\backsimeq 1377
\newsymbol\complement 107B
\newsymbol\intercal 127C
\newsymbol\circledcirc 127D
\newsymbol\circledast 127E
\newsymbol\circleddash 127F
\newsymbol\lvertneqq 2300
\newsymbol\gvertneqq 2301
\newsymbol\nleq 2302
\newsymbol\ngeq 2303
\newsymbol\nless 2304
\newsymbol\ngtr 2305
\newsymbol\nprec 2306
\newsymbol\nsucc 2307
\newsymbol\lneqq 2308
\newsymbol\gneqq 2309
\newsymbol\nleqslant 230A
\newsymbol\ngeqslant 230B
\newsymbol\lneq 230C
\newsymbol\gneq 230D
\newsymbol\npreceq 230E
\newsymbol\nsucceq 230F
\newsymbol\precnsim 2310
\newsymbol\succnsim 2311
\newsymbol\lnsim 2312
\newsymbol\gnsim 2313
\newsymbol\nleqq 2314
\newsymbol\ngeqq 2315
\newsymbol\precneqq 2316
\newsymbol\succneqq 2317
\newsymbol\precnapprox 2318
\newsymbol\succnapprox 2319
\newsymbol\lnapprox 231A
\newsymbol\gnapprox 231B
\newsymbol\nsim 231C
\newsymbol\ncong 231D
\newsymbol\diagup 231E
\newsymbol\diagdown 231F
\newsymbol\varsubsetneq 2320
\newsymbol\varsupsetneq 2321
\newsymbol\nsubseteqq 2322
\newsymbol\nsupseteqq 2323
\newsymbol\subsetneqq 2324
\newsymbol\supsetneqq 2325
\newsymbol\varsubsetneqq 2326
\newsymbol\varsupsetneqq 2327
\newsymbol\subsetneq 2328
\newsymbol\supsetneq 2329
\newsymbol\nsubseteq 232A
\newsymbol\nsupseteq 232B
\newsymbol\nparallel 232C
\newsymbol\nmid 232D
\newsymbol\nshortmid 232E
\newsymbol\nshortparallel 232F
\newsymbol\nvdash 2330
\newsymbol\nVdash 2331
\newsymbol\nvDash 2332
\newsymbol\nVDash 2333
\newsymbol\ntrianglerighteq 2334
\newsymbol\ntrianglelefteq 2335
\newsymbol\ntriangleleft 2336
\newsymbol\ntriangleright 2337
\newsymbol\nleftarrow 2338
\newsymbol\nrightarrow 2339
\newsymbol\nLeftarrow 233A
\newsymbol\nRightarrow 233B
\newsymbol\nLeftrightarrow 233C
\newsymbol\nleftrightarrow 233D
\newsymbol\divideontimes 223E
\newsymbol\varnothing 203F
\newsymbol\nexists 2040
\newsymbol\Finv 2060
\newsymbol\Game 2061
\newsymbol\mho 2066
\newsymbol\eth 2067
\newsymbol\eqsim 2368
\newsymbol\beth 2069
\newsymbol\gimel 206A
\newsymbol\daleth 206B
\newsymbol\lessdot 236C
\newsymbol\gtrdot 236D
\newsymbol\ltimes 226E
\newsymbol\rtimes 226F
\newsymbol\shortmid 2370
\newsymbol\shortparallel 2371
\newsymbol\smallsetminus 2272
\newsymbol\thicksim 2373
\newsymbol\thickapprox 2374
\newsymbol\approxeq 2375
\newsymbol\succapprox 2376
\newsymbol\precapprox 2377
\newsymbol\curvearrowleft 2378
\newsymbol\curvearrowright 2379
\newsymbol\digamma 207A
\newsymbol\varkappa 207B
\newsymbol\Bbbk 207C
\newsymbol\hslash 207D
\undefine\hbar
\newsymbol\hbar 207E
\newsymbol\backepsilon 237F
%  Restore the catcode value for @ that was previously saved.
\catcode`\@=\csname pre amssym.tex at\endcsname
%\input mssymb      % these must be input at the BEGINNING of this file; some
% MAXIMUM LINE LENGTH IS 80 CHARACTERS
%
%\font\eurofont=times at 10pt% contains � symbol (as option-shift-2)
%
\font\fivebi=cmmib5
\font\fivebsy=cmbsy5
\font\sixrm=cmr6
\font\sixi=cmmi6
\font\sixbf=cmbx6
\font\sixsy=cmsy6
\font\sixmsa=msam5 at 6pt% This is needed since Bluesky's virtual fonts for
\font\sixmsb=msbm5 at 6pt% sizes 6.8.9 all come out too small (same below)
\font\sevenbi=cmmib7
\font\sevenbsy=cmbsy7
\font\eightrm=cmr8
\font\eightsl=cmsl8
\font\eightit=cmti8
\font\eighti=cmmi8
\font\eightbf=cmbx8
\font\eightsy=cmsy8
\font\eightmsa=msam7 at 8pt
\font\eightmsb=msbm7 at 8pt
\font\ninerm=cmr9
\font\ninesl=cmsl9
\font\nineit=cmti9
\font\ninei=cmmi9
\font\ninebf=cmbx9
\font\ninebi=cmmib10 scaled 900
\font\ninesy=cmsy9
\font\ninebsy=cmbsy10 scaled 900
\font\ninemsa=msam10 at 9pt
\font\ninemsb=msbm10 at 9pt
\font\tenbit=cmbxti10
\font\tenbsl=cmbxsl10
\font\tenbi=cmmib10
\font\tenbsy=cmbsy10
\font\twelvebf=cmbx12
\font\twelvebi=cmmib10 scaled 1200
\font\twelvebsy=cmbsy10 at 12pt

\let\sc=\sevenrm            % SMALL CAPS (in tenpoint)
\def\eightpoint{%
    % \font\eurofont=times at 8pt% contains � symbol (as option-shift-2)
     \def\rm{\fam0\eightrm}%         see p.415
     \textfont0=\eightrm \scriptfont0=\sixrm \scriptscriptfont0=\fiverm
     \textfont1=\eighti \scriptfont1=\sixi \scriptscriptfont1=\fivei
     \textfont2=\eightsy \scriptfont2=\sixsy \scriptscriptfont2=\fivesy
     \textfont3=\tenex \scriptfont3=\tenex \scriptscriptfont3=\tenex
     \textfont\itfam=\eightit \def\it{\fam\itfam\eightit}%
     \textfont\slfam=\eightsl \def\sl{\fam\slfam\eightsl}%
     \textfont\bffam=\eightbf \scriptfont\bffam=\sixbf
     \scriptscriptfont\bffam=\fivebf \def\bf{\fam\bffam\eightbf}%
     \textfont\msbfam=\eightmsb \textfont\msafam=\eightmsa
     \scriptfont\msafam=\sixmsa \scriptfont\msbfam=\sixmsb
     \scriptscriptfont\msafam=\fivemsa \scriptscriptfont\msbfam=\fivemsb
      \skewchar\eighti='177 \skewchar\sixi='177
      \skewchar\eightsy='60 \skewchar\sixsy='60
     \normalbaselineskip=10pt%(normally 9pt; must be =height+depth in next line)
     \setbox\strutbox=\hbox{\vrule height7pt depth3pt width0pt}%(normally 7&2pt)
     \let\sc=\sixrm \let\big=\eightbig \normalbaselines\rm}
\def\ninepoint{%
    % \font\eurofont=times at 9pt% contains � symbol (as option-shift-2)
     \def\rm{\fam0\ninerm}%         see p.415
     \textfont0=\ninerm \scriptfont0=\sixrm \scriptscriptfont0=\fiverm
     \textfont1=\ninei \scriptfont1=\sixi \scriptscriptfont1=\fivei
     \textfont2=\ninesy \scriptfont2=\sixsy \scriptscriptfont2=\fivesy
     \textfont3=\tenex \scriptfont3=\tenex \scriptscriptfont3=\tenex
     \textfont\itfam=\nineit \def\it{\fam\itfam\nineit}%
     \textfont\slfam=\ninesl \def\sl{\fam\slfam\ninesl}%
     \textfont\bffam=\ninebf \scriptfont\bffam=\sixbf
     \scriptscriptfont\bffam=\fivebf \def\bf{\fam\bffam\ninebf}%
     \textfont\msbfam=\ninemsb \textfont\msafam=\ninemsa
     \scriptfont\msafam=\sixmsa \scriptfont\msbfam=\sixmsb
     \scriptscriptfont\msafam=\fivemsa \scriptscriptfont\msbfam=\fivemsb
      \skewchar\ninei='177 \skewchar\sixi='177
      \skewchar\ninesy='60 \skewchar\sixsy='60
     \normalbaselineskip=11pt
     \setbox\strutbox=\hbox{\vrule height8pt depth3pt width0pt}%
     \let\sc=\sevenrm \let\big=\ninebig \normalbaselines\rm}
\def\tenpoint{%
    % \font\eurofont=times at 10pt% contains � symbol (as option-shift-2)
     \def\rm{\fam0\tenrm}%         see p.415
     \textfont0=\tenrm \scriptfont0=\sevenrm \scriptscriptfont0=\fiverm
     \textfont1=\teni \scriptfont1=\seveni \scriptscriptfont1=\fivei
     \textfont2=\tensy \scriptfont2=\sevensy \scriptscriptfont2=\fivesy
     \textfont3=\tenex \scriptfont3=\tenex \scriptscriptfont3=\tenex
     \textfont\itfam=\tenit \def\it{\fam\itfam\tenit}%
     \textfont\slfam=\tensl \def\sl{\fam\slfam\tensl}%
     \textfont\bffam=\tenbf \scriptfont\bffam=\sevenbf
     \scriptscriptfont\bffam=\fivebf \def\bf{\fam\bffam\tenbf}%
     \textfont\msbfam=\tenmsb \textfont\msafam=\tenmsa
     \scriptfont\msafam=\sevenmsa \scriptfont\msbfam=\sevenmsb
     \scriptscriptfont\msafam=\fivemsa \scriptscriptfont\msbfam=\fivemsb
     \normalbaselineskip=12pt
%     \normalbaselineskip=13pt% (normally 12pt; here 13 = 9 + 4 of next line)
%     \setbox\strutbox=\hbox{\vrule height9pt depth4pt width0pt}%
     \let\sc=\sevenrm \let\big=\tenbig \normalbaselines\rm}
\catcode`@=11        % This allows the use of `@' in the next line (p.344)
\def\tenbig#1{{\hbox{$\left#1\vbox to8.5pt{}\right.\n@space$}}}
\def\ninebig#1{{\hbox{$\textfont0=\tenrm\textfont2=\tensy
      \left#1\vbox to7.25pt{}\right.\n@space$}}}
\def\eightbig#1{{\hbox{$\textfont0=\ninerm\textfont2=\ninesy
      \left#1\vbox to6.5pt{}\right.\n@space$}}}
\catcode`@=12        % This restores the `inhibiting' catcode of `@'.
\def\bold{%
     \textfont0=\tenbf \scriptfont0=\sevenbf \scriptscriptfont0=\fivebf
     \textfont1=\tenbi \scriptfont1=\sevenbi \scriptscriptfont1=\fivebi
     \textfont2=\tenbsy \scriptfont2=\sevenbsy \scriptscriptfont2=\fivebsy
       \textfont\itfam=\tenbit \def\it{\fam\itfam\tenbit}%
       \textfont\slfam=\tenbsl \def\sl{\fam\slfam\tenbsl}%
       \textfont\bffam=\tenbf \scriptfont\bffam=\sevenbf
       \textfont\msbfam=\tenmsb \textfont\msafam=\tenmsa
   \fam0\tenbf}
\def\bigbold{%
     \textfont0=\twelvebf \scriptfont0=\ninebf \scriptscriptfont0=\sevenbf
     \textfont1=\twelvebi \scriptfont1=\ninebi \scriptscriptfont1=\sevenbi
     \textfont2=\twelvebsy \scriptfont2=\ninebsy \scriptscriptfont2=\sevenbsy
     \fam0\twelvebf}
%
%

%
% MAXIMUM LINE LENGTH BELOW IS 80 CHARACTERS
%
% Instructions on how to use \item, and how to surround items by space,
% are given in paperformat.
%
\let\plainitem=\item% \item itself is redefined in paperformat.
\let\plainitemitem=\itemitem% \itemitem itself is redefined in paperformat.
%

%
% � is option-shift-2 on the Mac
%
\def\openface{\Bbb}                %  needs mssymb
\def\N{{\openface N}}              %  in old files, \nat may be used instead

%
%
%                            SKIPS AND BREAKS
%
\def\g{\hskip.17em\relax}               %  breakable thin space
\def\th{\thinspace}                     %  non-breakable thin space
\def\nl{\hfil\break}
\newskip\Bigskipamount
   \Bigskipamount=2\baselineskip plus.5\baselineskip minus.3\baselineskip
\def\Bigbreak{\removelastskip\vskip0pt plus .1\vsize\penalty-1000
              \vskip0pt plus-.1\vsize\vskip\Bigskipamount}
% The following can be used after a display to remove the belowdisplayskip, e.g.
% when the display finishes a proof and is just followed by \endproof. If text
% follows, e.g. just a word like "and" leading to another display, one may wish
% to say \noskip\smallskip\noindent.
\def\noskip{\vskip-\lastskip\noindent}%use after display followed by short line
\def\Nobreak$$#1$${\postdisplaypenalty=10000$$#1$$\postdisplaypenalty=0}
%
%
%                             ABBREVIATIONS
%
     % \H will redefined below.

\let\doublebar=\| 
\def\|{\!\!\restriction\!\!}

  \let\sub=\sube

\def\supe{\supseteq}

\def\sm{\smallsetminus}
\def\es{\emptyset}

 % italics, like K^n, P^n etc. Cf. I FOR "Inflated"

\def\:{\colon}
\def\minor{\preccurlyeq} 
\def\Minor{\succcurlyeq}

\def\slt{\mathrel{\hbox{$\minor$\kern-.6em\lower.33ex\hbox{${}_s\;$}}}}
\def\sgt{\mathrel{\mathchoice                        %("simplicial minors")
   {\hbox{$\Minor$\kern-.5em\lower.3ex\hbox{${}_s$}}}
   {\hbox{$\Minor$\kern-.5em\lower.3ex\hbox{${}_s$}}}
   {\hbox{$\scriptstyle\Minor\kern-.43em\lower.28ex\hbox{$\scriptstyle{}_s$}$}}
 {\hbox{$\scriptstyle\Minor\kern-.43em\lower.28ex\hbox{$\scriptstyle{}_s$}$}} }}    

\def\ucl(#1){\lfloor #1 \rfloor}% up-closure
\def\dcl(#1){\lceil #1 \rceil}% down-closure
%
%

%

% Example: '\interior P' puts a circle over the P.
% (the "\relax" is important: otherwise a capital letter A-F instead of P
% is interpreted as part of the number 7017,
% giving a "wrong math code complaint!!)
%
% The following commands can be used to 'specify' a relation (= #1) by putting
% something (= #2) underneath or above.
% Examples:\specrel<T, \specrel\sim G, \specrel={(1)}, \Specrel\Rightarrow?.
%
\def\specrel#1#2{\mathrel{\mathop{\kern0pt #1}\limits_{#2}}}
\def\Specrel#1#2{\mathrel{\mathop{\kern0pt #1}\limits^{#2}}}
%
% Next a version of \specrel for use in aligned equations;
% here the specification text is put in an \hbox of width 0,
% so as not to interfere with the alignment.
%
\def\alignspecrel#1#2{\mathrel{\mathop{\kern0pt #1}\limits_{\hbox
   to0pt{\hss$\scriptstyle#2$\hss}}}}
\def\alignSpecrel#1#2{\mathrel{\mathop{\kern0pt #1}\limits^{\hbox
   to0pt{\hss$\scriptstyle#2$\hss}}}}
\def\invlim{\specrel\lim{\raise 2pt\hbox{$\longleftarrow$}}}% inverse limit, projective limit
%
%

%

% Syntax example for the following: \looseproof{Zweiter Beweis von Satz \xxxA}
%space follows in input

%
\newcount\refno
\def\ref#1#2\par{{\plainitem{[??]}#2\smallskip}}
\newtoks\thingtowrite %USED AGAIN LATER FOR INDEX
\long\def\writerefnumber#1{%
    \thingtowrite={#1}%
    \immediate\write\refnumbersfile{\the\thingtowrite}%
    }
\newwrite\refnumbersfile
\def\makerefnumbers{\immediate\openout\refnumbersfile=RefNumbers%
  \refno=0 \writerefnumber{\refno=0}
  \def\ref##1##2\par{\global\advance\refno by 1
    \writerefnumber{\global\advance\refno by 1 \newcounter##1 ##1=\the\refno}%
       % \newcounter replaces \newcount, which is forbidden inside a def.
       % When using the auto-generated file, say "\let\newcounter=\newcount"
       % before reading in that file.
    \plainitem{[\the\refno]}##2\smallskip% automatic numbering, ignoring ##1
    }%
  }
\def\autorefnumbers{\refno=0
  \def\ref##1##2\par{\advance\refno by 1\plainitem{[\the\refno]}##2\smallskip}
  }
\def\userefnumbers{\refno=0
  \def\ref##1##2\par{\advance\refno by 1\plainitem{[\the##1]}##2\smallskip}
  }
% TO USE, FIRST RUN WITH \makerefnumbers active (and no file named "RefNumbers" present
% in the directory to which such a file is to be written; it cannot overwrite).
% THEN DISABLE "\makerefnumbers", and make the following line active:
% \userefnumbers\let\newcounter=\newcount\input RefNumbers
%
%
\def\proclaimwithname #1. (#2) #3\par{{\bigbreak
  \clubpenalty=10000\noindent{\bf#1.\enspace}(#2)\nl
  {\sl #3}\par\bigbreak}}
\def\proposition (#1) #2\par{{\setbox0\hbox{(#1)\enspace}\bigbreak
   \sl\hangindent\the\wd0 \noindent\hskip\the\wd0
   \llap{\box0}\ignorespaces#2\par\bigbreak}}
\def\subsection #1\par{\vskip 2\medskipamount plus \smallskipamount minus \smallskipamount\leftline{\bold #1}
        \penalty10000\smallskip\noindent}
\def\longsubsection #1 #2\par{\vskip 2\medskipamount plus \smallskipamount minus \smallskipamount\varitem{\bold #1} {\bold  #2}
        \penalty10000\smallskip\noindent}
      \def\section #1\par{\Bigbreak\centerline{\bf #1} % TO BE PHASED OUT
              \penalty10000\bigskip\noindent}
%
%         The following versions of the beginsection macro are for
%         section headings immediately followed by \proclaim:
\def\beginpsection #1\par{\Bigbreak\centerline{\bold #1}
        \penalty10000\bigskip}
\def\psubsection #1\par{\bigbreak\leftline{\bold #1}\penalty10000\bigskip}
%
%         The following macro positions its argument flush right,
%         like the \endproof box or an equation number.

%
%  The following three items usually need an \enditem. They are for use inside
%  \proclaim or another \item. Note that these items don't require their own
%  \par at the end, but they have no finishing \smallskip. (There seems to be
%  no way around this: we can't delimit by \par, since in the input this would
%  end the \proclaim.) So, if the item is followed by some text within the same
% \proclaim (say), the \enditem should (and can) be replaced with \smallskip.
%        
\def\pitem#1{\smallskip\advance\parindent by 3mm
             \plainitem{\rm(#1)}\advance\parindent by-3mm}
\def\pitemitem#1{\smallskip\advance\parindent by 3mm
             \plainitemitem{\rm(#1)}\advance\parindent by-3mm}
\def\varitemitem#1{{\setbox0\hbox{\hskip\parindent#1\enskip}
           \smallbreak\hangindent\the\wd0 \noindent\hskip\the\wd0
           \llap{#1\enskip}\ignorespaces}}
%
%        Note that the following variations of \item must end with \par.
%
\newdimen\newparindent
\def\iitem#1#2\par{\newparindent=\parindent \advance\newparindent by 3mm
           \smallbreak \hangindent\newparindent \noindent\hskip\newparindent
           \llap{{\rm #1}\enspace}\ignorespaces#2\par\smallbreak}
\def\iitemitem#1#2\par{\newparindent=\parindent \advance\newparindent by 3mm
           \smallbreak \hangindent2\newparindent \noindent\hskip2\newparindent
           \llap{{\rm #1}\enspace}\ignorespaces#2\par\smallbreak}
\def\varitem#1#2\par{{\setbox0\hbox{{\rm #1}\enspace}
           \smallbreak \hangindent\the\wd0 \noindent\hskip\the\wd0
           \llap{{\rm #1}\enspace}\ignorespaces#2\par\smallbreak}}

% For use in Exercises and Hints
%
\def\Textindent#1{\par \advance\parindent by 3mm
                  \textindent{{\rm #1}} \advance\parindent by -3mm}
\def\indentedline#1{\advance\hsize by -\parindent \line{#1}
                   \advance\hsize by \parindent}
\def\iindentedline#1{\advance\parindent by 3mm
                     \advance\hsize by -\parindent
                     \line{#1}
                     \advance\hsize by \parindent
                     \advance\parindent by -3mm}
\newdimen\margin   % needed for macros \textdisplay & \ltextdisplay
%  The following macro takes 3 arguments, #1 and #3 in math-mode,
%  #2 as plain text. It displays #1\quad centered w.r.t. the whole \line
%  (if #2 leaves enough space), adds \quad#2 to the right of #1\quad, and puts
%  #3 flush right. The arguments must be separated by &'s. The argu-
%  ments themselves may be empty, but there must be 2 ampersands.
%  Examples: $$\textdisplay x=y &for all $x\in X$& (1')$$
%            $$\textdisplay x=y &for all $x\in X$&$$
\def\textdisplay#1&#2&#3$${\margin=\hsize
          \setbox1=\hbox{$\displaystyle#1\quad$}%
          \setbox2=\hbox{\quad#2\qquad$#3$}%
                     \advance\margin by-\wd1
                     \divide\margin by 2
   \ifdim\wd2 < \margin
      \box1\rlap{\quad#2}\eqno#3$$%
   \else
      \line{\qquad\hfil \box1\quad #2 \qquad $#3$}$$%
   \fi}
%
%  The following macro is the `\leqno' version of \textdisplay; argument
%  #1 is the \leqno, #2 is the formula to be displayed, and #3 is the
%  text following the formula.
\def\ltextdisplay#1&#2&#3$${\margin=\hsize
           \setbox2=\hbox{$\displaystyle#2\quad$}
           \setbox3=\hbox{\quad#3\qquad}
                     \advance\margin by-\wd2
                     \divide\margin by 2
   \ifdim\wd3 < \margin
      \line{$#1$\hfil\box2\hbox to \margin{\box3\hfil}}$$%
   \else
      \line{$#1$\qquad\hfil\box2\quad #3\qquad} $$%
   \fi}
%
%   The next macro displays and centres #1, typically a paragraph of text.
%   #2, typically an `eqno', is set flush right, vertically centred, and
%   has to be in math mode.
%   Example: \textno This is a ... lot of text&(3')\par
\def\textno#1&#2\par{%
   \margin=\hsize
   \advance\margin by -4\parindent
          \setbox1=\hbox{\sl#1}%
   \ifdim\wd1 < \margin
      $$\box1\eqno#2$$\endgraf%
   \else
      \bigbreak
      \line{\indent$\vcenter{\advance\hsize by -3\parindent
      \sl\noindent#1}\hfil#2$}%
      \bigbreak
   \fi}
%
%   The same with a left eqno:
\def\textlno#1&#2\par{%
   \margin=\hsize
   \advance\margin by -4\parindent
          \setbox1=\hbox{\sl#1}%
   \ifdim\wd1 < \margin
      $$\box1\leqno#2$$%
   \else
      \bigbreak
      \line{$#2\hfil\vcenter{\advance\hsize by -3\parindent
          \sl\noindent#1}\hskip\parindent$}%
      \bigbreak
   \fi}
%
% For use inside environments (such as \items, \solutions etc) that are delimited by \par:
%
\def\textnoNonPar#1&#2\endgraf{%
   \margin=\hsize
   \advance\margin by -4\parindent
          \setbox1=\hbox{\sl#1}%
   \ifdim\wd1 < \margin
      $$\box1\eqno#2$$\endg�raf%
   \else
      \bigbreak
      \line{\indent$\vcenter{\advance\hsize by -3\parindent
      \sl\noindent#1}\hfil#2$}%
      \bigbreak
   \fi}
%
%
%                        MARGINAL HACKS ETC.
%
\newcount\commentno
\def\COMMENT#1{$^{<\the\commentno>}$%
     \vadjust{\vbox to 0pt{\vss\vskip-8pt\rightline{%
     \rlap{\hbox{\hskip7mm \vbox{\pretolerance=-1
     \doublehyphendemerits=0 \finalhyphendemerits=0
     \hsize40mm\tolerance=10000\eightpoint
     \lineskip=0pt\lineskiplimit=0pt
     \rightskip=0pt plus16mm\baselineskip8pt\noindent
     \hskip0pt       %(without this, the first word is never hyphenated!)
     {$\langle$\the\commentno. #1$\rangle$}\endgraf}}}}\vss}}%
     \global\advance\commentno by1}%
\def\writecommentsasfootnotes{%
 \def\COMMENT{\global\advance\commentno by1\footnote{$^{<\the\commentno>}$}}%
 }
\def\nocomments{\def\COMMENT##1{}}
%
%
%   The \? macro puts
%   the argument #1 in the left margin. Examples: \??, \?{What nonsense!}.
\def\?#1{\vadjust{\vbox to 0pt{\vss\vskip-8pt\leftline{%
     \llap{\hbox{\vbox{\pretolerance=-1
     \doublehyphendemerits=0\finalhyphendemerits=0
     \hsize16truemm\tolerance=10000\eightpoint
     \lineskip=0pt\lineskiplimit=0pt
     \rightskip=0pt plus16truemm\baselineskip8pt\noindent
     \hskip0pt        %(without this, the first word is never hyphenated!)
     #1\endgraf}\hskip7truemm}}}\vss}}}
\def\d{}% DISABLES WHAT'S JUST BEEN DEFINED ABOVE
%
% The following \ds macro is the "silent" version of \d: it writes its argument
% in the margin (or index file), just as \d does, but not into the current text.
% Note the different syntax: no delimiting blank in input.
%
%  CHANGE: I got fed up with marginal reminders to define things, so this
%  feature is disabled right away by the line "\def\ds#1{}" below. Note that
%  \makeindex still works, because it defines \d and \ds anew.
%
%\def\ds#1{\ifmmode
%     \vadjust{\vbox to 0pt{\vss\vskip-8pt\leftline{%
%     \llap{\hbox{\vbox{\pretolerance=-1
%    \doublehyphendemerits=0\finalhyphendemerits=0
%     \hsize16truemm\tolerance=10000\eightpoint
%     \lineskip=0pt\lineskiplimit=0pt
%     \rightskip=0pt plus16truemm\baselineskip8pt\noindent
%     define $#1$!\endgraf}\hskip7truemm}}}\vss}}%
%   \else
%     \vadjust{\vbox to 0pt{\vss\vskip-8pt\leftline{%
%     \llap{\hbox{\vbox{\pretolerance=-1
%     \doublehyphendemerits=0\finalhyphendemerits=0
%     \hsize16truemm\tolerance=10000\eightpoint
%     \lineskip=0pt\lineskiplimit=0pt
%     \rightskip=0pt plus16truemm\baselineskip8pt\noindent
%     \hskip0pt        %(without this, the first word is never hyphenated!)
%     define ``#1''!\endgraf}\hskip7truemm}}}\vss}}%
%   \fi}
%
\def\ds#1{}% DISABLES WHAT'S JUST BEEN DEFINED ABOVE
%
% The following is a device (due to CET1) which allows to \write something to
% a file without expanding all the tokens completely. (This would result in
% the use of lots of "@"s, which make the file untexable. Another way around the
% problem would be to precede every argument of \write by "\catcode`@=11" (and
% to set it back to 12 after the argument). That would make the file texable,
% but since it would still be difficult to read, the method below is better.
% NOTE: This works well for a Remarks file (say), but NOT for an index file. 
% The reason is that in order to get the page numbers of entries from the first
% paragraph of a page right, one has to say "\write" rather than
% "\immediate\write". But in the def of \indexwrite one has to say "\immediate",
% since otherwise the token gets overwritten and rather
% than n index words from a page the index will
% contain the last index word of that page n times. But if the index words are
% written to file immediately while their page numbers get queued, then the two
% get separated from one another.
%
%\newtoks\thingtowrite % NOW EARLIER
\long\def\indexwrite#1{%
    \thingtowrite={#1}%
    \immediate\write\index{\the\thingtowrite}%
    }
%
%  The following macro \makeindex redefines \d and \ds.
%  It suppresses the appearance of \d's arguments in the margin
%  and writes them to a file called "index" instead.
%
\newwrite\index
\def\makeindex{\immediate\openout\index=index%
   \immediate\write\index{\catcode`@=11}%
   \def\d##1 {\ifmmode
     \write\index{$##1$, }%
     \write\index{\the\count0}\write\index{}% blank line for para
   \else
     \write\index{{##1}, }%
     \write\index{\the\count0}\write\index{}% blank line for para
   \fi {##1} }% It's important to put the text AFTER the def: otherwise
              % the blank is not merged with possible blanks following in
              % the input file, which may result in an additional blank line
      \def\ds##1{\ifmmode
     \write\index{$##1$, }%
     \write\index{\the\count0}\write\index{}% blank line for para
   \else
     \write\index{##1, }%
     \write\index{\the\count0}\write\index{}% blank line for para
   \fi}}
%
%   The \m macro puts its (direct!) argument in the right margin,
%   leaving it also in the text (in italics). 
%   Syntax example: this \m{word\/} is defined here
%   and will be shown in the margin; similarly the set $X$ in $A =: \m X$.
%   Intended use: in preprints, at places where something is first defined
%   (to enhance readability of the paper, in the absence of an index).
%   To disable, say "\disablems" at beginning of file.
\newdimen\gap% gap between text and hack
\gap=3truemm
\newdimen\hackwidth
\hackwidth=15truemm
\def\disablems{\def\mos##1{\strut}}% NOT simply "{}": when \m is active
 % it puts a strut in the line, so this should be here also when \m is disabled
 % to avoid a change in vertical space (and hence possibly in page breaks)
% \ds writes to index file, \mo in margin and text
\def\mo#1{\ifmmode {#1}\else {\it#1}\fi\mos{#1}}
% \ds writes to index file, \mos in margin
\def\mos#1{\ifmmode
     \strut\vadjust{\vbox to 0pt{\vss\kern-11pt\leftline{%
     \llap{\hbox{\vbox{\pretolerance=-1
     \doublehyphendemerits=0\finalhyphendemerits=0
     \hsize\hackwidth\tolerance=10000\eightpoint
     \lineskip=0pt\lineskiplimit=0pt
     \rightskip=0pt plus\hsize\baselineskip8pt\noindent
     $#1$\strut\endgraf}\hskip\gap }}}\vss}}%
   \else
     \strut\vadjust{\vbox to 0pt{\vss\kern-11pt\leftline{%
     \llap{\hbox{\vbox{\pretolerance=-1
     \doublehyphendemerits=0\finalhyphendemerits=0
     \hsize\hackwidth\tolerance=10000\eightpoint
     \lineskip=0pt\lineskiplimit=0pt
     \rightskip=0pt plus\hsize\baselineskip8pt\noindent
     \hskip0pt    %(without this, the first word is never hyphenated!)
     {\sl#1}\strut\endgraf}\hskip\gap }}}\vss}}%
   \fi}%
\newcount\remarkno
\def\REMARK#1{{\footnote{${}^{\the\remarkno}$}{{#1}}%
   \global\advance\remarkno by1}}
\def\noremarks{\def\REMARK##1{}}
%
%
%
%                              PICTURES
%
% In Textures, one can directly use
%
%    \special{illustration picfile scaled 0.7}
%
% etc, where picfile.eps is an eps input file that Textures must be able to find.
% The \picture macro below seems to be obsolete: it seems to look for a picture
% in what used to be Textures' "pictures window".
%
% With eps macros installed, one can also use \epsfbox{picfile} (with scaling as below)
%
\def\picture #1 by #2 (#3){
  \vbox to #2{
          \vfill
          \special{picture #3}
          \hrule width #1 height 0pt depth 0pt
           }}
\newdimen\topfiguremargin
   \topfiguremargin=0pt                                  % default
\newdimen\bottomfiguremargin
   \bottomfiguremargin=\medskipamount                    % default
\newdimen\normalpictureheight
\normalpictureheight=40mm
%   The following macro uses TeXtures pictures; these MUST be named
%   Fig.1, Fig.2.5 etc. (without a space after the '.'), as in the macro call
%itself. The width (#2) and height (#3) should have the original values
%of the TeXtures picture, to facilitate proper scaling. The heightfactor (#4)
%is divided by 1000 and used to scale the \normalpictureheight. Thus, height-
%factor 500 (2000) results in a picture of half (twice) the normalpicureheight.
%   Syntax example for a figure at the standard \normalpictureheight:
%\Fig.2 (538pt by 536pt; heightfactor: 1000; caption: This is a short caption)
\def\Fig.#1 (#2by#3; heightfactor:#4; caption:#5) {{%
   \dimen2=\normalpictureheight
   \dimen0=#2                          % computing width
      \divide\dimen2 by 1000
      \multiply\dimen2 by#4              % \dimen2 := intended pictureheight
   \count2=\dimen2                  % computing scalefactor
      \dimen1=#3                             % \dimen1 := actual pictureheight
   \count1=\dimen1
   \divide\count1 by 1000
   \divide\count2 by \count1          % \count2 := scalefactor (times 1000)
%      \message{scalefactor in Fig.#1 is \the\count2}%
   \divide\dimen0 by 1000
   \multiply\dimen0 by \count2      % \dimen0 := width
         \dimen1=\hsize
         \advance\dimen1 by -\dimen0
         \divide\dimen1 by 2               % \dimen1 := margin
   \midinsert
   \vbox to \topfiguremargin{\vfil}
   \noindent\hskip\dimen1
   \picture\dimen0 by \dimen2  (Fig.#1 scaled \the\count2)%
   \vskip\bottomfiguremargin                     % beginning caption
      \ninepoint
      \parindent=.1\hsize\narrower\narrower
      \setbox0\hbox{#5}
      \ifdim\wd0 < .6\hsize
           \centerline{F{\sc IGURE} #1.\hskip1em#5}
       \else
           \plainitem{F{\sc IGURE} #1. }#5\par
       \fi
   \vskip0pt\endinsert}}
%
%The following \textpicture macro is for inserting pictures in the current
%text line. The width (#2) and height (#3) should have the original values of
%the TeXtures picture, to facilitate proper scaling. #4 offers the opportunity
%to reserve vertical space for the picture by saying "height15pt depth10pt" or
%so; this will be the height of the \vbox containing the picture (default=0).
%#6 is the amount by which the bottom left corner is placed below the baseline
%(poss.neg.), #5 is the horizontal extension of the picture. Syntax examples:
%   \textpicture flower(538pt by 536pt; width5em lower-5pt)
%   \textpicture bug(538pt by 536pt; height20pt depth20pt width5em lower20pt)
\def\textpicture #1(#2by#3; #4width#5lower#6){{%
  % computing scalefactor
      \dimen0=#5\count2=\dimen0                    % desired width
      \dimen0=#2\count1=\dimen0                    % actual width
   \divide\count1 by 1000
   \divide\count2 by \count1                 % \count2 is now = scalefactor
  %\count3=11820\divide\count3 by \count2
  %\message{The vertices of #1 should have width \the\count3}
   \hbox{\vrule #4width0pt\vbox to 0pt{\vss\vskip#6%
      \special{picture #1 scaled \the\count2}\hrule width#5 height0pt\vss}}}}
%
%
%The following \figure macro builds on the \epsf macro and includes an EPS file.
% THE FOLLOWING LINE HAS TO OCCUR AT THE START OF THE MAIN FILE TO BE TEXED:
% \input epsf.def
%
% #1 is just a figure number to be used in the caption.
% #2 is the caption itself.
% #3 is the file name of the figure; this must an EPS file
%    (or PS with bounding box).
% #4 is for scaling, but disabled now (see below)
%
% Syntax example: \figure 3. The funny graph $G$ (Funny.graph.eps; 800)
%
\def\figure #1. #2 (#3; #4) {{%
   \def\bigskip{\par\ifdim\lastskip<\bigskipamount\removelastskip
                                              % eg. for abb after \endproof
      \vskip\bigskipamount\fi}% takes effect inside the def(!) of midinsert
   \midinsert\vskip\topfiguremargin
   \dimen0=\normalpictureheight
      \divide\dimen0 by 1000
      \multiply\dimen0 by#4        % \dimen0 := intended pictureheight
   \centerline{\epsfbox{#3.eps}}%                  good placing - use this!
   \vskip\bottomfiguremargin                     % beginning caption
      \ninepoint
      \parindent=.1\hsize\narrower\narrower
      \setbox0\hbox{#2}
      \ifdim\wd0 < .6\hsize
           \centerline{F{\sc IGURE} #1.\hskip1em#2}
       \else
           \plainitem{F{\sc IGURE} #1. }#2\par
       \fi
  \endinsert}}
%
% The following is for notes for talks to be hidden on transparancies.
% Say \hide{bla} to make "bla" disappear (but it will take up space);
% \showhidden will show all hidded text in grey (eg for paper notes).
%

%
%                               REFERENCES
%
\def\Abh#1 {{\sl Abh.\g Math.\g Sem.\g Univ.\g Hamburg\penalty100\ \bf#1\ }}
\def\AMASH#1 {{\sl Acta Math.\g Acad.\g Sci.\g Hung.\penalty100\ \bf#1\ }}
\def\Advances#1 {{\sl Adv.\g Math.\penalty100\ \bf#1\ }}
\def\Annals#1 {{\sl Ann.\g Math.\penalty100\ \bf#1\ }}
\def\AnnComb#1 {{\sl Ann.\g Comb.\penalty100\ \bf#1\ }}
\def\AMM#1 {{\sl Amer.\g Math.\g Monthly\penalty100\ \bf#1\ }}
\def\Archiv#1 {{\sl Arch.}\g {\sl Math.\penalty100\ \bf#1\ }}
\def\ArsComb#1 {{\sl Ars Comb.\penalty100\ \bf#1\ }}% FROMERLY \AC
\def\CJM#1 {{\sl Can.\g J.\th Math.\penalty100\ \bf#1\ }}
\def\Comb#1 {{\sl Com\-bi\-na\-to\-ri\-ca\penalty100\ \bf#1\ }}
\def\CPC#1 {{\sl Comb.\g Probab.\g Comput.\penalty100\ \bf#1\ }}
\def\Crelle#1 {{\sl J.}\th {\sl Reine Angew.}\g
    {\sl Math.\penalty100\ \bf#1\ }}
\def\DM#1 {{\sl Discrete Math.\penalty100\ \bf#1\ }}
\def\DAM#1 {{\sl Discrete Appl.\g Math.\penalty100\ \bf#1\ }}
\def\EJC#1 {{\sl Eur.}\g{\sl J.}\g{\sl Comb.\penalty100\ \bf#1\ }}
\def\EJ#1 {{\sl Electronic.}\g{\sl J.}\g{\sl Comb.\penalty100\ \bf#1\ }}
\def\GC#1 {{\sl Graphs Comb.\penalty100\ \bf#1\ }}
\def\IJ#1 {{\sl Isr.\g J.\th Math.\penalty100\ \bf#1\ }}
\def\Inv#1 {{\sl In\-vent.\g math.\penalty100\ \bf#1\ }}
\def\JAlg#1 {{\sl J.}\th {\sl Algorithms\penalty100\ \bf#1\ }}
\def\JCTA#1 {{\sl J.}\th {\sl Comb.}\g {\sl Theory~A\penalty100\ \bf#1\ }}
\def\JCTB#1 {{\sl J.}\th {\sl Comb.}\g {\sl Theory~B\penalty100\ \bf#1\ }}
\def\JGT#1 {{\sl J.}\th {\sl Graph Theory\penalty100\ \bf#1\ }}
\def\BLMS#1 {{\sl Bull.\g Lond.\g Math.\g Soc.\penalty100\ \bf#1\ }}
\def\JLMS#1 {{\sl J.\g Lond.\g Math.\g Soc.\penalty100\ \bf#1\ }}
\def\PLMS#1 {{\sl Proc.\g Lond.\g Math.\g Soc.\penalty100\ \bf#1\ }}
\def\Order#1 {{\sl Order\ \bf#1\ }}
\def\Random#1 {{\sl Random Struct.\g Alg.\penalty100\ \bf#1\ }}
\def\MA#1 {{\sl Math.}\g {\sl Ann.\penalty100\ \bf#1\ }}
\def\MN#1 {{\sl Math.}\g {\sl Nachr.\penalty100\ \bf#1\ }}
\def\MPCPS#1 {{\sl Math.\g Proc.\g Camb.\g Phil.\g Soc.\penalty100\ \bf#1\ }}
\def\MS#1 {{\sl Math.}\g {\sl Scand.\penalty100\ \bf#1\ }}
\def\MZ#1 {{\sl Math.}\g {\sl Zeit.\penalty100\ \bf#1\ }}
\def\BAMS#1 {{\sl Bull.\th Amer.\g Math.\g Soc.\penalty100\ \bf#1\ }}
\def\JAMS#1 {{\sl J.\th Amer.\g Math.\g Soc.\penalty100\ \bf#1\ }}
\def\MAMS#1 {{\sl Mem.\g Amer.\g Math.\g Soc.\penalty100\ \bf#1\ }}
\def\PAMS#1 {{\sl Proc.\g Amer.\g Math.\g Soc.\penalty100\ \bf#1\ }}
\def\SIAM#1 {{\sl SIAM J.\g Discrete Math.\penalty100\ \bf#1\ }}
\def\SLNM#1 {{\sl Springer Lecture Notes in Mathematics\penalty100\ \bf#1\ }}
\def\TAMS#1 {{\sl Trans.\g Amer.\g Math.\g Soc.\penalty100\ \bf#1\ }}
\def\TCSA#1 {{\sl Theor.\g Comput.\g Sci.~A\penalty100\ \bf#1\ }}
%
%
%
%
%
%                     ONLY FOR TeXtures without mssymb.tex:
%
%\def\N{{\rm \rlap I{\kern.18em}N}}   
%  \catcode`@=11        % This allows the use of `@' in the next line (p.344)
%\def\not#1{\mathrel{\mathpalette\c@ncel#1}} % use (only) with Imagewriter.
%  \catcode`@=12        % This restores the `inhibiting' catcode of `@'.
%\def\subsetneqq{\mathrel %{\mathchoice
%   {\vcenter{
%      \hbox{\lower6pt\hbox{$\scriptstyle\subset$}}
%      \hbox{\raise3pt\hbox{$\flatneq$}}}} }
%   \def\flatneq{\rlap {$\scriptstyle =$} {\kern1.5pt} {\scriptscriptstyle /}}
%\def\square{\Square53}            % See def. of \Square below.
%\def\Square#1#2{{\vbox{\hrule height.#2pt
%       \hbox{\vrule width.#2pt height#1pt \kern#1pt
%          \vrule width.#2pt}
%       \hrule height.#2pt}}}
%\def\nexists{\hbox{\rm\raise1.2pt\rlap/$\exists$}}
%
%
%
%
%
%\def\language=#1{}% For use with Textures versions < 1.3
%
%
%                  MODIFICATIONS TO PLAIN TeX
%
%\catcode`[=\active \catcode`]=\active
%  \def[{\thinspace\lbrack\thinspace}
%  \def]{\thinspace\rbrack}
%\def\{{\lbrace\thinspace}
%\def\}{\thinspace\rbrace}
%
\bigskipamount=1\baselineskip plus.3\baselineskip minus.3\baselineskip
\medskipamount=\bigskipamount\divide\medskipamount by 2
\smallskipamount=\medskipamount\divide\smallskipamount by 2 % (p.349)
\medmuskip = 3mu plus 2mu minus 1mu
\thickmuskip = 6mu plus 4mu minus 2mu % (for spacing in formulae; pp.168/170)
\def\smallbreak{\par \ifdim\lastskip<\smallskipamount
   \removelastskip \penalty-100 \smallskip \fi}
\def\medbreak{\par \ifdim\lastskip<\medskipamount
   \removelastskip \penalty-250 \medskip \fi}
\def\bigbreak{\par \ifdim\lastskip<\bigskipamount
   \removelastskip \penalty-500 \bigskip \fi}
\catcode`@=11        % This allows the use of `@' in the next line (p.344)
  \def\raggedbottom{\topskip10pt plus20pt \r@ggedbottomtrue} % The amount of
%                      permitted raggedness is controlled by the `plus' item.
\catcode`@=12        % This restores the `inhibiting' catcode of `@'.
\def\ge{\geqslant}% \geq remains available for the default version of \ge.
\def\le{\leqslant}% \leq remains available for the default version of \le.
\let\elt=\in
\def\in{\mathrel{\mathchoice
   {\raise .7pt \hbox{$\scriptstyle\elt$}}
   {\raise .7pt \hbox{$\scriptstyle\elt$}}
   {\raise .5pt \hbox{$\hskip .5pt\scriptscriptstyle\elt\hskip .5pt$}}
   {\raise.35pt \hbox{$\scriptscriptstyle\elt$}} }}
\let\hasaselt=\owns
\def\owns{\mathrel{\mathchoice
   {\raise .7pt \hbox{$\scriptstyle\hasaselt$}}
   {\raise .7pt \hbox{$\scriptstyle\hasaselt$}}
   {\raise .5pt \hbox{$\hskip .5pt\scriptscriptstyle\hasaselt\hskip .5pt$}}
   {\raise.35pt \hbox{$\scriptscriptstyle\hasaselt$}} }}
\let\exis=\exists
   \def\exists{\exis\>}
\let\nexis=\nexists
   \def\nexists{\nexis\>}
                            % To be phased out
\let\foral=\forall
   \def\forall{\foral\>}
\let\Rightarro=\Rightarrow
   \def\Rightarrow{\>\Rightarro\>}
\let\mi=\min
   \def\min{\mi\>}
\let\ma=\max
   \def\max{\ma\>}
\let\su=\sup
   \def\sup{\su\>}
\let\inff=\inf
   \def\inf{\inff\>}
\mathchardef\to="2221   % = \rightarrow, but of `binop' type (p.154)
\def\proclaim #1.#2 #3\par{\bigbreak
   \noindent{\bf#1.}#2\enspace{\sl#3}\par\bigbreak}
% The second argument above is optional and intended for references. Note that
% it is terminated by a space, which must therefore be present unsuppressed in
% input. Syntax examples: "\proclaim Thm {1.3}. This is the theorem.\par" or
% "\proclaim Thm \Euler.\five{} This is Euler's theorem.\par", where \Euler
%  expands to {1.3} and \five to \th [5], say. (Note that \proclaim Thm 1.3.
% would treat the 3. as argument #2 (incorrectly) and fail to set it in bold.)
% Or directly: "\proclaim Theorem \xxxVTop.~[\the\ref] blabla" (note the ~).
%
\newskip\sectionheadlineskipamount
\sectionheadlineskipamount=8pt plus 2pt minus 1pt
\def\beginsection #1\par{\Bigbreak\centerline{\bold #1}
        \penalty10000\vskip\sectionheadlineskipamount\noindent}
\let\ffootnote=\footnote
\def\footnote#1#2{\ffootnote{#1}{\eightpoint#2\vskip-12pt}}
%                  (The \vskip inserts an implicit \par, which has two
%                    effects: first, the desired effect of wrapping
%                    up the last paragraph of the footnote giving it the
%                    correct linespacing (that of \eightpoint), secondly
%                    the undesired effect of starting a new paragraph with
%                    a strut. To counteract the arising blank vertical
%                    space, the skip is chosen negative.)
%
\newcount\footnoteno% currently reset to 0 in chapter macro
\def\Footnote#1{{\footnote{${}^{\the\footnoteno}$}{#1}%
   \global\advance\footnoteno by 1}}
%
% Note that the following redefinition of \item finishes with a smallbreak.
% Since two smallbreaks in sequence result in only a single smallbreak, this
% gives a smallskip at the beginning, between any two items, and at the end.
% if additional space is desired before and after a series of items, say
% \medskip before and \par\smallskip (not \smallbreak) after the series of
% items.
%
\def\item#1#2\par{\parindent=10mm\smallbreak\hang\indent
                  \llap{{\rm #1}\enspace}\ignorespaces#2\par\smallbreak
                  \parindent=7mm}
\def\itemitem#1#2\par{\parindent=10mm\smallbreak
                  \indent\hangindent2\parindent\indent
                  \llap{{\rm #1}\enspace}\ignorespaces#2\par\smallbreak
                  \parindent=7mm}
%
%
%                     FORMAT PARAMETERS
%
\pretolerance=0 %This prevents line breaks in maths formulas - I don't know why.
\tolerance=2000
%\fontdimen2\tenrm=3.8pt
%\fontdimen2\tensl=3.8pt
%\fontdimen2\tenit=4pt
\baselineskip=13pt                 %(DEFAULT IS 12pt)
\vsize=200mm                   % (used to be 240truemm; changed 06/98)
\hsize=120mm                   % (used to be 140truemm; changed 06/98)
\hoffset=9mm                   % (used to be 9truemm; changed 06/98)
\parindent=7mm
\relpenalty=2000 \binoppenalty=5000  % DISCOURAGES BREAKS IN FORMULAS
\hyphenpenalty=100
\abovedisplayskip=12pt plus3pt minus4pt
\belowdisplayskip=12pt plus3pt minus4pt    % (p.348)
%
%   Unfortunately, TeX is unable to adjust the skip following a display to the
%   length of the line **below** the display, in the way in which it chooses
%  between \abovedisplayskip and \abovedisplayshortskip depending on the length
%   of the line above it. Thus, such adjustment has to be done by hand: setting
%
\belowdisplayshortskip=9pt plus3pt minus3pt
    % (formerly 12pt, like \belowdisplayskip)
%
%   reduces the effect of the standard 'short' version, while
%
%                  \def\noskip{\vskip-\lastskip\noindent}
%
%   (which is contained in macros.tex) removes any belowdisplay skip
%   (long or short) altogether.
%   One may want to say \smallskip\noindent just after it.
%
 \hyphenation{ac-cess-ible ana-log-ous ana-log-ous-ly ana-lyze ana-lyse
ana-ly-sis answer answers aver-age axio-mat-ic bundle bundles Buch-ge-sell-schaft col-our
col-ours col-oured col-our-ing col-our-ings con-struct-ible con-struct-ive
con-struct-ive-ly co-rol-lary Co-rol-lary des-cend des-cend-ing Deut-sche elem-ent elem-ents
end-li-cher de-fi-ni-tion de-fi-ni-tions De-fi-ni-tion equi-val-ent
equi-val-ence Euler-ian exist-ence every Gra-phen Hamil-ton-ian homeo-mor-phic
homeo-mor-phism homeo-mor-phisms hy-po-thesis hy-po-theses in-ac-cess-ible
ir-regu-lar ir-regu-lar-ity method methods modi-fi-ca-tion mono-chro-matic ori-ent ori-ent-ed par-ticu-lar
popu-la-tion popu-la-tions Popu-la-tion Popu-la-tions pro-po-si-tion pro-po-si-tions Pro-po-si-tion regu-lar regu-lar-ity regu-lar-ly rele-vant
sig-ni-fi-cant sig-ni-fi-cant-ly sig-ni-fi-cance speci-fi-ca-tion Speci-fi-ca-tion situ-ation situ-ations Situ-ation Situ-ations to-po-lo-gical to-po-lo-gical-ly
ubi-qui-tous ubi-quity un-at-tached un-end-li-cher using Using Wis-sen-schaft-li-che}
\input btxmac% for bibtex
\disablems
\noremarks
\nocomments

\def\secIdea{1}
\def\secInformalExamples{2}
\def\secFormalSetup{3}
\def\secThms{4}
\def\secFormalExamples{5}

\def\xxxToTfixedS{1}
\def\xxxToTallSk{2}
\def\xxxTTD{3}

\def\longsubsection #1 #2\par{\vskip 2\medskipamount plus \smallskipamount minus \smallskipamount\varitem{\bold #1} {\bold  #2}
        \penalty10000\smallskip\noindent}

\hbox{}\vskip2cm

\centerline{\bigbold Tangles in the Social Sciences}
\medskip\centerline{\it A new mathematical model to identify types and predict behaviour\rm%
   \footnote*{This is a draft paper written, by a mathematician, for a readership in the social sciences (including economics) interested in the mathematical underpinnings of quantitative methods in their fields. Collaboration with such readers is actively sought and will be most welcome.\endgraf
   This is part of a broader project, {\sl Tangles in the empirical sciences}~\cite{TanglesEmpirical}, which includes applications in the natural sciences too. It adds to the account given here a description of the same theory of tangle applications in the language of clustering in large data sets, which I have avoided here in favour of a less technical descriptive style. The complementary account given in~\cite{TanglesEmpirical} can also help with forming an intuitive geometric picture of the tangles discussed here.\looseness=-1}}
\vskip 5mm
\centerline{Reinhard Diestel}
\disablems
\noremarks
\nocomments

\footnoteno=1

\def\lowfwd #1#2#3{{\mathop{\kern0pt #1}\limits^{\kern#2pt\raise.#3ex \vbox to 0pt{\hbox{$\scriptscriptstyle\rightarrow$}\vss}}}}
\def\fwd #1#2{{\lowfwd{#1}{#2}{15}}}
\def\lowbkwd #1#2#3{{\mathop{\kern0pt #1}\limits^{\kern#2pt\raise.#3ex
\vbox to 0pt{\hbox{$\scriptscriptstyle\leftarrow$}\vss}}}}
\def\vp{\lowfwd p{1.5}1}
\def\pv{\lowbkwd p{1.5}1}
\def\vq{\lowfwd q{1.5}1}
\def\qv{\lowbkwd q{1.5}1}
\def\vr{\lowfwd r{1.5}1}
\def\rv{\lowbkwd r{1.5}1}
\def\vS{\vec S}
\def\vSi{{\hskip-1pt{\fwd {S_i}3}\hskip-1pt}}
\def\vSk{{\hskip-1pt{\fwd {S_k}3}\hskip-1pt}}
\def\vSell{{\hskip-1pt{\fwd {S_\ell}3}\hskip-1pt}}
\def\vSdash{{\mathop{\kern0pt S\lower-1pt\hbox{${}% logically \v(e')
     \scriptstyle'$}}\limits^{\kern2pt\raise.1ex
     \vbox to 0pt{\hbox{$\scriptscriptstyle\rightarrow$}\vss}}}}
\def\vs{\lowfwd s{1.5}1}
\def\sv{\lowbkwd s{1.5}1}
\def\vsdash{{\mathop{\kern0pt s\lower.5pt\hbox{${}% logically \v(e')
     \scriptstyle'$}}\limits^{\kern0pt\raise.02ex
     \vbox to 0pt{\hbox{$\scriptscriptstyle\rightarrow$}\vss}}}}
\def\svdash{{\mathop{\kern0pt s\lower.5pt\hbox{${}% logically \v(e')
     \scriptstyle'$}}\limits^{\kern0pt\raise.02ex
     \vbox to 0pt{\hbox{$\scriptscriptstyle\leftarrow$}\vss}}}}

\def\svone{\lowbkwd {s_1}01}

\def\svfive{\lowbkwd {s_5}01}
\def\vsi{\lowfwd {s_i}11}
\def\svi{\lowbkwd {s_i}11}
\def\vt{\lowfwd t{1.5}1}

\def\vT{\vec T}
\def\vfwd{\lowfwd v{1.5}1}
\def\vbwd{{{\lowbkwd v{1.5}1}\hskip-1pt}}

\def\vV{\vec V}

\def\vv{\lowfwd v{1.5}1}
\def\vvback{\lowbkwd v{1.5}1}

\def\vvi{\lowfwd {v_i}11}
\def\vviback{\lowbkwd {v_i}11}

\def\vx{\lowfwd x{1.5}1}
\def\xv{\lowbkwd x{1.5}1}
\def\vX{\vec X}
\def\vY{\vec Y}
\def\vy{\lowfwd y{1.5}1}
\def\yv{\lowbkwd y{1.5}1}
\def\vz{\lowfwd z{1.5}1}

\def\Abar{{\overline A}}
\def\Bbar{{\overline B}}
\def\Cbar{{\overline C}}

\def\F{{\cal F}}
\def\TT{{\cal T}}

\bigskip\bigskip
{\narrower\narrower\noindent
   Traditional clustering identifies groups of objects that share certain qualities. Tangles do the converse: they identify groups of qualities that often occur together. They can thereby identify and discover {\it types\/}: of behaviour, views, abilities, dispositions.\endgraf
   The mathematical theory of tangles has its origins in the connectivity theory of graphs, which it has transformed over the past 30 years. It has recently been axiomatized in a way that makes its two deepest results applicable to a much wider range of contexts.\looseness=-1\endgraf
   This expository paper indicates some contexts where this difference of approach is particularly striking. But these are merely examples of such contexts: in principle, it can apply to much of the quantitative social sciences.\endgraf
   Our aim here is twofold: to indicate just enough of the \hbox{theory} of tangles to show how this can work in the various different contexts, and to give plenty of different examples illustrating this.
   \par}

\vfill\eject

\beginsection Contents

\medbreak
\def\n#1\par {\noindent #1\smallskip}
\def\nn#1\par {#1\par}

\n Introduction\medbreak

\n 1. The idea behind tangles

\nn 1.1 \ Features that often occur together

\nn 1.2 \ Consistency of features

\nn 1.3 \ From consistency to tangless

\nn 1.4 \ The predictive power of tangles

\nn 1.5 \ Witnessing sets and functions\medbreak

\n 2. Examples from different contexts

\nn 2.1 \ Sociology: discovering mindsets, social groups, and character traits

\nn 2.2 \ Psychology: understanding the unfamiliar

\nn 2.3 \ Politics and society: appointing representative bodies

\nn 2.4 \ Education: combining teaching techniques into methods

\nn 2.5 \ Linguistics and analytic philosophy: how to determine meaning

\nn 2.6 \ Economics: identifying customer and product types\medbreak

\n 3. The formal setup for tangles

\nn 3.1 \ Tangles of set partitions

\nn 3.2 \ The evolution of tangles

\nn 3.3 \ Stars, universes, and submodularity

\nn 3.4 \ Hierarchies of features and order

\nn 3.5 \ Duality of set separations\medbreak

\n 4. Tangle theorems and algorithms

\nn 4.1 \ Tangle-distinguishing feature sets

\nn 4.2 \ Tangle-precluding feature sets

\nn 4.3 \ Algorithms

\nn 4.4 \ Tangle-based clustering

\nn 4.5 \ Predictions\medbreak

\n 5. Applying tangles: back to the examples

\nn 5.1 \ Sociology: from mindsets to matchmaking

\nn 5.2 \ Psychology: diagnostics, new syndromes, and the use of duality

\nn 5.3 \ Political science: appointing representative bodies

\nn 5.4 \ Education: devising and assigning students to methods

\nn 5.5 \ Meaning: from philosophy to machine learning

\nn 5.6 \ Economics: customer and product types\bigbreak

\noindent
References

\beginsection Introduction

Suppose we run a survey $S$ of fifty political questions on a population~$P$ of a thousand people. If there exists a group of, say, a~hundred like-minded people among these, there will be a `typical' way of answering the questions in~$S$ in the way most of those people would. Quantitatively, there will exist a subset~$X$ of~$P$, not too small, and an assignment~$\tau$ of answers to all the questions in~$S$ such that, for most%
   \COMMENT{}
   questions $s\in S$, some 80\% (say) of the people in~$X$ agree with the answer to~$s$ given by~$\tau$. (Which 80\% of~$X$ these are will depend on the choice of~$s$.) We call this collection~$\tau$ of views~-- answers to the questions in~$S$~-- a {\it mindset\/}. Note that there may be more than one mindset for~$S$, or none.\looseness=-1

Traditionally, mindsets are found just intuitively: they are first guessed, and only then established by quantitative evidence from surveys designed to test them. For example, we might feel that there is a `socialist' way $\sigma$ of answering~$S$. To support this intuition, we might then check whether any sizable subset~$X\sub P$ as above exists for this particular $\tau :=\sigma$.

Tangles can do the converse: they will identify both $X$ and~$\tau$ without us having to guess them first. For example, tangle analysis of political polls in the UK in the years well before the Brexit referendum might have established the existence of a mindset we might now,%
   \vadjust{\penalty-500}
   with hindsight, call the `Labour-supporting non-socialist Brexiteer': a~mindset whose existence few would have guessed intuitively when Brexit was not yet on the agenda. And similarly in the US with the MAGA%
   \Footnote{Make America Great Again; Donald Trump's 2016 presidential campaign slogan.}
   mindset before 2016, or that of a `conservative Green' in the early 1970s. Tangles can identify previously unknown patterns of coherent views or behaviour.

There are two main mathematical theorems about tangles. One of them extracts, in the context of the above example, from the set $S$ of all questions and the way the people in $P$ answered them a small set $T$ of `critical' questions%
   \COMMENT{}
   that suffice to distinguish most of%
   \COMMENT{}
   the mindsets that exist in the popu\-la\-tion surveyed by~$S$.%
   \Footnote{The questions in~$T$ will either themselves be in~$S$ or be combinations of questions from~$S$.}
   For every two such mindsets $\sigma$ and~$\tau$ there will be a question in~$T$ on which $\sigma$ and~$\tau$ disagree, and which can thus be used to distinguish $\sigma$ from~$\tau$. As an immediate application, $T$~would make a good small questionnaire for a larger study if $S$ was a pilot study designed to test a large number of questions on a smaller population.

More fundamentally, if $P$ is large enough to be representative for some larger population, the tangle-distinguishing quality of~$T$ means that for any individual represented by~$P$ (but not necessarily in~$P$), one can predict from their answers to the questions in~$T$ how they would answer the much larger set of questions in~$S$. Such predictions will be reliable as soon as the individual's views are aligned with {\it some\/} mindset~-- and the theory of tangles can determine how likely this is to be the case.%
   \COMMENT{}

The gain in conditioning predictions for the unknown answers to questions in~$S$ on the known answers to the questions in this particular set~$T$, rather than an arbitrary small subset of~$S$, is that $T$ is so special: the answers to $T$ already identify which mindset on~$S$, if any, an individual belongs to,%
   \COMMENT{}
   which enables us to base our predictions for this individual's views (on all of~$S$) on this mindset.

The other fundamental theorem about tangles is about cases when no mindset exists for the questions in~$S$. The algorithm for this theorem is designed to%
   \COMMENT{}
   find in such cases a small subset $T$ of~$S$%
   \COMMENT{}
   that has the following special property: for every possible consistent%
   \Footnote{Consistency is an easily checkable property of any mindset. It will be formalized later.}
   way~$\tau$ of answering~$T$ there will be a small set~$S_\tau\sub T$ of questions, perhaps three or four, whose answers by~$\tau$ are not shared by enough people for $\tau$ to be a mindset.%
   \COMMENT{}%
   \COMMENT{}%
   \COMMENT{}
   Hence $S_\tau$~certifies that this~$\tau$, although consistent, cannot be a mindset (on~$T$). Since mindsets are consistent, and hence would be such sets~$\tau$, this means that there are no mindsets on~$T$. But then there are no mindsets on~$S$ either, because they would include mindsets on~$T$. So our set~$T$, together with all those small sets~$S_\tau$, certifies that no mindset on~$S$ exists.%
   \COMMENT{}

This second theorem, therefore, furnishes `negative' poll returns with a verifiable proof that mindsets not only were not found but do not in fact exist.
As with the first theorem, such special sets $T$ and~$S_\tau$ of questions cannot simply be found by trial and error, since there are too many subsets~$T$ of~$S$, too many consistent ways~$\tau$ of answering~$T$, too many subsets $S_\tau$ of every such~$\tau$, and too many people in~$P$ to test on each of those on~$S_\tau$.

In both theorems, the sets $T$ and~$S_\tau$ are highly valuable in both senses of the word: knowing them gives us a lot of valuable information that we would not otherwise have, and finding them requires some nontrivial mathematics. The theory of tangles provides this, and in many instances $T$ and the~$S_\tau$ can be computed easily.

This paper is organized as follows. Chapter~\secIdea\ describes the basic idea behind tangles. It is still written informally, but in a slightly more abstract language than used so far, so that other examples can be expressed in it too.

Such examples will be introduced in Chapter~\secInformalExamples. These are by no means a complete cross-section of where tangles can be used. My aim in selecting these examples was to choose potential applications where using tangles can make a particularly striking difference. And not too many: readers are encouraged to take the examples discussed here simply as templates for possible applications in their own field. The fields touched upon in Chapter~\secInformalExamples\ are:
   {\smallskip\advance\parindent by 7mm
   \plainitem{$\bullet$} Sociology
   \plainitem{$\bullet$} Psychology
   \plainitem{$\bullet$} Political science
   \plainitem{$\bullet$} Education
   \plainitem{$\bullet$} Linguistics and analytic philosophy
   \plainitem{$\bullet$} Economics
   \smallbreak}

Further examples can be found in~\cite{TanglesEmpirical}. They include applications in text analysis, medical diagnostic expert systems, DNA~sequencing, image processing, and generic clustering in large data sets. As here, those examples are just intended as pointers: pointers to possible uses of tangles in the empirical sciences that include, but go beyond, the social sciences discussed here.%
   \looseness=-1

In Chapter~\secFormalSetup\ we describe a more formal setup for tangles. It is not the most general formal setup,%
   \Footnote{This is presented concisely in~\cite{AbstractSepSys}, which includes further references for the mathemically inclined to where the main tangle theorems are proved.}
    but it is general enough to express formally the various types of tangle discussed in this paper, and probably most  or all the types of tangles with applications in the social sciences.

Chapter~\secThms\ then states, in the formal language introduced in Chapter~\secFormalSetup, the two fundamental tangle theorems: the {\it tree-of-tangles\/} theorem and the {\it tangle-tree duality\/} theorem hinted at earlier. It discusses tangle algorithms corresponding to these theorems, and their potential applications for predictions and traditional clustering. Formality is kept to the minimum needed to enable readers to check precisely where tangles might, in principle, play a role in their own fields.%
   \Footnote{Once this has been established in principle, my hope is that such readers get in touch to see how this can be implemented. We have developed some generic algorithms that compute the solutions offered by those theorems on a wide range of potential input data. However for best performance it will be best to fine-tune them to the concrete application intended.}\looseness=-1

Chapter~\secFormalExamples\ then revisits the examples from Chapter~\secInformalExamples\ in the light of the additional precision gained in Chapters \secFormalSetup\ and~\secThms. The aim now is to re-cast these examples, whose essence is already understood from their informal treatment in Chapter~\secInformalExamples, in a sufficiently formal way that they can be used as templates for similar applications. We shall also be able to describe some deeper aspects of these examples that could not be expressed earlier at the informal level.

\beginsection \secIdea. The idea behind tangles

Consider a collection $V\!$ of objects and a set $\vS$ of features%
   \Footnote{Logicians may prefer to say `predicates' instead of `features' here. That would be correct,\penalty-200\ but I~am trying to avoid any (false) impression of formal precision at this stage.}
   that each of the objects in~$V\!$ may have or fail to have. Given such a (potential) feature~$\vs\in\vS$, we denote its negation by~$\sv$. The pair~$\{\vs,\sv\}$ of the feature together with its negation is then denoted by~$s$, and the set of all these~$s$ is denoted by~$S$.

For example, if $V\!$ is a set of pieces of furniture, then $\vs$ might be the feature of being made of wood. Then $\sv$~would be the feature of being made of any other material, or a combination of materials, and $s$ could be thought of as the question of whether or not a given element of~$V\!$ is made of wood.

In our example from the introduction, $V\!$~would be the population~$P$ of people polled by our survey~$S$ (which, for simplicity, we assume to consist of yes/no questions). Then $\{\,\vs\mid s\in S\,\}$ might be the set of `yes' answers to the questions in~$S$, while $\sv$ would denote the `no' answer to the question~$s$.%
   \Footnote{The arrow notation has nothing to do with vectors. For the moment, it is just a convenient way of associating with a question~$s$ its two possible answers `yes' and~`no', or with a feature~$\vs$ its negation~$\sv\!$. In Chapter~\secFormalSetup\ the arrow notation will acquire an additional layer of meaning that does have to do with pointing, but not in any geometric sense such as that of vectors.}

\subsection \secIdea.1\ \ Features that often occur together

Tangles%
   \Footnote{The word was first introduced by Robertson and Seymour~\cite{GMX} in their ground-breaking work on graph minors. They use it for a region of a graph that hangs together in an intricate way. Intricate in that, while being close-knit in the sense of being difficult to separate, it still does not conform to the usual graph-theoretic notions of being a highly connected subgraph, minor or similar. The abstract notion of tangles that underlies our discussions in this paper was first introduced in~\cite{AbstractSepSys}.}
   are a way to formalize the notion that some features typically occur together. They offer a formal way of identifying such groups of `typical' features, each `type' giving rise to a separate tangle. One particular strength of tangles is that the sets of {\it objects\/} (elements of~$V$) whose features largely agree with a given type need not be clearly delineated from each other. As in most real-world examples, these sets may be `fuzzy': in order for certain features to form a type, or tangle, it is {\it not\/} necessary that there exist a group of objects in~$V\!$ that have precisely these features~-- indeed there need not even exist a single such object.

Let us return to the example where $V\!$~is a set of pieces of furniture. Our list $\vS$ of possible features  (including their negations) consists of qualities such as colour, material, the number of legs, intended function, and so on~-- perhaps a hundred or so potential features. The idea of tangles is that, even though $\vS$ may be quite large, its elements may combine into groups that correspond to just a few types of furniture as we know them: chairs, tables, beds and so on.

The important thing is that tangles can identify such types without any prior intuition: if we are told that a container~$V\!$ full of furniture is waiting for us at customs in the harbour, and all we have is a list of items~$v$ identified only by numbers together with, for each number, a list of which of our 100~features this item has, our computer~-- if it knows tangles~-- may be able to tell us that our delivery contains furniture of just a few types: types that we  (but not our computer) might identify as chairs, tables and beds, perhaps with the beds splitting into sofas and four-posters.

\subsection \secIdea.2\ \ Consistency of features

To illustrate just how our computer may be able to do this, let us briefly consider the inverse question: starting from a known type of furniture, such as chairs, how might this type be identifiable from the data if it was {\it not\/} known?

A~possible answer, which will lead straight to the concept of tangles, is as follows. Each individual piece of furniture in our unknown delivery, $v\in V\!$ say, has some of the features from our list~$\vS$ but not others. It thereby {\it specifies\/} the elements $s$ of~$S$: as $\vs$ if it has the feature~$\vs$, and as $\sv$ otherwise.%
   \Footnote{For mathematicians: note that this is well defined even when both $\vs,\sv$ are in~$\vS$, as we assume.}
   We say that every $v\in V\!$ defines a {\it specification\/} of~$S$, a choice for each $s\in S$ of either $\vs$ or~$\sv$ but not both. We shall denote this specification of~$S$ as
 $$v(S) := \{\,v(s)\mid s\in S\,\},$$
 where $v(s):=\vs$ if $v$ specifies~$s$ as~$\vs$ and $v(s) := \sv$ if $v$ specifies~$s$ as~$\sv$.

Conversely, does every specification of~$S$ come from some $v\in V\!$ in this way?\penalty-200\ Certainly not: there will be no object in our delivery that is both made entirely of wood and also made entirely of steel. Thus no $v\in V\!$ will specify both $r$ as $\vr$ rather than~$\rv$, and $s$ as $\vs$ rather than~$\sv$, when $\vr$ and $\vs$ stand for being made of wood or steel, respectively. In plain language: no specification of~$S$ that comes from a real piece of furniture can contain both $\vr$ and~$\vs$, because these features are inconsistent.

Let us turn this manifestation in~$V\!$ of logical inconsistencies within $\vS$ into a definition of `factual' inconsistency for specifications of~$S$ in terms of~$V\!$. Let us call a specification of~$S$ {\it consistent\/}%
   \COMMENT{}
   if it contains no inconsistent triple, where an {\it inconsistent triple\/} is a set of up to three%
   \Footnote{It might seem more natural to say `two' here, as in our wood/steel dichotomy above. Our definition of consistency is a little more stringent, because the mathematics behind tangles requires it. Note that, formally, the elements of an inconsistent `triple' need not be distinct; an `inconsistent pair' of {\it two\/} features $\vr,\vs$ not shared by any $v\in V\!$, for example, also counts as an `inconsistent triple', the triple $\{\vr,\vr,\vs\} = \{\vr,\vs\}$.}%
   \COMMENT{}
   features that are not found together in any $v\in V\!$. Specifications of~$S$ that come from some $v\in V\!$ are clearly consistent. But $S$ can have many consistent specifications that are not, as a whole, witnessed by any $v\in V\!$.%
   \Footnote{Here is a simple example. Suppose some of our furniture is made of wood, some of steel, some of wicker, and some of plastic. Denote these features as $\vp$, $\vq$, $\vr$, $\vs$, respectively, and assume that $S = \{p,q,r,s\}$. Then the specification $\sigma = \{\pv, \qv, \rv, \sv\}$ of~$S$ is consistent, because for any three of its elements there are some items in~$V\!$ that have none of the three corresponding features: those that have the fourth. But no item fails to have all four of these features. So the consistent specification~$\sigma$ of~$S$ does not come from any one $v\in V\!$.}%
   \COMMENT{}

Tangles will be specifications of~$S$ with certain properties that make them `typical' for~$V\!$. Consistency will be a minimum requirement for this. But since any specification of~$S$ that comes from just a single $v\in V\!$ is already consistent, tangles will have to satisfy more than consistency to qualify as `typical' for~$V\!$.

\subsection \secIdea.3\ \ From consistency to tangles

It is one of the fortes of tangles that they allow considerable freedom in the definition of what makes a specification of~$S$ `typical' for~$V\!$~-- freedom that can be used to tailor tangles precisely to the intended application.%
   \vadjust{\penalty-200}
   We shall describe this formally in Chapter~\secFormalSetup. But we are already in a position  to mention one of the most common ways of defining `typical', which is just a strengthening of consistency.

To get a prior feel for our (forthcoming) formal definition of `typical', consider the specification of~$S$ in our furniture example that is determined by an `ideal chair' plucked straight from the Platonic heaven: let us specify each $s\in S$ as $\vs$ if a typical chair has the feature~$\vs$, and as $\sv$ if not.%
   \Footnote{Let us ignore for the moment the possibility that the question~$s$ may not have a clear answer for chairs, as would be the case, say, for questions of colour rather than function. This is an issue we shall have to deal with, but which tangles can indeed deal with easily.}\penalty-200
   \COMMENT{}
   This can be done independently of our delivery~$V\!$, just from our intuitive notion of what chairs are. But if our delivery has a sizable portion of chairs in it, then this phantom specification of~$S$ that describes our `ideal chair' has something to do with~$V\!$ after all.

Indeed, for every triple $\vr,\vs,\vt$ of features of our `ideal chair' there will be a few elements of~$V\!$, at least~$n$ say, that share these three features. For example, if $\vr,\vs,\vt$ stand for having four legs, a flat central surface, and a near-vertical surface, respectively, there will be~-- among the many chairs in~$V\!$ which we assume to exist~-- a few that have four legs and a flat seating surface and a nearly vertical back.

By contrast, if we pick twenty rather than three features of our ideal chair there may be no $v\in V\!$ that has all of those, even though there are plenty of chairs in~$V\!$. But for {\it every\/} choice of three features there will be several~-- though which these are will depend on which three features of our ideal chair we have in mind.

Simple though it may seem, it turns out that for most furniture deliveries and reasonable lists~$S$ of potential features this formal criterion for `typical' distinguishes those specifications of~$S$ that describe genuine types of furniture from most of%
   \COMMENT{}
   its other specifications.%
   \Footnote{\dots of which there are many: if $S$ has 100 elements, there are $2^{100}$ specifications of~$S$.}
   But in identifying such specifications as `types' we made no appeal to our intuition, or to the meaning of their features.%
   \Footnote{This is not to say that the use of tangles is free of all preconceptions, biases etc.: the choice of~$S$, for example, is as loaded or neutral as is would be in any other study that starts with a survey. The statement above is meant relative to the given~$S$ once chosen. In Chapter~\secFormalSetup\ we shall discuss how the deliberate use of preconceptions, e.g.\ by declaring some questions in~$S$ as more fundamental than others, can help to improve tangles based on such preconceptions. We shall also see how to do the opposite: how to find tangles that arise naturally from the raw data of~$S$ and~$V\!$, without any further interference from ourselves.}

So let us make this property of specifications of~$S$ that describe `ideal' chairs, tables or beds into our formal, if still ad-hoc, definition of `typical': let us call a specification of~$S$ {\it typical\/} for~$V\!$ if every subset of at most three of its elements is shared by at least $n$ elements of~$V\!$, where $n$ is now a fixed parameter on which our notion of `typical' depends, and which we are free to choose.%
   \COMMENT{}

Crucially, this definition of `typical' is purely intrinsic: it depends on~$V\!$, but it makes no reference to what a typical specification of~$S$ `is typical of'. Specifications of ideal chairs, tables or beds are all typical in this sense: they all satisfy the same one definition.

Equally crucially, a specification of~$S$ can be typical for~$V\!$ even if $V\!$ has no element that has {\it all\/} its features at once.%
   \COMMENT{}
   Thus, we have a valid and meaningful formal definition of an `ideal something' even when such a thing does not exist in the real world, let alone in~$V\!$.

Relative to the definition of `typical' we can now define tangles informally:

\medskip\indent\indent\indent {\it A~\underbar{tangle} of~$S$ is any typical specification of~$S$.\/}%
   \COMMENT{}
   \medbreak
   \noindent
Since our ad-hoc definition of `typical' is phrased in terms of small subsets of~$\vS$, sets of size at most~3 (of which there are not so many), we can compute tangles without having to guess them first.

In particular, we can compute tangles of~$S$ even when $V\!$ is `known' only in the mechanical sense of data being available (but not necessarily understood), and $S$ is a set of potential features that are known, or assumed, to be rele\-vant but whose relationships to each other are unknown. Tangles therefore enable us to {\it find\/} even previously unknown `types' in the data to be analysed: combinations of features that occur together significantly more often than others.

The two fundamental tangle theorems mentioned in the introduction then allow us to either structure these types by refinement and determine a small set of critical features that suffice to distinguish them, or else produce verifiable evidence that the data is too unstructured to contain significant types at all.

\subsection \secIdea.4\ \ The predictive power of tangles%
   \COMMENT{}

The general idea of trying to predict a person's likely behaviour in a future situ\-ation from observations of his or her actual behaviour in some past \hbox{situations} is as old as humanity: as we learn `how someone ticks', we are better able to make such predictions.%
   \Footnote{We are talking about `predictions' here at a quantitative level of, at best, weather forecasts; not about absolute predictions whose failure, even once, invalidates the entire theory.\looseness=-1}

If we have the chance to choose, or even design, those earlier situations, e.g.\ by selecting them carefully from a collection of past situations in which our individual's behaviour was observed, or  by devising a test consisting of hypothetical situations the response to which our individual is willing to share with us, we have a chance of improving our predictions by choosing particularly rele\-vant such past or hypothetical scenarios. Tangles can help to identify these.

Let us think of $V\!$ as a set of people in whose likely actions we are interested, and of $\vS$ as the set of potential such actions. To keep things simple, let us consider the example from the introduction, where $\vS$ is a set of potential views of the people in~$V\!$,%
   \Footnote{Thus, `holding the view$\,\vs$' would count as an `action' in our new context; if desired, think of it as the action of answering the question $s$ as~$\vs$.}
   i.e., we wish to predict how a person would answer a question~$s$ from~$S$. Let us assume that, as the basis for our prediction, we are allowed to quiz that person on some small set~$T$ of questions,%
   \COMMENT{}
   and base our prediction for their answer to~$s$ on their known answers to the questions in~$T$. Our aim is to choose $T$ so as to make these predictions particularly good.

As indicated in the introduction, the tree-of-tangles theorem will likely produce a particularly valuable%
   \COMMENT{}
   and much smaller set~$T$ of questions%
   \COMMENT{}
   whose answers entail more predictive power for~$S$%
  \COMMENT{}
   than an arbitrary subset of~$S$ of that size would. This is because $T$ consists of just enough questions to distinguish all the tangles of~$S$, and tangles represent the typical ways to answer~$S$. The answers of an individual to the questions in~$T$ thus determine exactly one such type: there is only one typical way of answering all the questions in~$S$ that includes these particular answers to~$T$. This typical way of answering~$S$, a~tangle of~$S$, is a particularly good prediction for the answers of our individual of whom only his or her answers to~$T$ are known, because he or she is more likely to answer in some, and hence this,%
   \COMMENT{}
   typical way than in an arbitrary other way.%
   \COMMENT{}%
   \looseness=-1

It may seem that there is a practical problem with this approach in that we are trying to base our predictions for the answers to questions in~$S$ to be expected from the people in~$V\!$ on being able to compute all the tangles of~$S$ first, which may be possible only if we already know what we are trying to predict: the answers of all the people in~$V\!$ to all the questions in~$S$. However, there is no such problem.%
   \COMMENT{}
    In any real-world application we would indeed compute these tangles based on full knowledge of how the elements of $V\!$ specify~$S$. But we would be interested in predictions about individuals outside the set~$V\!$. For this to be possible, $V\!$~would have to be representative for that larger population~-- which, however, is a standard assumption one always has to make (and justify).

A~classical application would be that $S$ is a pilot study run on a small subset~$V\!$ of a larger population~$P$, and $T$ is a study with fewer but particularly rele\-vant questions selected from~$S$, to be run on~$P$. Then the answers to $T$ of an individual in $P$ can justify predictions on how this individual would answer the rest of~$S$, and $T$ is a particularly well-chosen subset of~$S$ for this purpose.

\subsection \secIdea.5\ \ Witnessing sets and functions

In the example discussed in the introduction, where $S$ is a political survey and tangles are mindsets~-- collections of views that are often held together~-- we quantified the meaning of `often held together'. We did so in terms of a subset $X$ of all the people in~$V=P$ that should not be too small (but comprise at least, say, 1/10 of the population~$P$) and whose views are typified by the given tangle~$\tau$ in that, for every question $s\in S$, some 80\% of~$X$ answer~$s$ the way $\tau$ does.%
   \Footnote{Recall that~$\tau$, being a tangle of~$S$ and thus a specification of~$S$, `answers' every question $s\in S$ by choosing either $\vs$ or~$\sv$.}
   In our furniture example, the tangle of being a chair would be witnessed in this way by the set $X$ of chairs in~$V$: every feature of our `ideal chair'~$\tau$ will be shared by some 80\% of all the chairs in~$V\!$, though not all by the same~80\%.

Formally, we say that a set $X\sub V\!$ {\it witnesses\/} a specification $\tau$ of~$S$ if, for every $s\in S$, there are more $v$ in~$X$ that specify $s$ as~$\tau$ does than there are $v\in X$ that specify $s$ in the opposite way. If these majorities are greater than~2/3, then $\tau$ is likely to be a tangle for most definitions of `typical' (see Section~\secIdea.3).%
   \COMMENT{}

More generally, a `weight' function $w\colon V\!\to\N$%
   \COMMENT{}
   {\it witnesses\/}~$\tau$ if, for every $s\in S$, the collective weight of the $v\in V\!$ that specify $s$ as $\tau$ does exceeds the collective weight of the $v\in V\!$ that specify $s$ in the opposite way.%
   \Footnote{Formally, $w\colon V\!\to\N$ {\it witnesses\/}~$\tau$ if $\sum \{\,w(v)\mid v(s)=\tau(s)\,\} > \sum \{\,w(v)\mid v(s)\ne \tau(s)\,\}$  for every $s\in S$.}

Much of the attraction and usefulness of tangles stems from the fact that, in practice, most of them have such witnessing sets or functions~\cite{Deciders_k-conn}.%
   \COMMENT{}
   But it is important to note that the definition of a tangle, even our preliminary definition from Section~\secIdea.3, does not require that such sets or functions exist. It relies only on notions of consistency and%
   \Footnote{Our example definition of `typical' in Section~\secIdea.3 implied consistency. But there are other useful notions of type that do not, in which case we simply add consistency to the requirements made of a tangle.}
   type, which are defined by banning small subsets of $\vS$ deemed `inconsistent' or `atypical' from occuring together in a tangle.\looseness=-1

In some contexts, tangles of~$S$ can even be defined without any reference to~$V\!$ at all. In our furniture example we could have defined the consistency of a set of features, or predicates, about the elements of~$V\!$ in purely logical or linguistic terms that make no appeal to~$V\!$. If $\vr$ stands for `made entirely of wood' and $\vs$ stands for `made entirely of steel',\penalty-200\ then the set~$\{\vr,\vs\}$ is inconsistent. But we have two ways of justifying this. The reason we chose to give was that no object in~$V\!$ is made entirely of wood and also made entirely of steel. But we might also have said that these two predicates are logically inconsistent~-- which implies that there is no such object in~$V\!$ but which can be established without examining~$V\!$.

The way consistency and type are defined formally~\cite{AbstractSepSys} as part of the notion of abstract tangles is something half-way between these two options: it makes no reference to~$V\!$ but refers only to some axiomatic properties of~$\vS$ which reflect our notion that $\vS$ is a set of `features'.%
   \COMMENT{}
   In this way it also avoids any appeal to logic or meaning.

For the rest of this paper the only important thing to note about witnessing sets or functions is that while many tangles have them, tangles can be identified, distinguished by the first tangle theorem, or ruled out by the second tangle theorem without any reference to such sets or functions. {\it Pars pro toto:}\ the mindset of being\penalty-200\ socialist can be identified without having to find any actual socialists, let alone delineating these as a social group against others.

\beginsection \secInformalExamples. Examples from different contexts

The aim of this chapter is to illustrate the range of potential uses of tangles by a handful of further examples from different contexts, still at the level of informality adopted in Chapter~\secIdea. After introducing the formal setup for tangles and their two main theorems in Chapters \secFormalSetup\ and~\secThms, we shall then return to these examples in Chapter~\secFormalExamples\ and see how they are made precise in that setup.

The examples given below are mostly from the social sciences. Since I am not an expert in any of these, I~shall not attempt more than to indicate the kind of use that tangles might find in these disciplines. Any actual use will require genuine expertise in the respective field, but my hope is that the examples below indicate enough of their potential to entice the experts to pursue this further.\looseness=-1

Further examples including text analysis, diagnostic expert systems, image processing or compression, and DNA or protein sequencing can be found in~\cite{TanglesEmpirical}.

The default format for describing the examples below will be to\dots
   {\smallskip\advance\parindent by 7mm
   \plainitem{$\bullet$} name the set $V\!$ of objects} studied;%
      \Footnote{As indicated in Chapter~\secIdea.5, for tangles in their most general framework there is not even a need for objects: tangles can be made to work with just a set $\vS$ and some axiomatic assumptions about~$\vS$ that reflect properties of `features of objects' sufficiently to enable us to formally define consistency. Then tangles can be defined as before: as consistent specifications of~$S$, and the two tangle theorems can be proved even in this very abstract framework~\cite{AbstractSepSys,ProfilesNew,TangleTreeAbstract}.\looseness=-1}
   {\advance\parindent by 7mm
   \plainitem{$\bullet$} describe the set $S$ of their potential features;
   \plainitem{$\bullet$} discuss when a specification of~$S$ is deemed to be typical;
   \plainitem{$\bullet$} describe what tangles~-- typical specifications~-- of~$S$ mean in this setup;
   \plainitem{$\bullet$} indicate some uses and applications of tangles in this context.
   \smallbreak}

When we return to these examples in Chapter~\secFormalExamples, we shall further describe what the two tangle theorems mean in the context of each of these examples: how they help us organize the tangles into a global structure (e.g., how some tangles refine others), how tangles can be distinguished by some small subset of particularly crucial questions or features, and how to obtain quantifiable evidence of the non-existence of tangles whenever our data is too unstructured.

Consistency is always defined as in Chapter~\secIdea.2. The definition of when a specification of~$S$ is deemed to be typical is always made by declaring some small sets of elements of~$\vS$ as atypical, and then calling a specification of~$S$ {\it typical\/} if it is consistent and contains none of these `forbidden' subsets. The collection of these `forbidden' atypical sets will usually be denoted by~$\F$.

Tangles, then, are typical specifications of~$S$ (as earlier): consistent sets of features that contain exactly one of $\vs$ and~$\sv$ for every $s\in S$ and have no subset in~$\F$.

\subsection \secInformalExamples.1\ \ Sociology: discovering mindsets, social groups, and character traits

Mindsets were the example of tangles we discussed in the introduction, where we have a set $V\!$ of people polled with a questionnaire~$S$.%
   \Footnote{In the language developed since then, the mindsets discussed there were tangles with witnessing sets. From now on, we shall use the word `mindset' for any tangle in a questionnaire scenario, regardless of whether that tangle has a witnessing set.}
   We may assume for our model that $S$ consists of yes/no questions: if it does not, we can still use it, and simply translate its answers into answers of an equivalent imaginary questionnaire of yes/no questions, which we then use as the basis for computing our tangles.%
   \Footnote{For example, if a question $s$ in $S$ asks for a numerical value between 1 and~5 on a scale from `do not agree at all' to `agree entirely', we can replace~$s$ with five imaginary yes/no questions $s_1,\dots,s_5$ asking whether the value of $s$ is $1\,,\dots,5$, respectively, or with the four yes/no questions $s'_i$ asking whether the value of $s$ is greater than~$i$, for $i=1,\dots,4$.}%
   \COMMENT{}

Specifying the collection~$\F$ of `forbidden' sets of features is our main tool for determining how broad or refined the mindsets will be that show up as tangles. If we forbid no feature sets, i.e., leave $\F$ empty, then every specification of~$S$ returned by one of the people polled will determine a tangle~-- because it will be consistent and hence typical if $\F=\es$.%
   \Footnote{This is not the same as saying that all possible specifications of~$S$ will be consistent (and hence typical). Specifications that contain an inconsistent triple will be inconsistent even then, and hence not be tangle. But these specifications are not among those returned in the poll~-- recall how consistency and inconsistent triples were defined.}

Such tangles would be too `fine' to be helpful. If we do not wish to influence which tangles are returned by our algorithms other than by determining how `broad' the mindsets they define should be, we can define $\F$ as in Chapter~\secIdea.3: as the set of all triples $\{\vr,\vs,\vt\}$ shared by no more than $n$ people polled, where $n$ is a number we can experiment with at the computer and see how it influences the number (and fineness) of tangles it returns.

But we can also decide to add some sets of features to~$\F$ that we think of as inconsistent because of what these features mean. What these are will be up to us: we might add sets of features that are intuitively inconsistent, or sets of features that will occur together as answers on the same questionnaire only if someone tampered with it or tried to influence the survey.

We might even use $\F$ to deliberately exclude some types of mindsets from showing up as tangles, e.g., mindsets we are simply not interested in. In order to exclude such mindsets we simply add to $\F$ the feature sets that make them uninteresting:%
   \COMMENT{}
   then no specifications of~$S$ that include these feature sets will show up as tangles.%
   \Footnote{For example, we might be interested in the opinions about hooligans among football crowds, but want to exclude hooligans themselves from this survey. Since we may not be able to identify them when we hand out the questionnaires, but know some answer patterns they are likely to give, we can add these patterns to~$\F$ to ensure that the tangles found are mindsets of spectators that are not themselves hooligans.}
   Of course, subsets we forbid for this reason must be specific enough that they occur {\it only\/} in mindsets that are not of interest to us.

We already discussed in the introduction what tangles mean in this example~-- namely, mindsets~-- and what the two main tangle theorems offer: a small set of critical questions that suffices to distinguish all the existing mindsets and on which predictions can be based, and verifiable evidence that no mindsets other than those found exist, possibly none. In Chapter~\secFormalExamples.1\ we shall add another aspect of the tree-of-tangles theorem: it enables us also to structure the mindsets found hierarchically, into broader mindsets and more focused ones refining these.

Discovering social groupings (not groups, see below) is similar to discovering mindsets, except that $S$ is different now: it can still be thought of as a questionnaire, but the `questions' it contains may be answered by an observer rather than the subjects $v\in V\!$ themselves. Thus, its specifications~$\vs$ will be more like the `features' discussed in our furniture example in Chapter~\secIdea. What makes tangles special in this application is again that they can identify such groupings without identifying actual groups in terms of their members: tangles are patterns of behaviour (etc.)\ often found together, not groups of similar people.\looseness=-1

While this may sound like a truism, it does differ fundamentally from the classical approach that seeks to find groups as clusters of people. Being defined directly in terms of the phenomena that define social groups, rather than indirectly in terms of the people that display them, tangles bypass many of the usual problems that come with traditional distance-based clustering, such as the fact that an individual is likely to belong to more than one social group. However, tangles can help even in finding those groups~-- which takes us on to character traits.

Imagine we wish to use tangles for a matchmaking algorithm. As before, they can help us to identify, from the answers people have given to a questionnaire~$S$ of character-related questions, some combinations of traits that are typically found together.%
   \Footnote{We have done this with data designed to test for the `big five' personality traits, and got some rather suprising results; see~\cite{PsychologyTangles}.}%
   \COMMENT{}
   Thus, tangles in this context will be `types of character'.

However, for our matchmaking algorithm it is not enough to know such types: it will also have to match individuals. So we will need a metric that tells us which pairs of people are `close in character'~-- assuming that matching like individuals is our aim (which, of course, can be disputed). The simplest example of a distance function on the set of people that returned our questionnaire would be to consider two people as close if they answered many questions identically.%
   \Footnote{For mathematicians: this is the Hamming distance in the hypercube~$2^S$.}
   But tangles can be used to define more subtle distance functions. We shall discuss some of these in Chapter~\secThms.4, and refer to them when we revisit the topic of matchmaking in Chapter~\secFormalExamples.1.

To summarize, let us emphasize again what makes the tangles approach to the study (and discovery) of mindsets, social groupings, or traits of character different from traditional approaches. It is that we can mechanically find these mindsets (etc.) from observational data without any prior intuitive hypothesis of what they may be, and that we can find them without having to group the people observed accordingly. However, knowing the tangles can then help us also to determine these groups; this will be discussed further in Chapter~\secThms.4.

\subsection \secInformalExamples.2\ \ Psychology: understanding the unfamiliar

In a way, this example is no more than a special case of the mindset example from~(\secInformalExamples.1). But the special context may lend it additional relevance.

While it is interesting to search for combinations of, say, political views that constitute hitherto unknown mindsets that can have an impact on political developments, it is not only interesting but crucial to try to {\it understand\/} minds that work in ways very different from our own. This is particularly rele\-vant in doctor-patient relationships where the doctor seeks to offer an individual with such a different mind a bridge to society, or even just to their particular environment: to enable them to understand the people around them, and to help these people to understand them.

Tangles can already help bridge this gap at the most fundamental level, the level of notions into which we organize our perceptions. We all have a notion, for example, of `threat'. But some patients'%
   \COMMENT{}
    notion of threat may be different from ours: they may be scared by things we would not see as threatening. And what these are may well come in types: typical combinations of perceptions of everyday phenomena which, for people with a certain psychiatric condition, may combine to a perception of a threat, and thereby form their notion of `threat' that may well differ from ours.

Tangles are designed to identify such types. The set $\vS$ would consist of various possible perceptions that experience has shown are rele\-vant, and a tangle would identify the combinations of these that are typical in the sense of Chapter~\secIdea.3. The definition of `typical' can be made sufficiently flexible to allow tangles to capture notions, unfamiliar to us, that consist of perceptions that typically occur together in some patients. Once these have been made explicit~-- recall that every tangle of~$S$ will be one specific set of possible perceptions~-- we can train our intuition on them in an effort to {\it understand\/} our patients rather than just collect lists of unrelated symptoms. We shall get back to how tangles can help identify meaning in Section~\secInformalExamples.5.

At a higher level, tangles can help to identify psychological syndromes as such, and maybe discover hitherto unknown syndromes. Indeed, psychological syndromes appear to be exactly the kind of thing that tangles model: collections of features that often occur together.%
   \COMMENT{}

This would be true for any medical condition, or even for mechanical `conditions' that lead to the failure of a machine. And indeed, there is a corresponding application of tangles for such cases, where they are used to build expert systems for medical (or mechanical) diagnosis~\cite{TanglesEmpirical}. But what makes psychiatric conditions even more amenable to the use of tangles is that the symptoms of which they are combinations are so much harder to quantify. Tangles come into their own particularly with input data that cannot be expected to be reliably precise.\looseness=-1 %
   \COMMENT{}

So how would we formalize the search for hitherto unknown psychological syndromes? Our ground set~$V\!$ would be a large pool of patients in some database. The set~$\vS$ would be a set of possible symptoms. Tangles will be collections of symptoms that typically occur together: medical, or psychological, conditions or illnesses.

We shall see in Chapter~\secFormalExamples.2\ that the set $T$ of tangle-distinguishing features returned by the tree-of-tangles theorem will consist of `critical' symptoms or combinations of symptoms that can be tested on a patient in the process of diagnosis. Every consistent specification of~$T$ defines a unique tangle: a unique condition that has all the symptoms in this specification of~$T$. In most cases, this will be the correct diagnosis.

Note, however, that our emphasis here is on {\it finding\/} psychiatric conditions: our aim is to discover combinations of symptoms that constitute an illness. It is not on diagnosing a given patient with one of these conditions (although checking the symptoms in~$T$ can play an important part on this), but on establishing what are the potential conditions to look for in the diagnosis. It seems to me that a lot has happened in psychology here in recent years, in that conditions are now recognized as illnesses that were not even thought of as typical combinations of symptoms (i.e., tangles) not so long ago.

\subsection \secInformalExamples.3\ \ Politics and society: appointing representative bodies

As discussed in Chapter~\secIdea, the tangle approach to clustering is that we do not primarily seek to divide our set~$V\!$ into groups based on similarity between its elements. This is true even if similarity is measured by how the elements of~$V\!$ specify~$S$, in which case we might group $u,v\in V\!$ together if they specify many $s\in S$ in the same way (as $\vs$ or~$\sv$): if $u(s) = v(s)$ for many~$s\in S$.

Rather, determining tangles is more about grouping features: not in the simple sense that we find traditional distance-based clusters in the set~$\vS$, but in that we look at how the elements of~$V\!$ specify~$S$ and find particular such specifications that are typical for the elements of~$V\!$. In the mindsets application from Section~\secInformalExamples.1, for example, tangles are ways of answering all the questions in~$S$ that are typical for how the various $v\in V\!$ answered them.

Once we have found the tangles of~$S$, however, we can use them to group the elements of~$V\!$ after all. In terms of the mindsets scenario, on which we shall base our further discussion here for better intuition, we would thus be looking for ways to group the elements of~$V\!$ according to their mindsets.

This can be done in various ways, which we shall look at in Chapter~\secThms.4. For example, we might associate each $v\in V\!$ with the mindset~$\tau$ of~$S$ that represents its views best: the tangle $\tau$ of~$S$ for which the number of elements~$s$ of~$S$ that $\tau$ specifies as $v$ does is maximum. Conversely we may seek, for each mindset~$\tau$ found, the person $v\in V\!$ closest to~$\tau$ in this sense, and think of these~$v$ as best representing the views held amongst the members of~$V\!$.

If we add our assumption from Section~\secInformalExamples.1 that $V\!$ was chosen so as to represent some larger population~$P$ well, we shall then have found a small group of people that ideally represent the views held in our population~$P$ on matters explored by~$S$. This process might be used, then, to appoint delegates to a body whose brief is to make decisions likely to find maximum consensus in~$P$.\looseness=-1

This contrasts with the usual democratic process of electing the delegates to such a decision-making body by majority vote. In a first-past-the-post system with constituencies, this can generate parliaments with large majorities even when, in each constituency, the majority of the successful candidate was small but the population is homogeneous enough across constituencies that in most of these the successful candidate is of the same political hue.

In systems with proportional representation this is avoided, but such systems require the previous establishment of political parties. Even when these exist, they may have developed historically in contexts that are less rele\-vant today. Finding the tangles of a political questionnaire is like {\it finding\/} the political parties that ought to exist today for the elections at hand. People $v\in V\!$ representing these virtual `parties' as described above could then be delegated to our decision-making body straight away.

If we wish to appoint more delegates for tangles, or virtual parties, with a larger following, we could elect representatives from these virtual parties by a standard vote using proportional representation. Alternatively we could simply refine this tangle into smaller tangles (see Chapters~\secFormalSetup.2 and~\secFormalExamples.3), so that all tangles end up with roughly the same amount of support%
   \COMMENT{}
   and could thus be represented by one delegate per tangle.

This approach may be even more rele\-vant outside politics, where there are many situations on a smaller scale in which we seek to appoint a decision-making body. This might be the governors for a school, or a steering committee for a choir. At such a smaller scale there will be no constituencies, and there may be no established parties rele\-vant to the brief of that body. But appointing for every tangle $\tau$ of~$S$ the member of~$V\!$ (amongst those willing to stand) whose views on~$S$ are closest to that tangle would produce a committee likely to represent the views held in~$V\!$ well.

\subsection \secInformalExamples.4\ \ Education: combining teaching techniques into methods

This application starts from the assumption that different teaching techniques work well or less well with different students: that a given technique is not necessarily better or worse than an alternative for all students at once, but that each may work better for different sets of students.

Ideally, then, each student should be taught by precisely the set of techniques that happen to work best for him or her. Of course, this is impractical: there are so many possible combinations of techniques that most students would end up sitting in a class of their own. But our aim could be to group techniques into, say, four or five groups, to be used in four or five classes held in parallel, so that each student can then attend the one of these four or five classes whose techniques suit him or her best.

The question then is: how do we divide the various techniques into the groups that correspond to the classes? Just to illustrate the problem, consider a very simple example. Suppose first that some students benefit most from supervised self-study while others are best served by a lecture followed by discussion in class. Suppose further that some students understand a grammatical rule best by first seeing motivating examples that prepare their intuition, while others perfer to see the rule stated clearly to begin with and examples only afterwards. So there are four possible combinations of techniques, but maybe we can only have two classes. So we have to group our techniques into pairs.

How shall we select the pairs? Shall we have one class whose teacher lectures and motivates rules by examples first, or shall we group the lecture-plus-discussion class with the technique of introducing rules before examples?

This is where tangles come into their own. Think of $V\!$ as a large set of students evaluated for our study, and of $\vS$ as a set of teaching techniques. Every tangle of~$S$ will be a particular combination of techniques (a~specification of~$S$) which, in the language of our furniture example in Chapter~\secIdea\ (where $\vS$ consists of features), `typically occur together'. What does this mean in our context? 

What it means formally depends on the definition of `typical' that is determined by our choice of the set~$\F$ introduced in Section~\secInformalExamples.1. Informally it implies, for example, that if the tangle is witnessed by a set~$X\sub V\!$ (see Chapter~\secIdea.5) then this $X$ is a group of students that would benefit particularly well from the combination of techniques specified by the tangle: for every technique in that tangle, a majority of the students in~$X$ prefer that technique to its converse. Thus, $X$~would be an ideal population for the class in which the techniques from this tangle are used, while a set~$Y$ witnessing another tangle would consist of students best served by the teaching techniques specified by that tangle.

Once more, speaking about tangles in terms of witnessing sets helps to visualize their benefits, but it is not crucial as such: the benefit arises from finding the tangles as such, and setting up those four or five classes accordingly. We then {\it know\/} that this serves our students best collectively~-- and can happily leave the choice of class to them.

Let us, from now on, refer to tangles of teaching techniques as {\it teaching methods\/}: combinations of techniques that work particularly well together for some students.

When we return to this example in Chapter~\secFormalExamples.3, we shall have to address the question of how to set our parameters in such a way that, even if the set~$\vS$ of potential teaching techniques is large, we still end up with a desired number of tangles~-- four or five in our case~-- that may be dictated by the formal environment, school etc., in which the teaching is taking place. We shall also see how, if we do not wish to let the students choose their class themselves, the tree-of-tangles theorem can help us devise an entry test that assigns students to the class that benefits them best.

\subsection \secInformalExamples.5\ \ Linguistics and analytic philosophy: how to determine meaning

The aim of this section is to indicate how tangles can help us identify {\it meaning\/}. Inasmuch as meaning is constituted in a social process this discussion belongs within the scope of this paper. It also has direct implications for the teaching of languages: in education, but also in machine learning. In Chapter~\secFormalExamples.5 we shall discuss an application in which tangles identifying the meaning of a word are used to steer a user of an interactive thesaurus towards the desired word.

To keep the discussion focused we shall concentrate on the meaning of words. There are obvious wider analogues, of course, and indeed the relevance of social aspects to the constitution of meaning will be more rele\-vant still when we talk about the meaning of entire phrases, perhaps depending on their contexts, or the meaning of behaviour in non-verbal communication.

Meaning of words, however, is not only easier to talk about: we have already done much of this when we discussed the furniture example in Chapter~\secIdea. There we observed that tangles can identify types of furniture, such as chairs, tables and beds, from lists of features of concrete specimen.

This contrasts with the na\"\i ve approach of trying to define words by a list of predicates, as classical dictionaries used to do.%
   \COMMENT{}
   The idea there is that something warrants being referred to by that word if and only if it satisfies all the predicates on the list.%
   \Footnote{\dots or, more generally, satisfies some logical formula in terms of these predicates}
   In other words, the list should be long enough that no things other than those we want to use our word for have all these properties, and it should at the same time be short enough that, conversely, all the things we wish to name by our word do indeed have all the properties from the list.

Put another way, if we assume that every predicate on the list describes a well-defined set of objects,%
   \Footnote{This is another dilemma with this na\"\i ve approach: it requires that there exists a hierarchy of predicates, where some are defined `before' others and can hence be used in the list for their definition. Tangles require no such hierarchy.}
   we would like the set of objects to be referred to by the word we are trying to define to be exactly the intersection of those sets: no larger and no smaller. As the furniture example demonstrates well, such a list is unlikely to exist: it is probably impossible to come up with a list of potential features (of furniture) such that the things that have all these features are precisely all the chairs.

The same problem arises more generally if we try to define the meaning of a word by any logical formula of previously defined predicates.%
   \COMMENT{}

Wittgenstein~\cite{WittgensteinPU} recognized this dilemma and argued that meaning cannot be captured by a simple taxonomic approach. Instead, he argued, the things referred to by a word behave more like members of a family: individuals who resemble each other but cannot be identified by any list of features that precisely they share.

Tangles offer a way to quantify Wittgenstein's family resemblances.%
   \Footnote{Thanks to Nathan Bowler for pointing this out to me.}
   The chair example in Chapter~\secIdea\ describes how.

In particular, tangles show that the intrinsic extensional%
   \COMMENT{}
   imprecision in most attempts to define the meaning of a word in terms of known%
   \COMMENT{}
   predicates simply by way of a logical formula does not make it futile to look for precise alternative ways to define meaning (in other ways than by logical formulas).%
   \Footnote{It seems that this difference has occasionally been overlooked.}
   As we shall see, there are a number of parameters we can use to quantify the notion of a tangle.%
   \Footnote{For example, by choosing the parameter~$n$ from our discussion of the term `typical' in Chapter~\secIdea.3~-- we shall meet this again in the sets~$\F_n$ in Chapter~\secFormalSetup.1~-- or by choosing an order function as described in Chapter~\secFormalSetup.4.}
   So even if we take the view that the notion of `meaning' can be formalized as `tangle of predicates', as I would argue, there will not be one ultimate such notion.%
   \COMMENT{}

Rather, there are many precise ones, one for every choice of our quantitative tangle parameters. Then in every `context'~$S$, i.e., for every list $\vS$ of predicates deemed rele\-vant to the definitions we are trying to make, every tangle of~$S$ constitutes an instance of that notion of `meaning': a~potential definition of what one particular word means.%
  \Footnote{In our earlier furniture example, these were the meanings of `chair', `table', `bed' and so on. By our choice of quantitative tangle parameters we can influence how broad or narrow that meaning is going to be.}
   Each of these notions will be completely precise: {\it tangles are a~precise way to define extensionally imprecise meanings of words.}

After this rather abstract discussion let us return once more to our concrete example, and summarize how exactly the notion of `chair' is given by a tangle. We start with a list $S$ of potential features of furniture. Any tangle of~$S$ is one particular specification of~$S$: a choice of either $\vs$ or~$\sv$ for every $s\in S$, either confirming the feature~$\vs$ or confirming its converse~$\sv$. The combination of features that we might think of as decribing `the perfect chair' would be a tangle~-- no matter whether such a chair exists in the world or not.

Although every specific tangle in this context is a concrete list of predicates, the notion of `tangle' as such is purely formal and makes no reference to such predicates: a {\it tangle\/} is a specification of~$S$ that satisfies certain formal requirements~-- consistency and being typical~-- which say that it must not contain certain small sets of features (which in our context are predicates). These requirements ensure that tangles define `good' notions rather than `bad' ones. The settings of certain parameters in these requirments also determine whether the corresponding tangles capture broad or narrow notions.

Thus, whether or not a combination of potential features (a~specification of~$S$) is a tangle is decided by a precise definition. This does not imply, however, that it must be clear for every piece of furniture whether or not it `belongs to' a given tangle, e.g., whether it is a chair. Indeed this is {\it not\/} clear for real pieces of furniture in the real world, and it is one of the strengths of our approach that tangles offer a precise way to capture even extensionally imprecise notions such as this one.

\subsection \secInformalExamples.6\ \ Economics: identifying customer and product types

The context of economics offers a wide range of interesting choices for our set $V\!$ of `objects' studied together with some potential `features' listed in a set~$S$. This is illustrated particularly well by a pair of examples described below, which describes two complementary, or `dual', aspects of the same scenario.

In our first example, let $V\!$ be a set of customers of an online shop, and $S$ the set of items sold at this shop. Let us assume that each customer~$v$ makes a single visit to the shop, specifying $s\in S$ as~$\vs$ if $v$ includes~$s$ in his or her purchase, and as~$\sv$ if not. We shall think of the specification $v(S)$ of~$S$ as $v$'s `shopping basket', dividing as it does the set of items into those bought and those not bought.%
   \COMMENT{}

    An arbitrary specification of~$S$, then, is a hypothetical shopping basket, and a tangle of~$S$ is a {\it typical\/} shopping basket for this set of customers.%
   \Footnote{This is not the same as just a group of items `often bought together', as are already suggested by some online shops today. A~shopping basket defined by a tangle is typical in a more subtle way; for example, it is also typical in what it does {\it not\/} include.}%
   \COMMENT{}
   If we mentally identify\vadjust{\penalty-200} a customer~$v$ with his or her shopping basket~$v(S)$, we may also think of a tangle of~$S$ as a (hypothetical) {\it type of customer\/} whose `features' are his or her purchases $v(s)=\vs$ and non-purchases $v(s)=\sv$.%
   \Footnote{This is analogous to thinking of the chair tangle as a type of furniture, or of a mindset tangle as a type of person holding such views. In each case, the tangle is a set of features typical for the elements of~$V\!$, and we may think of it as a hypothetical `typical' element of~$V\!$ that has exactly these features.}

Alternatively, we may think of~$S$, the set of items of our shop, as our set of {\it objects\/} (which would normally be denoted by~$V\!$, but this letter is taken now), and of the customers in~$V\!$ as potential `features' of these items: a~customer $v$ becomes a feature~$\vv$ of precisely those items~$s$ that $v$~bought, and a feature~$\vvback$ of those items~$s$ that $v$ did not buy. Each item $s$ then defines a specification~$s(V)$ of~$V\!$, specifying those $v$ that bought it as $s(v):=\vfwd$, say, and those that did not buy it as~$s(v):=\vbwd$. We may think of $s(V)$ as something like the `popularity score' of the item~$s$ with the customers in~$V\!$.

Note that the information encoded here is no different from that described by our earlier setup, where every customer specified each item~$s$ as~$\vs$ or~$\sv\!$. We thus have two ways now of describing the same set of customer preferences.%
   \Footnote{For mathematicians: there is a formal duality here in that $v(s)=s(v)$ with the obvious interpretations of $v$ as a function $S\to\{0,1\}$ and $s$ as a function $V\!\to\{0,1\}$. The customer preferences described in two ways can be formalized as the edge set of the bipartite graph with vertex classes $V\!$ and~$S$ and edges $vs$ whenever $v$ bought~$s$. This edge set can be described alternatively as a list of neighbourhoods of the vertices in~$V\!$ or of those in~$S$.}%
   \COMMENT{}

What are the tangles in our second setup? They are tangles, or typical specifications, of~$V$: hypothetical popularity scores with the customers in~$V\!$ that are typical for the items in~$S$. If we mentally identify an item~$s$ with its popularity score~$s(V)$, we may think of a tangle of~$V\!$ as a {\it type of item\/} whose `features' are its fans $s(v)=\vv$ and non-fans~$s(v)=\vvback$.

Let us look at some examples of both kinds of tangle, and of how they might be used.

Tangles of~$S$ are typical (if hypothetical) shopping baskets, collections of items that are typically bought or avoided together. For example, there might be a tangle, or `typical shopping basket', full of ecological items but containing no environmentally harmful ones. Another might contain mostly inexpensive items and avoid luxury ones.%
   \Footnote{In our current informal setup, these would be tangles of subsets of~$S$, not of $S$ itself~-- for example, of the subset of the ecologically critical or the unusually priced items. When we revisit this example in Chapter~\secFormalExamples.6, we shall find a way of capturing all groups of items, such as the ecological or the inexpensive ones, by tangles of the same set.}%
   \COMMENT{}

Note that `inverting' these tangles will not normally give a tangle. Indeed, a shopping basket full of luxury goods and avoiding cheap items is unlikely to be typical, because customers that like, and can afford, luxury good will not necessarily shun inexpensive items. Similarly, inverting the ecological tangle will not produce another tangle, since there is no unifying motivation amongst shoppers to buy environmentally damaging goods, and probably no unifying motivation to avoid `green' products either.

If we personify these tangles as indicated earlier, and think of them as (hypothetical) people putting together these hypothetical shopping baskets, we could think of our second tangle as the budget-oriented customer type, one that prefers inexpensive brands and cannot afford or dislikes unnecessarily expensive ones, and of the first as the ecological type that prefers organic foods and degradable detergents but avoids items wrapped in plastic. And, crucially, our tangle analysis might throw up some unexpected customer types as well~-- perhaps one that prefers items showing the picture of a person on the packaging.

Applications of tangles of~$S$ might include strategies for grouping goods in a physical shop, or running advertising campaigns targeted at different types of customer.

Tangles of~$V\!$, on the other hand, are typical popularity scores, ways of dividing $V\!$ into fans and non-fans that occur for many items simultaneously. If what interests us about the items in~$S$ is mainly how they appeal to customers, we may think of tangles of~$V\!$ as types of items.

For example, assume that the ecological goods in our shop are substantially more expensive than non-green competing goods. Then $V\!$ splits neatly into a `green fan set' of customers, those likely to buy these green items even though they are pricier, and the rest of all customers, who will avoid them because they are more expensive. This division of~$V\!$, or hypothetical popularity score, is borne out in sufficient numbers%
   \Footnote{This is a reference to the number~$n$ used in Chapter~\secIdea.3 to define tangles informally. It will be referred to again when we define tangles formally in Chapter~\secFormalSetup.1.}
   by the green items in~$S$ to form a tangle of~$V\!$, since these items are both liked by the `green' customers and disliked by the others for their price.%
   \COMMENT{}

An application of finding tangles of~$V\!$ could be to set up a discussion forum for each tangle and invite the fans for that tangle%
   \COMMENT{}
   to join. Since they have shown similar shopping tastes, chances are they might benefit more from hearing each others views than could be expected for an arbitrary discussion group of customers.

One last remark, something the reader may be puzzling over at this point. Since our two types of tangle are `dual' to each other, it so happens that witnessing sets for the first type are the same kinds of sets as tangles of the second type: they can both be thought of as sets of customers. But they are not the same.
A set $X\sub V\!$ witnessing a tangle of~$S$ is a set of customers such that for every item $s\in S$ either most of the people in~$X$ bought~$s$ (if the tangle specifies~$s$ as~$\vs$) or most of the people in~$X$ did not buy~$s$ (if the tangle specifies~$s$ as~$\sv$). By contrast, a tangle of~$V\!$ (thought of as the set of those~$v$ which it specifies as~$\vfwd$)%
   \Footnote{Unlike in other contexts, every $v$ in this example has a default specification~$\vv$: the set of purchases of~$v$ rather than his non-purchases. The set $\{\vv\mid v\in V\}$, therefore, is well defined.\looseness=-1}
    is a group of customers such that every three of them jointly bought some sizable set of items (at least~$n$).%
   \Footnote{Indeed, more is true: for every three customers~$v$, regardless of how the tangle~$\tau$ specified them, there exists a set of at least $n$ items that are each bought by $v$ if $\tau(v)=\vv$ {\it and\/} not bought by~$v$ if $\tau(v)=\vvback$.}%
   \COMMENT{}

Put more succinctly: in the first case, most people in some set of customers (not too few) agree about every single item, while in the second case the tastes of every few (e.g., three) customers are witnessed by some common set of items (not too few), regardless of whether their tastes coincide on these items or not.%
   \COMMENT{}

\beginsection \secFormalSetup. The formal setup for tangles

The aim of this chapter is to give a more formal definition of tangles: just formal enough to enable readers to try them out in their work if they wish.%
   \Footnote{If you do, please feel free to contact me for support at any time.}
   The definitions here will be less general than is possible, but general enough to model all the examples of tangles we shall be looking at here. More general definitions, for which all the theorems still hold, can be found in~\cite{AbstractSepSys}.

We shall make precise the following notions introduced informally so far:
   {\smallskip\advance\parindent by 7mm
   \plainitem{$\bullet$} {\it (potential) features\/} $s$ of (elements $v$ of)~$V\!$
   \plainitem{$\bullet$} {\it specifications\/} $\vs,\sv$ of such $s$
   \plainitem{$\bullet$} {\it consistency\/} of specifications
   \plainitem{$\bullet$} {\it types\/} (or {\it typical\/} specifications of~$S$)
   \plainitem{$\bullet$} {\it tangles\/}
   \smallbreak}
\noindent
   Readers consulting the mathematical papers on tangles referenced here should be aware that, for historical reasons, these things have different names there:
   {\smallskip\advance\parindent by 7mm
   \plainitem{$\bullet$} Features of (elements of) $V\!$ are called {\it separations\/} of~$V\!$.
   \plainitem{$\bullet$} Specifications of features are {\it orientations\/} of those separations.
   \plainitem{$\bullet$} {\it Consistency\/} is defined more narrowly than here: specifications of~$S$ that are consistent in the sense defined formally below are treated in~\cite{ProfilesNew} as {\it robust profiles\/}, and in~\cite{AbstractTangles} as {\it abstract tangles\/}.
   \plainitem{$\bullet$} Typical specifications of $S$ are those that {\it avoid\/} some given collection $\F$ of subsets of~$\vS$, in that they have no subset that is an element of~$\F$.
   \plainitem{$\bullet$} {\it Tangles\/} are the same as here: they are typical consistent specifications of~$S$. If reference is desired to the set $\F$ underlying the notion of `typical' used for a given tangle, it is called an $\F$-tangle.%
   \Footnote{When tangles are mentioned without any $\F$ being specified, then this $\F$ is either arbitrary, or it is some specific $\F$ traditionally considered in the given context.}
   \plainitem{$\bullet$} The operations $\land$ and~$\lor$ have their meanings swapped.
   \smallbreak}
\noindent
   In the interest of readability, these traditional names will not be used in the rest of this paper.

The two main tangle theorems will be stated more formally in the next chapter, in terms of the notions introduced below. The discussion of {\it predictions\/} is also deferred to the next chapter, since predictions based on tangles make reference to the first of these theorems, the {\it tree-of-tangles\/} theorem.

Also included in this chapter will be some basic facts about tangles that can help the reader develop an intuition for them. We shall see how some tangles evolve out of others, and die again to spawn new tangles~-- just as the detailed notions of `armchair' and `dining chair' refine the more basic notion of `chair', which at their level of detail recedes and becomes less meaningful.

Importantly, we discuss how to {\it influence\/} this `evolution of tangles' according to our needs and aims, by organizing features into hierarchies of `orders'.

We shall also make precise the notion of {\it duality\/} for tangles, of which we already caught a glimpse in Chapter~\secInformalExamples.6. This, too, adds substantially to our toolkit for applications.

\subsection \secFormalSetup.1\ \ Tangles of set partitions

A {\it partition\/} of a set~$V\!$, for the purpose of this paper, is an unordered pair $\{A,B\}$ of disjoint subsets of~$V\!$, possibly empty, whose union is~$V\!$. Its two elements, the complementary sets $A$ and~$B$, are its {\it specifications\/}.%
   \Footnote{For example, $V\!$ might be a set of people whom we asked a yes/no question. Then $V\!$ is made up of the set $A$ of the people that answered yes and the complementary set $B$ of the people that answered no. The partition $\{A,B\}$ still records the information of how the vote was split across~$V\!$, but `forgets' which of its two elements corresponds to which vote.\looseness=-1}
   The sets $A$ and $B = V\sm A$, which we also denote by~$\Abar$, are {\it inverses\/} of each other.

When $S$ is a set of partitions of~$V\!$, we denote the specifications of its elements~$s$ as $\vs$ and~$\sv$,%
   \Footnote{It is meaningless here to ask which is which. For example, if $s=\{A,B\}$ as above, whether $\vs=A$ and $\sv = B$ or $\vs=B$ and $\sv = A$. This is because the two elements of $s=\{A,B\}$ are `born equal': the expressions `$\{A,B\}$' and `$\{B,A\}$' denote the same set or unordered pair. If this did not worry you, let it not worry you now.}
   and write $\vS$ for the set of all the specifications of its elements. A~subset of $\vS$ is a {\it specification of~$S$\/} if for every $s\in S$ it contains either $\vs$ or~$\sv$ but not both.

The elements of~$S$ are called {\it potential features\/} of (the elements of)~$V\!$, those of $\vS$ are their {\it features\/}. An element of~$V\!$ {\it has\/} a feature~$A$ if it lies in~$A$. Given two features $\vr = A$ and $\vs = C$, we write
 $$\vr\lor\vs := A\cup C \quad\hbox{and}\quad \vr\land\vs := A\cap C.$$
 Note that $\vr\lor\vs$ and~$\vr\land\vs$ need not be elements of~$\vS$, i.e., features of~$V\!$ in terms of our given set $S$ of potential features. But they become features if we add to $S$ the partitions $\{A\cup C, \overline{A\cup C}\}$ and $\{A\cap C, \overline{A\cap C}\}$ as potential features.

We read $\lor$ as `or' and $\land$ as `and'. This corresponds to the fact that the elements of~$V\!$ that have the feature $\vr\lor\vs$ (if it is a feature) are those in~$A\cup C$, which are precisely those elements of~$V\!$ that have the feature~$\vr$ or the feature~$\vs$ (or~both), and similarly with~$\land$.

A subset of~$\vS$ is {\it inconsistent\/} if it contains three%
   \Footnote{These features $\vr,\vs,\vt$ need not be distinct. For example, our subset of~$\vS$ will already be inconsistent if it contains two features $\vr,\vs$ such that $\vr\land\vs=\es$; this is included in the definition given, simply by taking $\vt = \vs$ or $\vt=\vr$.}
   features $\vr,\vs,\vt$ such that no element of~$V\!$ has them all, i.e., such that $\vr\land\vs\land\vt = \es$. Otherwise the subset is {\it consistent\/}. For example, every feature $\vs=A$ is inconsistent with its inverse $\sv=\Abar$, because $\vs\land\sv = A\cap \Abar = \es$. Similarly, $\vr$, $\vs$ and the inverse of~$\vr\land\vs$ form an inconsistent triple of features.%
   \Footnote{Indeed, if $\vr=A$ and $\vs=C$, say, then the inverse of $\vr\land\vs$ is $\overline{A\cap C}$. Let us denote this as~$\vt$. Then $\{\vr,\vs,\vt\}$ is inconsistent, because $\vr\land\vs\land\vt = A\cap C\cap (\overline{A\cap C}) = \es$.}

Given a set $\F$ of subsets of~$\vS$, a~consistent specification of~$S$ that has no subset in~$\F$ is called {\it $\F$-typical\/}, or an {\it $\F$-tangle\/} of~$S$. Given an integer~$n$, we write $\F_n$ for the set of all sets $\{\vr,\vs,\vt\}$ of up to three features that no $n$ elements of~$V\!$ have in common, i.e., which satisfy $\vr\land\vs\land\vt = U$ for some $U\sub V\!$ with $|U| < n$:
 $$\F_n = \{\,\{A,B,C\}\sub \vS : |A\cap B\cap C| < n\,\}.$$
 Thus, a~specification $\tau$ of~$S$ is an $\F_n$-tangle if for every triple of features in~$\tau$ there are at least $n$ elements of~$V\!$ that share these three features.%
   \Footnote{In some contexts%
   \COMMENT{}
   it may be more appropriate to require this only on average: that the number of elements of~$V\!$ sharing a given triple of features in~$\tau$ is at least~$n$ on average, taken over all the triples in~$\tau$. We call such $\tau$ {\it averaging $\F_n$-tangles\/}. We might further require that $\tau$ is {\it stable\/}: locally optimal in the sense that replacing any one of its elements with its inverse reduces the average size of the intersections of the triples in~$\tau$.}%
   \COMMENT{}

A~{\it typical\/} specification, or {\it tangle\/}, of~$S$ is an $\F_n$-tangle for some $n$ specified in the context, possibly as a function of~$|V|$ or~$|S|$.%
   \COMMENT{}
   If no $n$ is given explicitly it is assumed to be~1, in which case `typical' is no more than `consistent'.%
   \Footnote{Indeed, $\F_n$-tangles for any~$n$ are also $\F_m$-tangles for all $m<n$, since $\F_m\sub\F_n$. In particular, they are $\F_1$-tangles, but $\F_1$ is just the set of inconsistent triples in~$\vS$.}

Let us look at a couple of very simple tangles. The simplest tangles are the $\F_1$-tangles that are `focused on' some fixed $v\in V$:
 $$\tau_v = \{\,\vs\in\vS \mid v\in\vs\,\}.$$
   These tangles exist for every~$S$ and~$v$, so their existence does not tell us anything about~$S$ (which, of course, is our aim in studying its tangles). Thus, $\F_n$-tangles are interesting only for $n\ge 2$.

Every $\F_n$-tangle of~$S$ contains all the features~$\sv$ for which $\vs$ has fewer than $n$ elements, because it cannot contain $\{\vs\}\in\F_n$ as a subset.

If $S$ is the set of {\it all\/} bipartitions of~$V\!$, then the only tangles $\tau$ of~$S$ are the focused tangles~$\tau_v$.%
   \Footnote{For the mathematically inclined, here is the easy and pretty proof. If we can find any $v\in V\!$ such that $\{v\}\in\tau$, then $\tau=\tau_v$, because any $A\in\tau$ not containing~$v$ would be inconsistent with~$\{v\}\in\tau$. To find~$v$, consider any $A\in\tau$. Then $A\ne\es$, since $\tau$ is consistent. If $|A| = 1$, we have found~$v$. If $|A|\ge 2$, we partition $A$ into two non-empty sets $B$ and~$C$. One of these sets must be in~$\tau$, as otherwise $\tau$ would contain the inconsistent set $\{A,\Bbar,\Cbar\}$. We then partition this set and proceed as we did earlier with~$A$, until we arrive at some $\{v\}\in\tau$.}%
   \COMMENT{}
   Since $\{\{v\},\overline{\{v\}}\}\in S$, this means that $\{v\}\in\tau$, and hence that $\tau$ is not an $\F_2$-tangle. The set $S$ of all bipartitions of~$V\!$, therefore, has no $\F_n$-tangles for $n\ge 2$. For this reason, we shall not be interested in tangles of this entire set~$S$.

\subsection \secFormalSetup.2\ \ The evolution of tangles

However, we do often consider the set $S$ of all partitions of~$V\!$ in increments: we study tangles of subsets $S_1\sub S_2\sub\dots$ of~$S$ that ultimately exhaust~$S$. These $S_k$ have to be carefully chosen for this to make sense, and we look at how to do this in Section~\secFormalSetup.4. But if they are, the tangles of the various~$S_k$ can interact in interesting ways, and studying this interaction adds greatly to our understanding of how they all together describe the features of~$V\!$.

To see how this works, note that for $i<k$ every tangle $\tau_k$ of~$S_k$ {\it induces\/} a tangle $\tau_i = \tau_k\cap\vSi$ of~$S_i$. Indeed, since $S_i$ is a subset of~$S_k$, each of its elements $s$ is specified by~$\tau_k$ as $\vs$ or~$\sv$, and the set of all these specifications of elements of~$S_i$ is consistent (because its superset $\tau_k$~is). Similarly, if $\tau_k$ is an $\F$-tangle for some~$\F$ then so is~$\tau_i$, since any subset of~$\tau_i$ in $\F$ would also be a subset of~$\tau_k$ in~$\F$ (which by assumption does not exist).

Conversely, however, a tangle $\tau$ of~$S_i$ need not {\it extend to\/} a tangle of~$S_k$ in this way: there may, but does not have to, be a tangle of~$S_k$ whose elements in~$\vSi$ are exactly~$\tau$. This is because extending $\tau$ to a tangle of~$S_k$ requires that we specify all the $s\in S_k$ that are not already in~$S_i$ as well; and there may simply not be a consistent way of doing that, or one that avoids creating a subset that lies in our chosen~$\F$.

\goodbreak

Similarly, a tangle of~$S_i$ can `spawn' several tangles of~$S_k$ that all extend it.%
   \Footnote{For example, if a tangle of order~$k$ identifies the chairs in our set $V\!$ of furniture, we might find that $S_{k+1}$ has three tangles that induce this chair tangle of order~$k$: one for armchairs, another for dining chairs, and perhaps a third for garden chairs.}
   We can thus study how tangles are `born' at some time~$i$ (when they are one of several tangles of~$S_i$ spawned by the same tangle of~$S_{i-1}$), then `live on as' (extend to unique) tangles of~$S_k$ for some $k>i$, and eventually `die' at some time $\ell > k > i$ (when they spawn several tangles of~$S_\ell$, or extend to none).\looseness=-1

\subsection \secFormalSetup.3\ \ Stars, universes, and submodularity

A set of features is a {\it star\/} if their inverses are disjoint subsets of~$V\!$. The sets $\F$ that we use to define tangles will often consist of stars. They also arise naturally as in the following example.

Consider a survey $S$ whose questions allow answers on a scale from 1~to~5. To compute its tangles, we first need to convert $S$ into a hypothetical questionnaire $S'$ of yes/no questions that is equivalent to~$S$ in the sense that its answers can be computed from the (known) answers recorded for~$S$, and whose (computed) yes/no answers conversely imply the original answers to~$S$. For each $s\in S$ we let $S'$ contain five questions $s_i$ (where $i=1,\dots,5$) asking whether the given answer to~$s$ is the $i$th of the five potential answers, yes or~no. Let $s_i = \{A_i,\Abar_i\}$, where $\vsi=A_i$ consists of the yes-answers. These sets~$A_i$ are disjoint, because every $v\in V\!$ had to choose only one of the five possible answers to~$s$, which makes him an element of only one~$A_i$. So the $\Abar_i$ have disjoints complements, which means that $\{\svone,\dots,\svfive\}$ is a star.

We shall be particularly interested in sets $S$ such that $\vS$ contains, for any given $\vr,\vs\in\vS$, also the features $\vr\lor\vs$ and~$\vr\land\vs$. If $S$ has this property, then $\vS$ will, in fact, contain {\it all\/} Boolean expressions of features it contains already:%
   \Footnote{Recall that~$\vS$, the set of all specifications of elements of~$S$, already contains the inverses of all its elements.}
   combinations of features built from others by using the symbols $\lor$ and~$\land$ and taking inverses.%
   \COMMENT{}%
   \Footnote{For example, $\rv\land (\vr\lor\vs)$ and its inverse are such `Boolean expressions'.}
   We then call $\vS$ a {\it universe\/} of features.

If we are given a set~$S$ of potential features that lacks this property~-- for example, when $S$ is given as a questionnaire%
   \Footnote{From now on, we shall assume in all examples where $S$ `is a questionnaire' that the people in~$V\!$ have already answered all its questions. Only then do the elements of~$S$ become potential features in our new formal sense, i.e., partitions of~$V$: these partitions are, for every `question'~$s$, given as the sets $A$ and $B$ of the people in~$V\!$ that answered $s$ as yes or as no, respectively, so that $s$ becomes the partition~$\{A,B\}$.}~--
   we can extend~$S$ by adding to~$\vS$ any missing features and their inverses until it does have this property. The set $S'$ of potential features generated from $S$ in this way will be much larger than~$S$ and may well contain all bipartitions of~$V\!$.%
   \COMMENT{}
   We call $\vSdash$ the feature universe {\it generated by\/}~$S$.

Finally, $S$~and~$\vS$ are {\it submodular\/} if for any two features $\vr,\vs\in\vS$ at least one of $\vr\lor\vs$ and $\vr\land\vs$ is also in~$\vS$. Submodularity is a technical `richness' condition on~$S$ that is needed in the proofs of our two tangle theorems, and which helps with our tangle algorithms. If $S$ is not submodular, we can make it so by adding the required potential features, as earlier when we generated a universe from~$S$.%
   \Footnote{However, the submodular set $S'\supe S$ is often smaller than the universe generated by~$S$.}%
  \COMMENT{}
   Any tangle of $S'$ then also has to specify the newly added features. Tangles of sets $S$ that are not submodular do not automatically extend to such larger tangles.%
   \COMMENT{}

\subsection \secFormalSetup.4\ \ Hierarchies of features and order

In Section~\secFormalSetup.2 we described how tangles evolve from each other when they are tangles not of the entire set $S$ of partitions of~$V\!$ we are considering,%
   \Footnote{This $S$ may be the set of all partitions of~$V\!$, or just of some that are of particular interest.\looseness=-1}%
   \COMMENT{}
   but of increasing subsets $S_1\sub S_2\sub\dots$ of~$S$. Whether or not this formal evolution of tangles corresponds to any meaningful refinement of typical features, however, will depend on how we set up this stratification $S_1\sub S_2\sub\dots$ of~$S$.

Technically, we shall do this as follows. To every $s\in S$ we shall assign an integer, to be called the {\it order\/} of~$s$ (and of $\vs$ and of~$\sv$) and denoted as $|s|$, so that%
   \Footnote{Thus, when we choose~$S_1$ we assign its elements the order~0. When we increase $S_1$ to~$S_2$, we assign the new elements~-- those in~$S_2\sm S_1$~-- the order~1. And so on, as $k=1,2,\dots\,$.}
 $$S_k = \{\,s\in S: |s| < k\,\}.$$
 Informally, the tangles of~$S_k$ are also often said to have `order~$k$'.

The function $s\mapsto |s|$ is called an {\it order function\/}. We sometimes require it to be {\it submodular\/} on~$\vS$, which means that $|\vr\lor\vs| + |\vr\land\vs| \le |\vr|+|\vs|$ for all $\vr,\vs\in\vS$. Note that this implies that the sets $S_k$ are submodular as defined in Section~\secFormalSetup.3.

Ideally, then, the typical features encoded in a tangle of any given order should `refine' those of tangles of lower order in some natural sense given by the application at hand. We might hope to achieve this by choosing for $S_1$ only the most basic potential features of the objects in~$V\!$, for~$S_2$ some slightly less basic features, and so on until the $S_k$ for large~$k$ consist of quite specific potential features that are only needed to distinguish objects not distinguished by any of the more basic features.%
   \Footnote{For example, the potential feature of being collapsible or not may serve to distinguish some garden or dining chairs from others. It might then split the corresponding tangles in two, giving rise to a tangle of camping chairs and another of rigid garden chairs, and similarly for the dining chairs. The armchair tangle would be unaffected by this additional potential feature, since armchairs are all non-collapsible, and so the feature of not being collapsible will simply integrate into the existing armchair tangle.}
    In a questionnaire environment, $S_1$~would consist of the most basic questions in the survey, and as $k$ increases, the questions become more and more specific.

In a more subtle approach, we might defer the choice of which potential features to include in~$S_{k+1}$ until we have computed the tangles of~$S_k$.%
   \Footnote{Our standard algorithms do indeed compute the tangles of~$S_k$ successively as $k$ increases.}
   In our questionnaire scenario, for example, the question of whether someone likes Tchaikovsky may be rele\-vant when we are looking for how to split a tangle of love for classical music, but less rele\-vant when we analyse a tangle of sport enthusiasm. With this approach, there will be several sets~$S_{k+1}$, one for every choice of a $k$-tangle (which in turn may be tangles of different sets~$S_k$), and the same feature may be given different orders depending on the choice of the tangle relative to which this order is assigned.%
   \Footnote{Our question about Tchaikovsky, for example, might receive order~$k$ when we consider the $k$-tangle of classical music enthusiasm, while also receiving order $k+5$ relative to some $(k+5)$-tangle extending the $k$-tangle for sports, such as a tangle for the love of ballet.}

While defining the order of potential features explicitly appears to be the most natural option, there are a few issues we should be aware of. One is that sometimes we simply do not know which potantial features are `basic' and which are not. This may be because our set $S$ is not a hand-designed survey but a collection of technical properties of the objects in~$V\!$ resulting from some measurements whose significance we cannot judge.

A more subtle issue is that by declaring some potential features as more fundamental than others we are in danger of influencing the outcome of our study more than we would like to. Indeed, which tangles arise in the~$S_k$ has everything to do with how these~$S_k$ are chosen. If we are searching for hitherto unknown types of features, we should keep our influence on how to group them to a minimum, unless we know what we are doing and do it deliberately.

Finally, there is a technical problem with choosing the $S_k$ explicitly. This is that we would like them to be submodular (see Section~\secFormalSetup.3): for any $\vr,\vs\in\vSk$ we would like at least one of $\vr\lor\vs$ and $\vr\land\vs$ to be in~$\vSk$ too. We might simply add one of them to~$\vSk$ to achieve this, but this can cause other problems. One is that $\vr\lor\vs$ and $\vr\land\vs$, while being well-defined partitions of~$V\!$, need not lie in~$S$ and therefore may lack the `meaning' that comes with the element of~$S$, making it difficult to decide which of them is the appropriate one to add to~$S_k$ because it is more fundamental, or basic, than the other.

We shall return to this problem later in this section. Choosing the order of potential features explicitly remains a valuable aim in many scenarios, and there is a way to deal with the technical problem just described.

There are also some generic ways of defining an order function, ways that make no reference to the interpretation of the potential features but which do take into account the abstract way in which these features structure the set~$V\!$, i.e., how the elements of~$S$ partition~$V\!$. Moreover, just as the set~$S$ of all our given potential features can generally determine the order of a particular partition $s\in S$ (as we shall see below), it can determine the order of an arbitrary partition $\{A,B\}$ of~$V\!$.

In the rest of this section we shall look at an example of such a generic way of assigning an order $|\{A,B\}|$ (which we shall abbreviate to~$|A,B|$) to all the partitions of~$V\!$ by appealing to~$S$ in one way or another. It is only an example intended to demonstrate the idea; another example%
   \COMMENT{}
   is given in~\cite{TanglesEmpirical}. On the whole, there is scope for some substantial and interesting mathematics here, some of which is explored in~\cite{TangleOrder}.

To understand better our task of assigning an order $|A,B|$ to an arbitrary partition $\{A,B\}$ of~$V\!$, let us look at cases where we already have an aim for how to do this: when $\{A,B\}$ is itself in~$S$. In that case we said~-- appealing to an interpretation of the elements of~$S$, which we shall now seek to avoid~-- that $|A,B|$ should be small if its complementary features $A$ and~$B$ are `basic', or `fundamental'. Let us try to capture this notion without interpretation, but allowing ourselves to refer to how the other elements of~$S$ partition~$V\!$.

One obvious aspect of $s=\{A,B\}$ being `fundamental' is that it does not divide $V\!$ in such a way that large groups of `similar' objects are split more or less evenly: that about half of them lie in~$A$ and the other half in~$B$.%
   \Footnote{For example, the feature of having (long) legs is fundamental in this sense for furniture: most chairs have legs, most tables have legs, and most beds do not have legs. The feature of being made almost entirely of wood is not fundamental in this sense: there are about as many chairs and tables that have this feature as there are chairs and tables that do not have it.\looseness=-1}
   This can be expressed in an abstract way~-- with reference to~$S$, but without interpretation~-- if being `similar' can be. And indeed it can: we may think of two objects in~$V\!$ as similar if they share many of the features in~$S$, if for many $s\in S$ they both specify $s$ as~$\vs$ or they both specify it as~$\sv$.

So let us define a similarity function $\sigma\: (u,v)\mapsto \sigma(u,v)\in\N$ on~$V\times V\!$, the set of pairs of objects, by counting common features in this way, letting
  $$\sigma(u,v) := \big|\{\,s\in S : u(s) = v(s)\,\}\big|\,,$$
   where
 $$ v(s) := \cases{ \vs & if $v$ specifies $s$ as~$\vs$\cr
                    \sv & if $v$ specifies $s$ as~$\sv$\cr}$$
 and $u(s)$ is defined analogously. Thus, $u$ and~$v$ have large {\it similarity\/}~$\sigma(u,v)$ if they share many of the features in~$\vS$.

Note that there is ample scope to adjust this definition of a similarity function~$\sigma$ to a concrete situation. One natural adjustment can be made if the potential features $s\in S$ have default specifications~$\vs$ that are more distinguishing than their inverses~$\sv$. In our furniture example, $\vs$~might denote the feature of being collapsible, which is a lot more distinguishing than not being collapsible. In this case we might want to count only these distinguishing features in the definition of~$\sigma$, adjusting it to
  $$\sigma(u,v) := \big|\{\,s\in S : u(s) = v(s) = \vs\,\}\big|\,,$$
 reflecting the fact that we consider two collapsible pieces of furniture as more like one another than two arbitrary non-collapsible pieces of furniture.

Based on such a similarity measure on~$V\!$, we can now define an order function on the partitions of~$V\!$ that assigns low order to partitions that do not split many pairs of similar objects:
 $$|A,B| := \sum_{a\in A}\sum_{b\in B} \sigma (a,b)\,.$$
Note that, as desired, $|A,B|$ is large if there are {\it many\/} pairs $(a,b)$ with {\it high\/} similarity~$\sigma(a,b)$, and small otherwise. It is not hard to show that this order function is submodular.%
   \Footnote{Here is a sketch of the easy proof. Given $\vr=A$ and~$\vs=C$, the sum $\sigma(A)$ of all $\sigma(u,v)$ with $u\in A$ and $v\in \Abar$ is~$|r|$, and the analogously defined sum $\sigma(C)$ is~$|s|$. Now computing these sums for $\vr\lor\vs = A\cup C$ and~$\vr\land\vs = A\cap C$, we notice that the terms $\sigma(u,v)$ that count towards $\sigma(A\cup C) + \sigma (A\cap C)$ also count towards~$\sigma(A) + \sigma(C)$, and equally often. Thus, $|\vr\lor\vs| + |\vr\land\vs| = \sigma(A\cup C) + \sigma (A\cap C) \le \sigma(A) + \sigma(C) = |r|+|s|$ as desired.}

Although order functions as above can be defined on all partitions of~$V\!$, they can equally well be used also when we quite deliberately restrict our attention to the original partitions of~$V\!$ from~$S$~-- for example, because they come with `meaning' which we still need in order to interpret our tangles once they are found.%
   \Footnote{Tangle-based teaching methods, discussed in Chapters~\secInformalExamples.4 and~\secFormalExamples.4, are a good example.}

If desired, we can even tweak this method to mimick explicit definitions of order functions on~$S$ in a way that circumvents the problem with submodularity discussed earlier. We wish to assign low order to important features. Instead of doing this explicitly, we can give important features greater weight in the definition of our similarity function~$\sigma$.

To do this formally, we choose a {\it weight function\/} $w\colon S\to\N$ that assigns large values to important potential features. Our first definition of~$\sigma$ then turns into
 $$\sigma_w(u,v) := \sum\{\,w(s): s\in S, \ u(s) = v(s)\,\}\,,$$
 while our second definition turns into
 \Nobreak$$\sigma_w(u,v) := \sum\{\,w(s): s\in S, \ u(s) = v(s) = \vs\,\}\,.$$
We leave the definition of $|A,B|$ as $\sum_{a\in A}\sum_{b\in B} \sigma_w (a,b)$ unchanged.

Note that increasing the weight of a potential feature $s=\{A,B\}$ does not affect the value of $|A,B|$ itself, since $w(s)$ does not occur in the sum $\sigma_w(a,b)$ for any $a\in A$ and~$b\in B$, because $a(s)\ne b(s)$ for such $a$ and~$b$. But raising the weight of~$s$ increases the order of other $r\in S$: of every $r\in S\sm \{s,\{\es,V\}\}$ in the case of the first definition of~$\sigma_w$, and of every $r\in S$ that splits the set $\vs\sub V\!$ in the case of the second definition of~$\sigma_w$.%
   \Footnote{Indeed, let $r=\{C,D\}$. As $s\ne r$, there exist $c\in C$ and $d\in D$ that both lie in~$A$ or both lie in~$B$.%
   \COMMENT{}
   Hence $c(s)=d(s)$, and so $w(s)$ occurs in the sum of~$\sigma_w(c,d)$ and therefore in the sum of~$|C,D|$.}
   Hence if we increase~$w(s)$, the weight of~$s$, the order of~$r$ will rise: a little bit if $r$ and~$s$ are similar as partitions of~$V\!$, and a lot if the they are very different.%
   \COMMENT{}
   Thus, indirectly, assigning $s$ large weight amounts to assigning it, and elements of~$S$ similar to it, low order relative to the order of the other elements of~$S$.
   
However, there is an important difference to choosing orders explicitly: order functions based on a weight function $s\mapsto w(s)$ in this way are defined on all partitions $\{A,B\}$ of~$V\!$, not just on $S$ itself. In particular, they are defined on all Boolean expressions of features in~$S$.%
   \Footnote{We shall need these for our tangle algorithms later.}
   Moreover, these order functions are submodular, with the same easy proof as in the unweighted case.

Finally, tangles of sets~$S_k$ defined with respect to such an order function are usually%
   \COMMENT{}
   witnessed by a weight function on~$V\!$ (see Chapter~\secIdea.5). Indeed, since our order function is submodular, we can extend~$S_k$ to a submodular set $S'_k$ of partitions of~$V\!$ which still all have order~$<k$.%
   \COMMENT{}
   Then all tangles of~$S'_k$ are witnessed by a weight function on~$V\!$, see~\cite{Deciders_k-conn}. Hence any tangle of~$S_k$ that extends to a tangle of~$S'_k$ will also be witnessed by the weight function on~$V\!$ for that tangle of~$S'_k$, since $S_k\sub S'_k$, and most tangles of a set of partitions of~$V\!$ do extend to tangles of the submodular extensions of that set.%
   \COMMENT{}

We thus have an indirect way of assigning the potential features in~$S$ low or high order, explicitly if indirectly, which avoids the technical problems this would cause if we did it directly. When we later say in our discussion of applications that we might choose the order of some features explicitly in one way or another, then doing it in this indirect way will often be what is meant.

Tangles based on the above type of order function have been used, for example, to classify texts by topic. Even if $S$ is a set of very simple potential features of texts describing such basic properties as the frequency of words, one can set the parameters so that the tangles found correspond to the topics that these texts are about. See~\cite{TextTangles} for more.

\subsection \secFormalSetup.5\ \ Duality of set separations

In Chapter~\secInformalExamples.6 we saw an intriguing example of a pair $V,S$ of sets each of which could be regarded as a set of potential features of the other. Each of them had its own tangles, and these tangles described different aspects of the relationship between $V\!$ and~$S$. Now that we have defined tangles and their ingredients more formally, we can see easily that such a `duality' between $V\!$ and~$S$ exists always.

This means that whenever we encounter a natural setting for studying tangles, there is always also a dual setting that may have tangles. These may shed additional light on the situation we are studying, and may be worth exploring even when this was not initially envisaged.

In this spirit, this section is intended to throw an additional light on the applications of tangles we describe in this paper (all of them), but will not be formally needed for the rest of this paper except in Chapters \secThms.4\ and~\secFormalExamples.2.

Consider two disjoint sets, $X$ and~$Y\!$. Assume that some elements $x$ of~$X$ are in a `special relationship' with some elements $y$ of~$Y\!$. We can express this formally by way of a set $E$ whose elements are sets $\{x,y\}$ with $x\in X$ and ${y\in Y}$:\penalty-200\ not all of them, but some.%
   \Footnote{For mathematicians: $E$~is the edge set of a bipartite graph with vertex classes $X$ and~$Y\!$.}

Our choice of~$E$ defines for every $x\in X$ a subset $\vx$ of~$Y$ consisting of just those~$y$ with whom $x$ has this special relationship. Similarly, for every $y\in Y$ we have the set~$\vy$ of all $x\in X$ with whom $y$ has this special relationship:
  $$\vx := \{\,y\in Y\mid \{x,y\}\in E\,\}\quad\hbox{and}\quad\vy := \{\,x\in X\mid \{x,y\}\in E\,\}.$$
 Then for all $x\in X$ and $y\in Y$ we have
 $$x\in\vy\ \Leftrightarrow\ y\in\vx\eqno (*)$$
as both these are true if $\{x,y\}\in E$ and false otherwise.

Let us define $\xv$ as the complement of the set~$\vx$ in~$Y\!$, and $\yv$ as the complement of the set~$\vy$ in~$X$. Further, let
 $$\eqalign{\vX &:= \{\,\vx\mid x\in X\,\}\cup \{\,\xv\mid x\in X\,\}\cr
    \vY &:= \{\,\vy\mid y\in Y\,\}\cup \{\,\yv\mid y\in Y\,\}.\cr} $$
 Then $\vX$ becomes a set of features of~$Y\!$, and $\vY$ is a set of features of~$X$, for the partitions $\{\vx,\xv\}$ of~$Y$ and $\{\vy,\yv\}$ of~$X$ as potential features. We say that these feature sets are {\it dual\/} to each other.%
   \Footnote{For mathematicians: this is well defined, stated symmetrically as above, since everything we have said has been symmetrical in $X$ and~$Y\!$.}%
   \COMMENT{}

Note that every $x$ determines the set~$\{\vx,\xv\}$ and vice versa, and similarly every~$y$ determines the set~$\{\vy,\yv\}$ and vice versa. We may therefore choose to ignore the formal difference between $x$ and~$\{\vx,\xv\}$, and between $y$ and~$\{\vy,\yv\}$. Then $X$ itself becomes a set of potential features of (elements of)~$Y\!$, and $Y$ a set of potential features of (elements of)~$X$.

Let us now see how in general, given any set $S$ of partitions of a set~$V\!$, we can define a set $\vV\!$ of features of~$S$ that is dual in this sense to the set~$\vS$ of features of~$V\!$.%
   \COMMENT{}
   We start by picking for every $s\in S$ a default specification, which we denote as~$\vs$ (rather than~$\sv$). This determines for every $v\in V\!$ the sets
 $$\vv = \{\,s\in S\mid v\in\vs\,\}\quad\hbox{and}\quad\vvback = \{\,s\in S\mid v\in\sv\,\},$$
 which form a partition of~$S$. If, as we shall assume, these partitions $\{\vv,\vvback\}$ differ for distinct $v\in V\!$, they determine their $v$ uniquely and we may think of each $v$ as shorthand for~$\{\vv,\vvback\}$. This%
   \COMMENT{}
   makes $V\!$ into a set of partitions $v = \{\vv,\vvback\}$ of~$S$ and
 $$\vV := \{\,\vv\mid v\in V\}\>\cup\> \{\,\vvback\mid v\in V\}$$
 into the set of all specifications of elements of~$V\!$. Thus $\vV$ and~$\vS$ form an instance of a pair $\vX,\vY$ of dual feature sets as defined earlier.%
   \Footnote{The `special relationship' exists between elements $v\in V$ and $s\in S$ for which $v\in\vs$ and, equivalently, $s\in\vv$.}

With this setup in place, we can now investigate how the tangles of~$S$%
   \COMMENT{}
   are related to those of~$V\!$, how they or their distinguishers~$T$ compare with subsets of~$S$ that witness tangles of~$V\!$, and so on. But this is not the place to do this; readers interesting in these connections may find some more inspiration in~\cite{TangleDualityHomology}.%
   \COMMENT{}

However we have seen enough now of this duality to apply it in practice: whenever we have a set $V\!$ of `objects' at hand, and a set $\vS$ of features of these objects, then together with investigating the tangles of~$S$ we can also investigate the tangles of~$V\!$ viewed as a set of potential features of the elements of~$S$. Both together are likely to paint a more comprehensive picture of the situation we are trying to analyse than just working with one of the two types of tangle.

\beginsection \secThms. Tangle theorems and algorithms

In this chapter we present the two most fundamental theorems about tangles, the {\it tree-of-tangles theorem\/} and the {\it tangle-tree duality theorem\/}, in their simplest forms. Yet the versions presented here are strong enough to apply to all the application scenarios we have discussed as examples, and are likely to apply equally to most other application scenarios that arise in the empirical sciences.

$\!\!$The tree-of-tangles theorem, of which we present two versions in Section~\secThms.1, finds in a submodular set $S$ of potential features a small subset~$T$ that suffices to distinguish all the tangles of~$S$: for every two tangles of~$S$ there is a potential feature in~$T$, not only in~$S$, which these two tangles specify differently.

In Section~\secThms.2 we present the tangle-tree duality theorem. If $S$ has no tangle, this theorem finds a small subset $T$ of~$S$ that suffices to demonstrate this fact conclusively. We say that $T$ demonstrably `precludes' the existence of a tangle. This can help to prove that some given data is, in a precise and quantitative sense, unstructured or contaminated.

In Section~\secThms.3 we look at some generic algorithms that compute tangles, tangle-distinguishing features sets, and tangle-precluding feature sets.

In Section~\secThms.4 we address the question of how the elements of~$V\!$ can be grouped, or clustered, `around' the tangles of~$S$.

In Section~\secThms.5, finally, we discuss why the tangle-distinguishing property of the set $T$ found by the two theorems of Section~\secThms.1 makes it particularly valuable: we can predict from how a given $v\in V\!$ specifies the potential features in~$T$ how it will probably specify most of the other $s\in S$.

\medbreak

The following definitions will be needed throughout this chapter.
A {\it feature system\/} is a set $\vS$ obtained from a set $S$ of partitions of a set~$V$: it consists of all the specifications $\vs$ and~$\sv$ of elements $s$ of~$S$. A~subset $T$ of~$S$ is a {\it tree set\/}%
   \Footnote{See~\cite{TreeSets} for what the structure of such sets~$T$ has to do with trees.}
   if $\{\es,V\}\notin T$ and for every two partitions in~$T$ one side of each partition is contained in a side of the other partition.%
   \Footnote{Thus, for all $\{A,\Abar\}, \{B,\Bbar\}\in T$ we have $A\sub B$ or $A\sub\Bbar$ or $\Abar\sub B$ or $\Abar\sub\Bbar$, and correspondingly $\Bbar\sub\Abar$ or $B\sub\Abar$ or $\Bbar\sub A$ or $B\sub A$. Such sets $A,B$ are also called {\it nested\/}.}

While $S$ can be exponentially large in terms of~$|V|$, tree sets in~$S$ are much smaller, at most linearly large in terms of~$|V|$.%
   \COMMENT{}
   The tangle-distinguishing tree sets $T$ we find in Section~\secThms.1 will be much smaller still, no larger than the number of tangles.

\subsection \secThms.1\ \ Tangle-distinguishing feature sets

Consider two different tangles of~$S$. As both are specifications of~$S$, the fact that they are different means that there exists an $s\in S$ which one of our two tangles specifies as~$\vs$, and the other as~$\sv$. We say that $s$ {\it distinguishes\/} the two tangles. More generally, we say that a subset $T$ of~$S$ {\it distinguishes all the tangles\/} of~$S$ if for every two tangles of $S$ we can find in this set~$T$ an $s$ that distinguishes them.

Our first theorem says that, if $S$ is submodular,%
   \Footnote{If it is not, we may wish to make it submodular by adding some elements obtained as Boolean expressions of existing elements, as outlined in Chapter~\secFormalSetup.3.}
   we can find such a tangle-distinguishing subset $T$ of~$S$ that is a tree set:

\proclaimwithname Theorem \xxxToTfixedS. (Tree-of-tangles theorem for fixed~$S$~\cite{AbstractTangles})
 For every submodular feature system~${\vS}$ there is a tree set $T\sub S$ that distinguishes all the tangles of~$S$.

Viewed in terms of~$V\!$, the elements of~$T$ are nested partitions of~$V\!$. If $T$ is chosen minimal, these partitions of~$V\!$ cut it up into disjoint parts%
   \COMMENT{}
   which correspond bijectively to the tangles of~$S$. Indeed, every tangle $\tau$ of~$S$ induces a tangle of~$T$, the tangle $\tau\cap\vT$, and the sets
 $$V_\tau := \bigcap \,\{\, A\mid A\in \tau\cap\vT\,\}$$%
   \COMMENT{}
   partition~$V\!$ as $\tau$ varies over the tangles of~$S$.%
   \COMMENT{} Each of the non-empty parts~$V_\tau$ into which $T$ cuts the set~$V\!$ is thus `home' to exactly one%
   \COMMENT{}
   tangle~$\tau$, but note that $V_\tau$ can be empty for some~$\tau$.%
   \COMMENT{}

   An important theoretical consequence of Theorem~\xxxToTfixedS\ is that $S$ has only few tangles, even of the most general kind (with $\F=\es$).%
   \COMMENT{}
   Indeed, as $T$ distinguishes all the tangles of~$S$, these induce distinct tangles of~$T$.%
   \COMMENT{}
   Now a tree set of size~$\ell$ has at most $\ell+1$ tangles,%
   \Footnote{For graph theoristis: it has exactly $\ell+1$ consistent specifications, all of which are tangles of the most general kind. This is because the consistent specifications of a tree set~$T$ correspond to the nodes of a tree with $|T|$~edges~\cite{TreeSets}, and a tree with $\ell$~edges has exactly $\ell+1$ nodes.\looseness=-1}
   and all tree sets in~$S$ are small.%
   \Footnote{They have size less than~$2|V|$, independently of the size of~$S$.%
   \COMMENT{}
   So the number of all tangles of~$S$ is at most~$|T|+1\le 2|V|$.%
   \COMMENT{}
   By contrast, $S$~itself can have size exponential in~$|V|$, and the number of specifications of~$S$ is exponential even in~$|S|$.}

In practice, the implication also goes the other way. We can control directly about how many tangles we get by specifying~$\F$, e.g., by choosing the parameter~$n$ for $\F=\F_n$. Since~$T$, if chosen minimal, is automatically smaller than the number of tangles of~$S$,%
   \COMMENT{}
   it will then be small as a consequence.

This aspect is particularly valuable if $S$ is not submodular (as Theorem~\xxxToTfixedS\ requires). To distinguish $\ell$~tangles we only need $\ell-1$ elements of~$S$, and these are quick to find.%
   \COMMENT{}
   However, they may not form a tree set. Our algorithm will then build from these elements of~$S$ a tree set~$T$, also of size~${\ell-1}$, which once more distinguishes all our $\ell$ tangles of~$S$. This $T$ will not be a subset of~$S$. But the elements of~$\vT$ will be Boolean expressions of specifications of those $\ell-1$ elements of~$S$ we started with,%
   \COMMENT{}
   obtained by iteratively adding features such as $\vr\lor\vs$ or ${\vr\land\vs}$ combined from existing or previously added features $\vr$ and~$\vs$.%
   \COMMENT{}
   So this $T$ will be small whichever way we choose to measure it: in absolute terms, as~$|T|$, or in terms of how many features from~$\vS$ are needed to build it.

There is no general rule of how many tangles to aim for when choosing~$\F_n$: this will depend on the application.  Choosing $n$ large will give us only the most pronounced tangles, while choosing it lower will give us more refined tangles.%
  \Footnote{This a different kind of refinement than in the hierarchy of tangles discussed earlier.} 
But if $S$ satisfies a little more than submodularity, we can dramatically improve Theorem~\xxxToTfixedS\ to find a tree set in~$S$ that distinguishes even tangles of varying degrees of refinement, all at once. This will be our next theorem.

Assume now that $\vS$ is not just a submodular feature system but even a universe of features of~$V\!$ with a submodular order function $s\mapsto |s|$.%
   \COMMENT{}
   Any tangle of some $S_k = \{\,s\in S: |s| < k\,\}$ will be called a tangle~{\it in\/}~$S$. Two tangles in~$S$, not necessarily of the same order, are {\it distinguishable\/} if there exists an $s\in S$ which they specify differently.%
   \Footnote{Such an $s$ will always exist for distinct tangles of the same order, i.e., for distinct tangles of the same set~$S_k$. But a tangle of order~$k$ is not distinguishable from the tangles of order $\ell<k$ that it induces on the sets~$S_\ell$, as defined in Section~\secFormalSetup.2. Our assumption that $s$ distinguishes two tangles, of orders $k$ and~$\ell$, say, thus means in particular that $|s|<\min\{k,\ell\}$, and it implies that the tangle of higher order does not induce the tangle of lower order.}
   Such an $s$ {\it distinguishes\/} the two tangles, and it does so {\it efficiently\/} if no $r\in S$ of lower order distinguishes them. A~set $T\sub S$ {\it distinguishes all the tangles in~$S$ efficiently\/} if for every two distinguishable tangles in~$S$ there is an $s\in T$ that distinguishes them efficiently.

\proclaimwithname Theorem \xxxToTallSk. (Tree-of-tangles theorem for tangles of variable order~\cite{ProfilesNew})
   For every feature universe $\vS$ with a submodular order function there is a tree set $T\sub S$ that distinguishes all the tangles in~$S$ efficiently.%
   \COMMENT{}

\noindent
   See~\cite{ProfilesNew,FiniteSplinters} for more general versions of Theorem~\xxxToTallSk.

\medbreak

As earlier, every tangle $\tau$ in~$S$ induces a tangle~$\tau\cap\vT$.%
   \COMMENT{}
   If $T$ is chosen min\-imal then the sets
 $$V_\tau := \bigcap \,\{\, A\mid A\in \tau\cap\vT\,\}$$
 partition~$V\!$ as $\tau$ varies over the maximal%
   \Footnote{Maximal in terms of inclusion as subsets of~$\vS$. The maximal tangles in~$S$ are those not induced by any other tangle of larger order.}
    tangles in~$S$.

The relative structure between the tree set $T\sub S$ in Theorem \xxxToTfixedS\ or~\xxxToTallSk\ and the tangles it distinguishes can be made visible by a graph-theoretic tree~$\TT$ whose edges correspond to the elements of~$T$ and whose nodes correspond to the (maximal) tangles $\tau$ of or in~$S$.%
   \Footnote{If $V_\tau\ne\es$,%
   \COMMENT{}
   then $\tau$ orients all the edges of~$\TT$ towards `its' node by specifying their corresponding $t=\{A,B\}\in T$ as the $A$ for which $V_\tau\sub A$.}
   In the case of Theorem~\xxxToTallSk, the tree $\TT$~also displays the non-maximal tangles in~$S$ and their relative structure (such as which refines which): not by its nodes, but by some of its subtrees.%
   \Footnote{Given~$k$, the tangles of~$S_k$ correspond to the nodes of the tree~$\TT_k$ for the tree set~$T\cap S_k$, which in turn correspond to the components of the forest obtained from~$\TT$ by deleting the edges that correspond to elements of~$T$ in~$S_k$. These subtrees are nested in a way that reflects the nestedness of the tangles in~$S$ as subsets. See~\cite{TreeSets,AbstractTangles,ProfilesNew}.}

\subsection \secThms.2\ \ Tangle-precluding feature sets

The tree-of-tangles theorems from Section~\secThms.1 offer a geometric picture of how the various (maximal) tangles of~$S$ are organized in a tree-like way. Viewed in terms of the set~$V\!$, they produce a tree set~$T$ of partitions which cut $V\!$ up into disjoint parts that correspond exactly to the (maximal) tangles of~$S$.

Our next theorem, the tangle-tree duality theorem, produces a similar tree set $T\sub S$ for the case that $S$ has no tangles~-- understood to be $\F_n$-tangles for some $n$ to be specified, as before.%
   \Footnote{See Chapter~\secFormalSetup.1 for the definition of $\F_n$-tangles.}
   Once more, $T$~partitions~$V\!$ into disjoint parts. This time these parts are not `home to a tangle', but the opposite: they are so small that they {\it cannot\/} be `home to a tangle'. Indeed if `tangle' means $\F_n$-tangle, they have size~$<n$. Let us make this more precise.

Consider a consistent specification $\tau$ of a tree set~$T\sub S$. Let $\sigma\sub\tau$ be the set of those elements $\vt$ of~$\tau$ for which $\vt=C$ is minimal under inclusion. One can show that~$\sigma$ is a star of features (see Chapter~\secFormalSetup.3). These sets~$\sigma$, one for every consistent specification of~$T$, are the {\it splitting stars\/} of~$T$.%
   \Footnote{As remarked earlier, the `tree-shape' of a tree set can be made visible by a graph-theoretic tree whose edges correspond to the elements of~$T$ and whose nodes correspond to the consistent specifications of~$T$. The splitting stars of~$T$ correspond to the stars at these nodes~\cite{TreeSets}.}%
   \COMMENT{}
  A~tree set $T\sub S$ is {\it over\/}~$\F$ if all its splitting stars are elements of~$\F$.

\proclaimwithname Theorem \xxxTTD.
 (Tangle-tree duality theorem~\cite{AbstractTangles})
  Let ${\vS}$ be a submodular feature system, and let $\F=\F_n$ for some integer~${n\ge1}$.%
  \COMMENT{} 
  Then exactly one of the following two statements holds:
     			\pitem{i} $S$ has an $\F$-tangle.%
   \COMMENT{}
     			\pitem{ii} $S$ contains a tree set over~$\F$.%
   \COMMENT{}

\noindent
   See~\cite{AbstractTangles,TangleTreeGraphsMatroids,TangleTreeAbstract} for more general versions of Theorem~\xxxTTD.

\penalty-500\medbreak

It is easy to see why (i) and (ii) cannot happen at once, why the existence of a tree set over~$\F$ precludes the existence of an $\F$-tangle.%
   \Footnote{Indeed, suppose there are both a tangle~$\tau$ as in~(i) and a tree set $T$ as in~(ii). As ${T\sub S}$, the tangle $\tau$ induces a consistent specification (indeed a tangle) of~$T$. As $T$ is over~$\F$, the splitting star of this specification of~$T$ is an element of~$\F$. Since it is also a subset of~$\tau$, this contradicts the fact that $\tau$ is an $\F$-tangle.}
   The harder part of the theorem is that at least one of (i) and~(ii) must always happen: that if $S$ has no $\F$-tangle then it contains a tree set as in~(ii), one which demonstrably precludes the existence of an $\F$-tangle (as just noted).

This is valuable, and it is the way in which Theorem~\xxxTTD\ is usually applied. For example, suppose that some tangle-finding algorithm fails to find a tangle. This does not mean, yet, that no tangle exists. But if the algorithm also finds a tree set as in~(ii) of Theorem~\xxxTTD, we can be sure that indeed no tangle exists.%
   \Footnote{It is easy to check that a given tree set that is claimed to be over~$\F$ is indeed over~$\F$: we just have to compute its splitting stars, which is easy, and check that they are all in~$\F$.}
   The theorem thus assures us that whenever no tangle exists there is conclusive and easily checkable evidence for this.

\subsection \secThms.3\ \ Algorithms

There are a number of algorithmic aspects concerning tangles. These include the computation of the following objects:%
   \COMMENT{}

\medskip
\plainitem{$\bullet$} the set $S$, or sets~$S_k$, in whose tangles we are interested;
\plainitem{$\bullet$} the tangles of~$S$ or of the~$S_k$;
\plainitem{$\bullet$} if tangles are found, a set tree $T$ of potential features distinguishing them (as in Theorems~\xxxToTfixedS\ or~\xxxToTallSk, respectively);
\plainitem{$\bullet$} if no tangles are found, a tree set $T$ of potential features that certifies their non-existence (as in Theorem~\xxxTTD).

\medbreak

In applications where $S$ is given directly by a questionnaire, the first task does not arise.%
   \COMMENT{}
   It may happen that, in the course of computing the tangle-distinguishing feature set $T$ for the third task, the set $\vS$ is enlarged by newly added features $\vr\land\vs$ or $\vr\lor\vs$ for already existing features $\vr$ and~$\vs$, but any such enlargement of~$S$ will happen only then. In particular, we do not need to generate a feature universe from~$\vS$ before we look for or process its tangles, or even to extend $\vS$ to a submodular feature system to which our tangle theorems can be formally applied.

In other applications,%
   \Footnote{including those where we start from a questionnaire but do not take this to be our set~$S$}
   where we consider as $S$ some submodular extension of a questionnaire, or where we take $S$ to be the entire set of partitions of a set~$V\!$ but only look for tangles of its subsets $S_1\sub S_2\dots$ as in Chapter~\secFormalSetup.4, we do have to compute $S$ or these sets~$S_k$.

When we extend a given feature system~$\vSdash$ that is not submodular to a submodular system $\vS\supe \vSdash$ to which we wish to apply our theorems (which assume submodularity for~$S$ or the order function that determines the~$S_k$), it can happen in rare cases that some tangles of~$S'$, those we are interested in, fail to extend to tangles of~$S$. In such cases the tangle-distinguishing set~$T$ found for~$S$ in the third task will distinguish most of the tangles of~$S'$, those that extend to tangles of~$S$, but maybe not all of them.

In the case of the fourth task the tangle-tree duality theorem applied to~$S$ might, theoretically, produce a tree set~$T$ certifying that $S$ has no tangle even when $S'$ did have tangles but these did not extend to tangles of~$S$. In order to avoid this, we always try to compute the tangles of the feature system~$S'$ we are interested in directly, even if $S'$ is not submodular. If no tangle is found, we then run the algorithm for the fourth task on~$S'$. This will usually produce a tree set $T\sub S'$ certifying that $S'$ has no tangle. If it does not, the algorithm will output a list of features that must be included in any tangle of~$S'$. If $S'$ is submodular then such a tangle will in fact exist. But even if $S'$ is not submodular, as in our hypothetial case here, this list will be a good start to restart our search for a tangle of~$S'$, which in most cases is also found.

Although the~$S_k$ may not be large for the small values of~$k$ relevant to us, determining them precisely can take longer than is practicable. For example, determining the order of every partition of a (large) set~$V\!$ takes time exponential in~$|V|$, since there are $|V|\choose 2$ such partitions. But often there are known heuristics for the application in question that can reduce the set of rele\-vant partitions, or the low-order partitions can be computed by known methods.%
   \Footnote{For example, if our application is to capture fuzzy but discernable regions in a picture, we might run a standard clustering algorithm to find a much smaller set~$S$ of partitions of the pixel set that correspond to natural lines%
   \COMMENT{}
   and do not cut right through a big cluster. If the order of our partitions of~$V\!$ is determined by a similarity measure~$\sigma$ or~$\sigma_w$, as in Chapter~\secFormalSetup.4, we can use network flow theory to find the partitions of low order quickly.}%
   \COMMENT{}
   If that is not the case, one can use standard approximation methods, such as finding local minima close to random partitions,%
   \COMMENT{}
   to find not all but enough features of order~$<k$ for each rele\-vant~$k$ to compute large enough subsets of the~$S_k$ to make their tangles significant.

Once we have decided for which $S$ or $S_k$ we wish to compute its $\F$-tangles, we can do this fairly quickly in an ad-hoc way: we go through its elements~$s$ one by one and add $\vs$ or~$\sv$ to every tangle previously found for the subset of~$S$ already considered, as long as this does not create a forbidden subset. But checking whether a given subset of~$\vS$ is forbidden, i.e., in~$\F$, is usually quick.%
   \Footnote{For example, for $\F=\F_n$ as in Chapter~\secFormalSetup.1\ we just have to count for every two elements $A$ and~$B$ in a partial tangle already computed how many elements of~$V\!$ lie in $A\cap B\cap C$ or in $A\cap B\cap \Cbar$, where $s=\{C,\Cbar\}$ is the new element of~$S$ under consideration.}

Computing a tangle-distinguishing feature set $T$ in our third task is more subtle. Grohe and Schweitzer~\cite{GroheSchweitzerTangleAlg} have turned the proof of Theorem~\xxxToTallSk\ given in~\cite{ProfilesNew} into an algorithm with polynomial running time, which however is not intended for practical applications.%
   \Footnote{However it is of theoretical significance for computational isomorphism problems. This takes us too far afield,%
   \COMMENT{}
   but any reader interested is encouraged to consult~\cite{ProfilesNew} and~\cite{GroheSchweitzerTangleAlg}.}
   Elbracht, Kneip and Teegen~\cite{TangleAlgorithms} have described a different algorithm, also intricate, which has an exponential worst-case running time but has shown excellent performance in practice.%
   \Footnote{Typical running times for a few dozen tangles of a few thousand features measure in minutes or seconds on an ordinary laptop.}%
   \COMMENT{}
   Even if $S$ is not submodular, it uses only $\ell-1$ elements of~$S$ in total%
   \COMMENT{}
   to build a tree set~$T$ of size $\ell-1$ to distinguish $\ell$ tangles.%
   \COMMENT{}

If we do not require our tangle-distinguisher~$T$ to be a tree set but just want it to be small, it can be found inside $S$%
   \COMMENT{}
   and still of size~$|T|<\ell$.%
   \COMMENT{}
   The fact that $T$ is a tree set, however, is crucial to a number of tangle applications. One consequence is that its tangles extend (uniquely) to tangles of~$S$.%
   \COMMENT{}
   Another is that it enables us to associate elements of~$V\!$ with tangles of~$S$; see Section~\secThms.4 below. Both these properties together form the basis for tangle-based predictions, discussed in Section~\secThms.5.

\goodbreak

Computing a tangle-precluding feature set~$T$ as in our fourth task may or may not be fast,%
   \Footnote{Our current algorithm has a worst-case runnig time of $O(|S|^4)$ for submodular~$S$.}%
   \COMMENT{}
   depending on the application; see~\cite{TangleAlgorithms} for details.

\subsection \secThms.4\ \ Tangle-based clustering

As we have seen, tangles of a set $S$ of features of a set~$V\!$ of objects, including tangles of views prevalent in a population of people, are best at discovering what {\it kinds\/} of groupings exist in~$V\!$, rather than at dividing~$V\!$ into such groups. However, once found, tangles of~$S$ can help with this too.

Before we explore this further, let us briefly look at two other ways of defining clusters in~$V\!$ that are also determined by tangles, but which we shall not pursue further here. The first of these are witnessing sets of tangles, as mentioned briefly in Chapter~\secIdea.5. These may indeed be interesting in practice, and they helped motivate tangles when they were invented. The question of when tangles have witnessing sets has not been studied sufficiently so far. But most of the tangles we have studied here do have witnessing functions~\cite{Deciders_k-conn}.

Another approach we shall not pursue in detail is duality: to obtain clusters in~$V\!$ from the tangles of~$V\!$ (rather than of~$S$) that arise when we interpret~$V\!$ as a set of potential features of the elements of~$S$, as explained in Chapter~\secFormalSetup.5. Informally, just as tangles of~$S$ are collections of features shared by many of the objects in~$V\!$, tangles of~$V\!$ are sets of points sharing many of the features in~$\vS$, which on the face of it looks exactly what we would expect of a cluster.

More formally, a specification~$\sigma$ of~$V\!$ is an $\F_n$-tangle of~$V\!$ if for all triples $\vx,\vy,\vz\in\sigma$%
   \COMMENT{}
   we have $|\vx\cap\vy\cap\vz|\ge n$. By definition of~$\vV$ (see Chapter~\secFormalSetup.5) this means that for all triples $v_1,v_2,v_3\in V\!$ there are $n$ elements $s$ of~$S$ such that ${v_i\in\!\vs}$ if $\sigma(v_i) = \vvi$ and $v_i\in\sv$ if $\sigma(v_i) = \vviback$, where $\vs$ is some fixed default specification of~$s$ chosen in the process of definition the feature system~$\vV\!$.

Given such an $\F_n$-tangle~$\sigma$ of~$V\!$, we can associate with it the set
 $$V^+_\sigma := \{\,v\in V\mid \sigma(v) = \vv\,\}$$
 of those $v\in V\!$ which $\sigma$ specifies as $\vv$ rather than as~$\vvback$.%
   \COMMENT{}
   This is a cluster in~$V\!$ in the sense described informally above: for every three points $x,y,z\in V_\sigma$ there are at least $n$ elements $s$ of~$S$ such that $x,y,z\in\vs$, i.e., such that $x$, $y$ and~$z$ share these $n$ features $\vs\in\vS$. Similarly,
 $$V^-_\sigma := \{\,v\in V\mid \sigma(v) = \vvback\,\}$$
 is a set of points sharing many of the features in~$\vS$: every three of them share at least $n$ of the features of the form~$\sv$, or equivalently, jointly lack at least $n$ features of the form~$\vs$ (in our default specification of~$S$). In fact, neither of these sets captures the full information encoded in the tangle~$\sigma$ of~$V\!$.%
   \COMMENT{}

More importantly, both clusters, $V^+_\sigma$ and~$V^-_\sigma$, depend not only on $V\!$ and~$S$ but also on the default specification of~$S$ we chose initially.%
   \COMMENT{}
   These kinds of clusters in~$V\!$, therefore, are significant only if this default specification of~$S$ is particularly natural~-- as is the case in our example of shopping tangles in Chapter~\secInformalExamples.6,%
   \Footnote{There, specifying an item $s$ as purchased rather than as not purchased is clearly the natural choice.}
   but cannot be expected in general.

Let us get back to our original topic of how to define clusters in~$V\!$ associated with tangles of~$S$, a~set of potential features of~$V\!$. A~simple way to do this is to make one group for each tangle of~$S$ and associate each $v\in V\!$ with the tangle closest to it under the first similarity function~$\sigma$ defined in Chapter~\secFormalSetup.4: the tangle~$\tau$ of~$S$ for which the number of elements $s\in S$ with $\tau(s) = v(s)$ is maximum.

However, in Chapter~\secFormalSetup.4 we also saw how this similarity function can be used to define an order function on~$S$,%
   \COMMENT{}
   indeed on the set of {\it all\/} partitions of~$V\!$, to obtain a hierarchy $S_1\sub S_2\sub\dots\sub S$ of sets of partitions of~$V\!$ and their associated tangles.%
   \Footnote{Recall that we use the term `tangle {\it in\/}~$S$' for all tangles of the sets~$S_k\sub S$, and that the tangles of~$S_k$ are said to be tangles (in~$S$) {\it of order~$k$.}}
   Using this hierarchy, we can do something more sophisticated to associate the elements of~$V\!$ with tangles in~$S$.

For example, we could make a group in~$V\!$ for every tangle in~$S$, and associate a given $v\in V\!$ with the tangle $\tau$ of largest order whose features are among those of~$v$, in that $\tau(s) = v(s)$ for all $s\in S_k$ if $\tau$ has order~$k$.%
   \COMMENT{}
   Thus, we would first look which tangle~$\tau_1$ of~$S_1$, if any, is a subset of $v(S) = \{\,v(s)\mid s\in S\,\}$; then which tangle~$\tau_2$ of~$S_2$ refining~$\tau_1$ is still a subset of~$v(S)$; and so on until we have found a tangle $\tau_k\sub v(S)$ such that no tangle of~$S_{k+1}$ refining~$\tau_k$ is a subset of~$v(S)$. Then the views of~$v$ differ from any tangle of~$S_{k+1}$ on at least one~$s\in S$ not specified by~$\tau_k$, and we associate~$v$ with this~$\tau_k$.%
   \COMMENT{}

Note, however, that this approach is not very error-tolerant: even if some $v$ agrees on all of~$S$ with a tangle~$\tau$ of highest order, except in one~$s$ of relatively low order~$k$, then $v$ will not be associated with~$\tau$ but only with the much broader tangle $\tau_k = \tau\cap\,\vSk\,$.

To mend this, we might combine the two approaches above by letting our similarity function take account of all $s\in S$ but attach more importance to those of low order. For example, we could measure the similarity of two specifications $f,g$ of~$S$ by a vector
 $$\sigma(f,g) = \big(\sigma_0(f,g),\sigma_1(f,g),\dots\big)$$
 where $\sigma_k(f,g)$ denotes the number of $s\in S$ of order~$k$ such that $f(s)=g(s)$, and order these vectors lexicographically.%
   \Footnote{For example, an element $v\in V\!$ is more similar in this ordering to a tangle~$\rho$ of~$S$ than to a tangle $\tau$ of~$S$ if $\rho$ and~$\tau$ agree with~$v$ on equally many $s$ of orders 0, 1 and~2, respectively but more $s$ of order~3 satisfy $\rho(s)=v(s)$ than $\tau(s)=v(s)$.}%
   \COMMENT{}

More radically,%
   \COMMENT{}
   we might define the similarity of $f$ and~$g$ as
 $$\sigma(f,g) := \sum \{\,N-|s| : f(s) = g(s)\,\}\,,$$
 where $N$ is the largest order of any $s\in S$.%
   \Footnote{As in our first approach, this counts the $s\in S$ which $f$ and~$g$ specify in the same way, but elements $s$ of low order are given greater weight in this count.}

   Such similarity functions could be used directly for the matchmaking problem discussed in Chapter~\secInformalExamples.1, or as input for a standard distance-based clustering algorithm on~$V\!$.

However, we can also use it to group $V\!$ around the tangles of~$S$, as is our aim here, as long as we say how $\sigma(v,\tau)$ should deal with the fact that if $\tau$ is a tangle of~$S_i$, say, it does not specify the entire~$S$. In our last definition of $\sigma(f,g)$ extended to such partial specifications $f,g$ of~$S$ we might treat $s\in S$ that are not specified by both $f$ and~$g$ like those~$s$ for which $f(s)\ne g(s)$, i.e., ignore them in the sum that computes the similarity between $f$ and~$g$. In our previous definition of~$\sigma(f,g)$, the one with the lexicographic ordering of similarity vectors, we would just set $\sigma_k(v,\tau):=0$ if $\tau$ has order~$i\le k$. We could then make one group in~$V\!$ for every tangle in~$S$,%
   \COMMENT{}
   as before, by associating a given $v\in V\!$ to the tangle most similar to it, the tangle~$\tau$ for which $\sigma(v,\tau)$ is maximum.

Finally, we can use the tree set $T\sub S$ from Theorem \xxxToTfixedS\ or \xxxToTallSk\ to partition~$V\!$ into the `clusters'
 $$V_\tau := \bigcap \,\{\, A\mid A\in \tau\cap\vT\,\}\,,$$
 where $\tau$ varies over the tangles of~$S$ (in the case of Theorem \xxxToTfixedS) or the maximal tangles in~$S$ (in the case of Theorem~\xxxToTallSk). As remarked after these theorems, the sets~$V_\tau$ do indeed form a partition of~$V\!$ if $T$ is chosen minimal, which we may assume.%
   \COMMENT{}

\subsection \secThms.5\ \ Predictions

Let us now discuss how knowing a tangle-distinguishing set~$T$ as in%
   \COMMENT{}
   Section~\secThms.1 enables us to make predictions. For simplicity, we shall do this for the scenario of Theorem~\xxxToTfixedS: we assume that we have some fixed set~$S$ of potential features of a set~$V\!$,%
   \COMMENT{}
   and that we have computed the tangles of~$S$ for some~$\F$ of our choice.%
   \COMMENT{}

To make our discussion more intuitive, let us use the language of the questionnaire example from the introduction. Thus, $S$~is a questionnaire, and we are interested in how an unknown person might answer it. This can be formalized, if desired, by assuming that $V\!$ is a sample of some larger population $P$ of people in whose views on~$S$ we are interested. We know only the answers given by the people in~$V\!$, but have reason to believe that these represent~$P$ well.%
   \Footnote{This is a heuristic assumption that needs validation, but that is a routine matter which is not our concern here. Mathematically, the `unknown person' is simply a random specification of~$S$: a point in our probability space of all the specifications of~$S$. The probability distribution on this space is determined roughly (but not induced exactly)%
   \COMMENT{}
   by how $S$ partitions~$V\!$, in that we require it to assign small measure to specifications of~$S$ that include a forbidden set of features for one of the tangles that we observed were induced on~$S$ by~$V\!$.%
   \COMMENT{}
   Applying this distribution to the abstract space~$\{0,1\}^S$, i.e., `forgetting'~$V\!$, allows us to reinterpret $S$ as a set of `potential features' also of the elements of~$P$, even if these potential features are not known individually as partitions of~$P$.}%
   \COMMENT{}
  
Our contention is that knowing our unknown person's answers to just the questions in~$T$ enables us to predict their answers to the other questions in~$S$ with greater confidence than knowing their answers to an arbitrary subset of~$S$ of the same size as~$T$, and certainly better than choosing their answers ($\vs$~or~$\sv$) to the various questions~$s$ at will.%
   \Footnote{This latter success rate would be~$1/2^{|S|}$ if we assume, for want of any information to the contrary, that all yes/no answers are equally likely and independent of each other for different questions. In mathematical terms, this is the probability of a single specification of~$S$ in the space of all its specifications with the uniform distribution. Our contention is that the existence of a tangle implies that the specifications of~$S$ are not distributed uniformly, and that knowing a tangle tells us enough about the actual distribution to make our prediction better.}%
   \COMMENT{}

Our emphasis will lie on predicting the views of our unknown person correctly with high probability for many of the questions from~$S$: not necessarily for all of them,%
   \COMMENT{}
   nor just for single questions~-- which is easier and does not require tangles.%
  \COMMENT{}
   But note that good predictions for many single~$s$ do not necessarily add up to a good collective prediction for one person's views on many~$s$. For example, the collection of the most popular views on single questions~$s$ may be inconsistent%
   \COMMENT{}
   and thus, in its entirety, not be held by anyone.%
   \Footnote{Here is a simple example in the notation of Capter~\secFormalSetup.1. Consider a tripartition of~$V$ into sets~$A_1,A_2,A_3$ of equal size, and let $S$ consist of the three partitions $s_i=\{A_i,\overline{A_i}\}$ of~$V\!$\penalty-200\ ($i=1,2,3$). For each~$i$, two thirds of the people in~$V\!$ answered~$s_i$ as~$\overline{A_i}$, so this will be our best prediction for how an unknown person would answer~$s_i$. But these three answers together are inconsistent: as $\bigcap_i\overline{A_i} =  \overline{\bigcup_i A_i} = \es$, not a single person in~$V\!$ answered them in this way for all~$i$, and hence this should not (and will not) be our prediction of how an unknown person would answer~$S$.}%
   \COMMENT{}

An aspect important for applications is that although we need to compute the tangles of~$S$ and the tangle-distinguishing set~$T$ in order to make such predictions, we do not need to `understand' the tangles we have computed, to `make sense' of them: the tangle-based process of predicting the answer to~$S$ of our unknown person can be carried out entirely mechanically on the basis of the data available to us, the way in which the people in~$V\!$ have answered~$S$.

So let us look at the question of how likely it is that our unknown person answers most of~$S$ in some specified way~$\tau$.%
   \COMMENT{}
   Knowing his or her answers to some of the questions does not necessarily place us in a better position to guess their remaining answers: if the events that $s$ is specified as $\vs$ or~$\sv$ are independent for different $s\in S$, which in the absence of any knowledge at all is a reasonable assumption, we learn nothing from such partial answers about the remaining questions.

However, if $S$ has tangles then these events are not independent: answers%
   \COMMENT{}
   that occur in the same tangle tend to be%
   \COMMENT{}
   positively correlated. (This inspired our informal description of tangles in \hbox{Chapter}~\secIdea\ as collections of features that `often occur together'.) Therefore, taking as our guess for the answers to~$S$ of our unknown person the specification of~$S$ given by {\it any\/} tangle is already better%
   \Footnote{\dots in the sense that it is more likely to be correct on many questions in~$S$}%
   \COMMENT{}
   than guessing arbitrarily.%
   \COMMENT{} But we can do better still.

For a start, let us restrict our guesses to consistent%
   \COMMENT{}
   specifications of~$S$. This is reasonable: recall that, by definition,%
   \COMMENT{}
   no $v\in V\!$ answered $S$ inconsistently, and so we may assume the same for our unknown person in~$P$. If $T$ in Theorem~\xxxToTfixedS\ is chosen minimal, which we may also assume, then every consistent specification of~$T$ extends to a unique tangle of~$S$.%
   \COMMENT{}
   In other words, once we know how our unknown person has answered~$T$ (consistently), there is exactly one tangle of~$S$ that makes a possible guess for all our person's other answers!
This tangle's set of answers to~$S$, then, is our prediction for how our unknown person would answer~$S$.%
   \COMMENT{}

To summarize, finding the tangles of $S$ gives us a set of better guesses than random: one for each tangle. Theorem~\xxxToTfixedS\ then narrows this down to a unique such guess (with an even better expected success rate) once we are able to quiz our unknown person on~$T$.

In a typical application in the social sciences, $S$~might be a large questionnaire used in a pilot study on a small set $V\!$ of people, and $T$ might be the main study run on~$P$, distilled from the findings for~$S$ on~$V\!$. In a typical application in economics, $T$~might be used as a set of test criteria for properties in~$\vS$ which we desire to know but cannot measure directly.

Let us remark that the prediction discussed above, which answers~$S$ as the unique tangle~$\tau$ of~$S$ does that induces the views of our unknown person $p\in P$ on~$T$, is not the only reasonable prediction for~$p(S)$ if we are in the scenario of Theorem~\xxxToTallSk. In this wider context, the specification $p(T)$ of~$T$ chosen by~$p$ will still extend to a unique maximal tangle in~$S$, but this may be a tangle of low order: a~tangle of $S_k$ for some small~$k$. In such a case we still have to decide what to predict as $p(s)$ for all $s\in S\sm S_k$. Making this choice is related to the problem, discussed in Section~\secThms.4, of how best to associate a given $v\in V$ with some tangle of the entire~$S$.

Quantitative estimates for the quality of these predictions can only be established depending on context. Since feature systems and their tangles are so general, good estimates cannot be expected to apply to worst-case scenarios but will have to be probabilistic. They will have to be based on a probability distribution on the set of all feature systems of a set~$V\!$, but which such distributions are reasonable to assume will depend on the intended application.%
   \COMMENT{}

\beginsection \secFormalExamples. Applying tangles: back to the examples

In this chapter we revisit the scenarios discussed in Chapter~\secInformalExamples\ of how tangles might typically be applied in the social sciences. Having seen a more formal definition of tangles in Chapter~\secFormalSetup\ and statements of the two main tangle theorems in Chapter~\secThms, we are now in a position to say more precisely for each of those scenarios how the tangles discussed there may be defined, and what the two theorems mean in those contexts.

Our aim is that every reader should be able to evaluate, by following these templates, how tangles might be defined in their own professional context, and gauge what their impact might be.%
  \Footnote{As pointed out before, my group and I will be happy to collaborate on any such effort, assist with mathematical or technical support, and provide some relevant software.}

For each of the contexts described in Chapter~\secInformalExamples\ we already discussed there how the sets $S$, $V\!$ and~$\F$ could be specified. All `tangles' mentioned in this chapter will be $\F$-tangles, typically%
   \Footnote{Other choices of~$\F$~-- e.g., meaning-driven choices as in Chapters~\secInformalExamples.1 and~\secFormalSetup.4~-- are possible and may be discussed too where rele\-vant.}
   with $\F=\F_n$ for some integer~$n$.%
   \Footnote{See Chapter~\secFormalSetup.1 for the definition of~$\F_n$.}
   Our aim now is to do at least%
   \COMMENT{}
   some of the following, separately for each context:

\medskip
\plainitem{$\bullet$} describe how Theorems \xxxToTfixedS\ and~\xxxTTD\ can be interpreted;
\plainitem{$\bullet$} discuss choices for an order function as introduced in Chapter~\secFormalSetup.4;
\plainitem{$\bullet$} describe the evolution of the resulting tangles (Ch.~\secFormalSetup.2), and Theorem~\xxxToTallSk;
\plainitem{$\bullet$} interpret any clustering in~$V\!$ based on the tangles of~$S$, as in Chapter~\secThms.4;
\plainitem{$\bullet$} discuss the use of the tangle-distinguishing tree set~$T$;
\plainitem{$\bullet$} compare the tangles of~$S$ with their dual tangles of~$V\!$ from Chapter~\secFormalSetup.5.
\medbreak

In order to avoid repetition, we shall discuss aspects that apply in most or all of these contexts only once, in the mindset example in Section~\secFormalExamples.1. This section, therefore, is recommended reading also for readers primarily interested in another of the fields discussed later.

\subsection \secFormalExamples.1\ \ Sociology: from mindsets to matchmaking

In the simplest%
   \COMMENT{}
   scenario for Chapter~\secInformalExamples.1 we think of $S$ as a questionnaire, which a sample~$V\!$ of a population $P$ of people in whose views we are interested have filled in.%
   \Footnote{To avoid repitition, we subsume under this also the slightly different scenario where $S$ is not a questionnaire but a list of aspects of behaviour that has been observed amongst the people in~$V\!$, and $S$ has been filled in by an observer.}
   The sample~$V\!$ is chosen so as to represent~$P$ well. A~tangle of~$S$ is a mindset or a typical pattern of behaviour, known or unknown.%
   \Footnote{While often we shall be naturally interested in understanding the mindsets found by our theorems and algorithms in terms of some intuitive interpretation, it is important to note that the theory does not require this. In particular, we can use it to make tangle-based predictions for the views or behaviour of people in~$P\sm V$ without such an intuitive understanding. Compare Chapter~\secThms.5.}
   The set~$\F$ may include subsets of~$\vS$ designed to exclude the views of people we are explicitly not interested~in.%
   \COMMENT{}

Theorem~\xxxToTfixedS, our tree-of-tangles theorem for a fixed feature set, finds a tree set~$T$ of questions that are critical for telling apart the mindsets found as tangles. Both the tangles found and the tree set~$T$ that distinguishes them will typically contain not only elements of~$S$ but also Boolean expressions of elements of~$\vS$: certain combinations of particular views on single matters queried by~$S$.\looseness=-1

The fact that $T$ is a tree set is likely to correspond to logical inference. For example, $S$~might have a tangle representing enthusiasm for skiing, and another for boxing. Since these are distinct mindsets, $T$~will have to contain a question that tells them apart. If we assume that no-one likes both skiing and boxing, this might be the question $s$ asking `do you like skiing', or the question $r$ asking `do you like boxing'. If the yes-answers to these questions are $\vs=B\sub V\!$ and $\vr=C$, respectively, we would have $B\sub\Cbar$ and $C\sub\Bbar$.

Suppose now that there are also mindsets that entail a loathing for all sports. Then $T$ must contain questions that distinguish these mindsets from the two sport-loving ones. Neither $r$ nor~$s$ achieves this: $r$~fails to distinguish those sport loathers from the skiers (since these also dislike boxing), while $s$ fails to distinguish them from the boxers. We thus need a further question $t=\{A,\Abar\}$ in~$T$ to do this. This might be the question `do you like any sports', with yes-answer $\vt=A$, say. Since skiing and boxing are types of sport,%
   \Footnote{This is an example of where the fact that $T$ is a tree set reflects logical inference.}
   it will satisfy $A\supe B$ and~$A\supe C$.

The splitting stars of~$T$, see Chapter~\secThms.2, correspond%
   \COMMENT{}
   to the various mindsets that $T$ distinguishes. Each splitting star consists of the logically strongest views in~$\vT$ consistent with the corresponding mindset. It thus represent this mindset `in a nutshell'. If, in the above example, we have one comprehensive mindset of loving sport but hating both skiing and boxing (rather than several mindsets for further types of sport), the splitting star of this mindset would be the set $\{\vt,\sv,\rv\}\sub\vT$. If there is no such mindset, but still no further sports mindsets, then $T$ will not contain all three questions $r,s,t$. It might contain $t$ and~$s$, in which case $\{\vt,\sv\}$ would be the splitting star of the mindset of loving sports but not skiing. This mindset would differ only immaterially from the boxing mindset, since there are not enough sports enthusiasts that dislike both skiing and boxing to form their own mindset. Similarly, $T$~might contain $t$ and~$r$ but not~$s$, in which the skiing fans would be indistinguishable from the more general sports fans that dislike boxing.

In Theorem~\xxxTTD, the splitting stars of~$T$ are in~$\F$, and are therefore each held only by a small set of the people in~$V\!$, or by a set of people deemed irrele\-vant.%
   \COMMENT{}
   The entire set~$V\!$ is then fragmented into these small subsets, which documents conclusively that no significant mindset in terms of~$S$%
   \COMMENT{}
   exists in the population~$V\!$ (and, by inference, in the population~$P$).

In a more subtle setup we have an order function on~$S$ that assigns%
   \Footnote{To ensure that our order function becomes submodular, as required by Theorem~\xxxToTallSk, we can make these assignments indirectly by way of a weight function $w\colon S\to\N$, as explained in Chapter~\secFormalSetup.4.}
   low order to questions of a fundamental or basic nature, and higher order to questions about more detailed aspects of mindsets expressable in terms of those basic views. Indeed, views held by a mindset $\tau$ of~$S_k$ refine those of the mindset $\tau\cap\vSell$ of~$S_\ell$ for each $\ell < k$.

In our earlier example, the question~$t$ about sports in general will probably have lower order than the questions $r$ and~$s$ about boxing and skiing, since it is less specific. We might then have a tangle of low order for sport enthusiasm, and tangles of higher order for the love of skiing and boxing. These latter tangles would both refine the general sportstangle.

The splitting stars of the tangles of lowest order in the tree set~$T$ found by Theorem~\xxxToTallSk\ are the collections of strongest views for the mindsets at this basic level. These mindsets branch out into more refined mindsets of which splitting stars of features of~$T$ of higher order are the sets of strongest views. Theorem~\xxxToTallSk\ thus enables us to extract from the data provided a hierarchy of mindsets of varying degree of complexity, which unfold from the most basic to the most specific. We also have in $T$ a small selection%
   \COMMENT{}
   of key questions on which these mindsets differ, so that we can tell them apart just by evaluating those few questions.%
   \Footnote{However, remember that the questions in~$T$ are not necessarily themselves in~$S$: they may be combinations of questions from~$S$. In order to distinguish all the mindsets on~$S$, it is enough to ask the questions from~$S$ needed to build the composite questions in~$T$.}

Finally, we have the option of not choosing the order of the questions in~$S$ explicitly, but to compute orders of all partitions of~$V\!$ as outlined in Section~\secFormalSetup.4, and then to consider for Theorem~\xxxToTallSk\ the sets $S_1\sub S_2\sub\dots$ arising from this order function. This approach has the advantage that it is entirely deterministic once the questionnaire $S$ has been set up, in that we do not interfere with the result~-- what mindsets are found as tangles~-- by declaring some questions as more basic, fundamental, or important than others. On the other hand, we have less influence on what kinds of possible mindsets are found.

Once we have determined the mindsets in terms of~$S$ that are prevalent in our population~$V\!$ we can, if desired, associate individuals $v\in V\!$ with mindsets that represent their views. As discussed in Chapter~\secThms.4, there are various ways of doing that, and which is best will depend on the context. Note, however, that the more sophisticated methods discussed there, which assign a person~$v$ to a tangle of~$S$ based on his or her entire set $v(S)$ of answers to~$S$, is not available for the larger population~$P$ whose views are known only for the questions in~$T$. For these individuals, the tangle of~$S$ that represents their views best is likely to remain the unique tangle of~$S$ that extends their views $v(T)$ on~$T$, assuming that these are consistent.

The tangle-based approaches to clustering outlined in  Chapter~\secThms.4 can also be used for a matchmaking algorithm as envisaged in Chapter~\secInformalExamples.1, by matching individuals%
   \COMMENT{}
   that have been assigned to the same tangle or to tangles that are  similar in terms of a similarity function~$\sigma$.%
   \COMMENT{}

The $T$-based clustering method mentioned last, for example, could be used to assign two individuals $u,v$ the distance $\sum_{t\in \rho T\tau} (1/|t|)$,%
   \COMMENT{}
   where $\rho T\tau$ is the set of questions $t$ in~$T$ that distinguish the tangle $\rho$ to which $u$ was assigned from the tangle~$\tau$ to which $v$ was assigned.%
   \Footnote{For graph theorists: these $t$ correspond to the edges of the path in the tree~$\TT$ between the node representing~$\rho$ and the node representing~$\tau$.}
   In particular, if $u$ and $v$ were assigned to the same tangle, because they answered the tangle-distinguishing questions from~$T$ identically, they are likely to be matched.

If the tangles $\rho$ and~$\tau$ that $u$ and~$v$ were assigned to are not the same, $u$~and $v$ will still be considered relatively close if they differ on only a few questions from~$T$ and these have relatively high order, i.e., are specific rather than fundamental questions. In the interpretation this means that $u$ and~$v$ are deemed to be close, and are thus likely to be matched, if they share the most fundamental character traits and differ only one some more specific traits that distinguish more refined types of character, while they may differ on many other questions in~$S$ that are irrele\-vant to distinguishing types of character.

\subsection \secFormalExamples.2\ \ Psychology: diagnostics, new syndromes, and the use of duality

In Chapter \secInformalExamples.2 we discussed two possible applications of tangles in psychology.

The first was the formation of idiosyncratic concepts%
   \Footnote{It is important here not to confuse concepts with words used to express them: since language is defined by interpersonal consensus, the individual's concepts and notions we are trying to understand here will typically not be associated with a word.}
   from a plethora of perceptions in the mind of a particular patient: perceptions that combine loosely into a notion in this patient's mind which cannot be naively understood by the therapist, because these perceptions are not commonly associated with each other to form a concept in other people's minds. While I believe that the idea of tangles can inspire the discovery of such idiosyncratic concepts in individual patients, the implementation of this may best be left to the therapist's informed intuition: there will rarely be enough data about our individual patient's perceptions to warrant a computer-based algorithmic search for tangles defining his or her unfamiliar concepts.

Our other application of tangles from Chapter \secInformalExamples.2, however, seems ideally suited to a large-scale computer-based approach. The idea there was to see which psychological symptoms, collected across many patients, combine into recognizable syndromes that can serve as a focus for more targeted research and the development of treatments.

An obvious possible emphasis for a large and fundamental study would lie on discovering previously unknown syndromes. However, any tangles-based definition of `syndrome' will also encompass known syndromes, which should likewise show up in such a study. This can be used in two ways. First, to validate the tangle-based model by checking that known syndromes and psychological conditions are reliably found. Secondly, tangles will offer some quantitative underpinnings also to the ongoing study of previously known syndromes.

The set $T$ of tangle-distinguishing features found by Theorems \xxxToTfixedS\ and~\xxxToTallSk\ on the basis of such a large study can then be used also for diagnostic purposes in individual patients. The idea would be to test~$T$, rather than the much larger~$S$, on a given patient first, and if this returns a consistent specification of~$T$ we have identified a tangle that is most likely to describe the patient's psychological condition: the unique tangle of~$S$ that extends this consistent specification of~$T$.

On the face of it, this procedure resembles the familiar classical diagnostic approach in medicine based on `exclusions': just as a doctor tests specific conditions in order to exclude potential illnesses as the cause of the symptoms observed, testing each potential feature in~$T$ will support some and exclude other possible tangles of~$S$. Indeed, $T$~may well be seen as a particularly efficient set of such tests, which evidence-based guidelines to doctors may suggest.

However we must not forget that each tangle, being a particular specification of~$S$, is just an `ideal' set of test results that corresponds to a syndrome. Diagnosing an individual patient with one of these syndromes is still a nontrivial process: we have to find the syndromes that best match this patient's own test results on~$S$. There are several ways to do this~-- see Chapter~\secThms.4~-- and this should be a matter for learned debate within the discipline as much as for the psychologist in charge of the patient.

Finally, there is an aspect of tangle duality that is particularly useful in a psychology context. Let us again think of $S$ as a set of questions answered by every person $v$ in~$V\!$, even though in our current context the $\vs\in\vS$ are more likely observations made by a psychologist than answers given directly by the client.\looseness=-1

In Chapter~\secIdea.5\ we briefly considered so-called witnessing sets for tangles: a~subset $X$ of~$V\!$ is said to {\it witness\/} a tangle~$\tau$ of~$S$ if, for every $s\in S$, the majority of the elements of~$X$ specify $s$ as $\tau$ does, rather than in the opposite way. If this majority subset of~$X$ has size at least~$p\,|X|$ for every~$s$, where $p\in [{1\over2},1]$, let us call $X$ a {\it $p$-decider\/} for~$S$.

Note that this definition no longer needs or refers to~$\tau$: a set $X\sub V\!$ is a $p$-decider for~$S$ as soon as every $s\in S$ is specified in the same way by at least $p|X|$ of the elements of~$X$. In particular, whether or not a set $X\sub V$ is a $p$-decider for~$S$ is not, as such, a property of tangles. Every $X\sub V$ is trivially a $1\over 2$-decider for~$S$, and every singleton $\{v\}$ is trivially a 1-decider for~$S$, but for given $p > 1/2$ there may be no $p$-decider for~$S$ in~$V\!$ of any non-trivial size. If there is, however, and $p$ is close to~1, then its elements are essentially indistinguishable in terms of~$S$: every question in~$S$ is answered by most of them in the same way.

Now recall from Chapter~\secFormalSetup.5\ that the elements $v\in V\!$ of the ground set of our feature system~$\vS$ are themselves potential features in a dual feature system~$\vV$ with ground set~$S$. \hbox{A~$p$-decider} for the dual system~$\vV\!$, then, is a subset~$S'$ of~$S$ such that every $v\in V\!$ answers most of the questions in~$S'$, at least $p\,|S'|$ many, in the same way: either yes or no.

The questions in~$S'$, therefore, can be used interchangably to some degree, which can be particularly useful in a psychology context. For example, if we are interested in how our patient would answer some particular question $s'\in S'$, but we cannot ask $s'$ directly because of some taboos or traumas, or because asking~$s'$ might influence his or her other answers, we can instead ask some more innocuous $s''\in S'$ and interpret its answer as an answer to~$s'$.

\subsection \secFormalExamples.3\ \ Political science: appointing representative bodies

The starting point for this example was that bodies elected to make decisions on behalf of the population that elected them should be composed of delegates that represent the typical views held in the electorate, perhaps in numbers reflecting the popularity of these views. Given that tangles provide a model for precisely this, views commonly held together, one might seek to constitute a representative body by appointing as delegates people whose views best resemble the tangles found for the set $S$ of rele\-vant issues determining the election.

This problem was addressed, and solved, in Chapter~\secThms.4: we discussed various similarity functions for specifications of~$S$, and since both the tangles of~$S$ and the sets of views on~$S$ of individuals~$v\in V\!$ are specifications of~$S$, we can compare them and delegate to our representative body those $v\in V$ whose views are closest to the various tangles found for~$S$.

Even if all the main views held (in the electorate~$V$) are now represented in our representative body, the question is still whether that representation is fair: more commonly held views should be represented by more delegates.

Theorem~\xxxToTallSk\ provides a neat way to solve this problem too. Recall that the algorithm which computes the tree of tangles for a hierarchy $S_1\sub S_2\sub\dots\sub S$ of feature sets (Chapter~\secThms.3) does so for the sets~$S_k$ in turn, as $k=1,2,\dots$  increases. At each stage~$k$ it evaluates whether or not a tangle found for~$S_k$ extends to one or several tangles of~$S_{k+1}$. At this time, the algorithm has already found the subset $T_k$ of~$T$ needed to distinguish all the tangles in~$S$ of order up to~$k$. The nested set~$T_k$ partitions~$V\!$ into the sets $V_\tau := \bigcap \,\{\, A\mid A\in \tau\cap\vT_k\}$, where $\tau$ varies over the consistent specifications (or tangles) of~$T_k$.%
   \COMMENT{}

At this stage we could establish our representative body by appointing for each~$\tau$ a number $d_\tau$ of delegates proportional to the number~$|V_\tau|$ of people whose views are best represented by~$\tau$.

However, Theorem~\xxxToTallSk\ enables us to improve on this still, by chosing the composition of the group of delegates for~$\tau$ not arbitrarily (specifying only its size) but according to the tangles of higher order into which $\tau$ `splits': the tangles of $S_{k+1}, S_{k+2},\dots$ that induce~$\tau$. In other words, when we run our algorithm we make a decision at stage~$k$ for every $\tau$ as above whether or not we seek to refine $T_k$ at the node representing~$\tau$: for those~$\tau$ whose $|V_\tau|$ is smallest we appoint one delegate, while for the $\tau$ with bigger~$|V_\tau|$ we keep the algorithm running.

Even then we may end up with sets~$V_\tau$ that are bigger than others: it can happen that a tangle $\tau$ is represented by many people but still does not split into higher-order tangles. So we shall still pick possibly different numbers~$d_\tau$ of representatives for different tangles~$\tau$. But this procedure uses such arbitrary choices of representatives only as a last resort, when no alternative choices informed by higher-order tangles can be used instead.

Unlike in a party-based representative vote, the composition of delegates obtained in this way reflects not only the relative support of the major parties in the electorate, but also within each of these blocks its various shades of more refined views, all proportionately to how these views live in the electorate.

\subsection \secFormalExamples.4\ \ Education: devising and assigning students to methods

In Chapter~\secInformalExamples.4 we discussed how tangles can help to group teaching techniques into `methods' so that different types of students can benefit best from techniques that work particularly well for them. A~`method' in this sense is simply a collection of techniques that are applied at the same time in one class. The idea is that, since at a large school it may be possible or necessary to offer the same content in several parallel classes, these classes could use different such `methods' to benefit different types of students.

Our approach was to define `methods' as tangles of techniques: groups of techniques that often work well together, where `often' means `for many potential students'.

Finding these tangles will require a fundamental study involving a large set~$V\!$ of students. Every student in such a study will have to be tested against various techniques, say those in a set~$S$, which requires a considerable effort over time.%
   \Footnote{Older students might, of course, simply be interviewed and asked which teaching methods they prefer. This may be less significant, but it may help make such a study feasible.}
   Once these tests have been performed we can define an order function on~$S$ to distinguish techniques that made a big impact (positive or negative) on many of the students tested: such techniques will be assigned a low order to ensure that they are considered even for the most basic tangles computed. If we wish not only to find tangles of~$S$ but to apply our theorems to them, we could view the $s\in S$ as partitions $\{A,B\}$ of~$V\!$ and define their order $|A,B|$ indirectly through weighted similarity functions, as in Chapter~\secFormalSetup.4, giving greater weight to the techniques found to be more divisive in our tests.

Either way, the $k$-tangles found will be concrete specifications of some techniques, ranging from the fundamental to the more subtle as $k$ increases. They each form a `method' in that they lay down which of the techniques they specify should be used or not when teaching `according to' that tangle.

All this cannot be done at a single school facing the practical problem of how to set up its classes and assign its students to these classes. Conversely, not every `solution' found by a large study, such as a set of tangles and a tangle-distinguishing set $T$ of critical techniques,%
   \COMMENT{}
   can be implemented at any given school. Indeed, local requirements or possibilities may dictate how many classes can be taught in parallel (and hence, into many tangles or `methods' the techniques available should be grouped), and which techniques are available may depend on the composition of the local teaching staff.

However, once we have the data from our large study, a~complete specification $v(S)$ of~$S$ for every student~$v$ tested, we can run our tangle algorithms with different parameters, parameters than can be adjusted to a school's given requirements or possibilities. For example, we can adjust the parameter~$n$ in the definition of~$\F_n$ so as to obtain a desired number of tangles to match the target number of parallel classes, or we can delete potential techniques $s$ from~$S$ before computing the tangles if these techniques cannot be implemented at the given school.

The result of this process, then, is a tangle-based concrete way to group the techniques available into as many `methods' as parallel classes can be offered.

But how do we assign students to these classes? One possibility is to simply describe the methods on which the various parallel classes are based, and let them choose. Another is to test, or at least ask, students which techniques from the specific set $T$ found by the algorithm for Theorem~\xxxToTfixedS\ works best for them: by definition, this is a small set of `critical' techniques, each a combination of the techniques tested directly (those in~$S$), on which the methods found differ. Every (consistent) specification of~$T$ defines a unique tangle that extends it~-- one of the teaching methods that correspond to the classes offered in parallel.

As with our earlier task of determining a set of methods to be offered, the set~$T$ that critically distinguishes them can be computed locally to match the concrete school's possibilities and requirements. In this way, the results of the same large study can be used by different schools to compute a bespoke set of particularly critical teaching techniques against which students can be tested in order to assign them to the classes corresponding to the various methods.

As the number of tangles involved in these applications will be small, the way to do this will be to assign each student $v$ to the class whose tangle~$\tau$ extends his or her specification $v(T)$ of~$T$, the class $V_\tau := \bigcap \,\{\, A\mid A\in \tau\cap\vT\,\}$ defined after Theorem~\xxxToTfixedS\ in Chapter~\secThms.1.%
   \COMMENT{}

\subsection \secFormalExamples.5\ \ Meaning: from philosopy to machine learning

Our notions%
   \COMMENT{}
   determine how we see the world. They are our way of cutting the continuum of our perceptions and ideas into recognizable chunks which, to some extent, persist in time, space, and across different people. Chunks of ideas or perceptions which we bundle together again and again, in different places, and which are similar to the recurring chunks in other people's minds.

Inasmuch as%
   \COMMENT{}
   we understand, remember, and communicate aspects of the world around us as structures of notions, the question of what these notions are, i.e., which ideas or perceptions we combine into notions, determines what we can understand, remember, or communicate about the world. Since in everyday life%
   \Footnote{Unlike, for example, in mathematics: there, decisions about which properties of the (mathematical) objects studied should be bundled into notions to facilitate further study are part of our daily bread. And they are as important as elsewhere in life, since they determine what structures of the objects under study become visible and can therefore be explored.}
   we do not choose our notions consciously, understanding and quantifying them~-- as tangles enable us to do~-- is perhaps significant for our own understanding of the world only when we seek to compare how people from different cultures or backgrounds understand the world differently.

But as soon as we try to teach a computer what our notions are, perhaps through some interactive process,%
  \Footnote{For example, the computer might learn by showing us pictures of groups of objects and asking which of them `do not belong' in this group.}
   we need some quantitative definition of a `notion'. As we have already seen, tangles may well offer that.

More ambitiously, if we seek to enable computers to `understand' the world by themselves, the question of what should be their own fundamental notions, not necessarily copying ours, will be on the agenda even more. This is because there are `good' notions and `bad' ones, judged by how they help us understand the world.%
   \COMMENT{}
   Good notions cut the continuum of ideas and perceptions along lines running between phenomena which we would like to distinguish,%
   \Footnote{Despite the wording here, this is not entirely a matter of taste. Very crudely, one might say that good notions facilitate the expression of theories that describe the world better than others~-- for example, in the sense of making more specific or substantial predictions.}
   whereas bad notions cut across such phenomena and are therefore less helpful for distinguishing them. As we saw in Chapter~\secFormalSetup, this aim of finding good notions, possibly unlike ours, is very similar to finding tangles seen as specifications of set partitions.

We discussed in Chapters \secIdea.3 and~\secInformalExamples.5 how the notion of `chair' can be captured as a tangle of potential features of furniture.  Note that we can choose to be more demanding regarding the quality of this notion by raising the parameter~$n$ in the definition of~$\F_n$-tangles, since a furniture type will show up as an $\F_n$-tangle $\tau$ only if every three of the features in~$\tau$ are shared by at least~$n$ items. Another way to be more demanding of an $\F_n$-tangle, in order to ensure it only identifies high-quality notions, would be to replace `three' with some higher number: to ask that not only all triples of features in tangle are shared by many items but, perhaps, all quintuples.%
   \Footnote{If we are too demanding, of course, we shall simply not get enough tangles to serve as notions, so we have to strike a balance here. But the quantitative parameters of $\F_n$-tangles allow for precisely that~-- not to mention the numerous possibilities of other choices of~$\F$.}

For the sake of our chair example, we considered as features only properties of furniture, such as being made of wood, which may or may not apply to a particular piece at hand. However when we try to train a computer to form or recognize notions, there is no reason to be so restrictive: we can also use other parameters that help us distinguish between different notions, such as context, or the time when it was fashionable.%
   \COMMENT{}%
   \Footnote{Imagine a party game where a person thinks of something and we have to guess what it is by asking yes-no questions of any kind. Every such question is a `potential feature'.}

In Chapter~\secFormalSetup.4 we saw how to use order functions to divide potential features into hierarchies, assigning low order to basic features and higher order to more specific ones. This is nowhere more rele\-vant than when we employ tangles to capture notions in our ideas and perceptions: some notions are clearly more fundamental than others, and some questions (potential features) are more basic than other questions and thus should have lower order.

In addition to, or as a basis for, assigning an order to potential features in terms of similarity functions for set partitions (see Chapter~\secFormalSetup.4), there are two ways of defining order explicitly that come to mind in the context of meaning: one is relevance, the other what one might call clarity or definiteness.

As an example for a relevance-based order function we could choose to assign all colour questions high order in our search for furniture types if we consider colour irrele\-vant to how furniture splits into types. As an example for a clarity-based order function, suppose we are trying to distinguish hairstyles. We are likely to get more helpful answers in the clarity sense (of delineating hairstyles from each other) if we ask when they were fashionable than whether they are pretty, so the former question should receive lower order on the clarity scale than the latter.

Note that these two considerations for how to define an order function may well conflict: being pretty or not is perhaps the most rele\-vant aspect of a hairstyle, but it may also be the least clear in that people have divided opinions about it. Note also that even if we wish to base our order function on relevance or clarity or both, we do not necessarily have to define the order~$|s|$, or a weight~$w(s)$, manually for every question~$s$: both relevance and clarity can be gleaned to some extent from how the people who took our survey answered its questions, and therefore computed mechanically.%
  \COMMENT{}

What we have described so far is how a computer can learn the meaning of words as we use them, or come up with notions of its own that are formed by observing the world. Both these are aspects of what one would call the formation of passive vocabulary. Once this has been achieved, our computer will also want to know which of its tangle-encoded notions best describe a given object presented to it. This is a problem we discussed in Chapter \secThms.4: the problem of how to match a given object to a tangle of potential features, the tangle that best captures its actual features. See there for further discussion.

Let us complete this section with a straightforward application of tangles of notions, just as an example of what may be possible: an interactive thesaurus for non-native speakers.

Roget's classical {\it Thesaurus\/} helps writers find the best word for what they are trying to say simply by grouping together some likely candidates. The rele\-vant group of words can be found by looking up a word whose meaning is close to the intended meaning, but maybe does not quite capture it. If a word has multiple meanings, it is linked to several such groups.

The value of this lies in offering the writer a rele\-vant choice, but the thesaurus does not help with making this choice. This is fine if the writer knows all the words on offer, and perhaps just could not think of the right one. For learners of a language, however, this falls short of their need: they, too, will know what they are trying to express, but need help in finding the best word for it. Let us see how tangles can help them.

In the simplest model, we could devise for every word field currently offered together as a group of choices a questionnaire~$S$ whose answers for any given context would enable us to choose the correct word from this field. Each of these words, then, is likely to correspond to a tangle of~$S$~-- a~tangle we would be able to compute before our thesaurus is published. We could then also compute the tangle-distinguishing small set $T$ of questions from Theorem~\xxxToTfixedS. These questions from~$T$ could then be put to the user trying to identify the right word for their intended notion, and answering just these will steer them to the word that best fits their intended notion.

Recall that the questions in~$T$ may be combinations of questions from~$S$, so the questions a user really has to answer may be a little larger. But the subset of~$S$ needed to form the questions in~$T$ is still likely to be smaller than $S$ itself. This is crucial for making a good thesaurus: it will matter to the user whether they need five questions or fifteen to be steered to the fitting word.

Devising a questionnaire~$S$ for each word field in the thesaurus, and answering all its questions for each word, may look like a lot of work. But this work has to be carried out only once, when the thesaurus is made: it is offset by a gain on the user side in that the questions asked are chosen specifically for each word field selected by the user, and are chosen particularly well for this word field.

In a more sophisticated model, we could grade the questions from~$S$ as more or less rele\-vant for the corresponding word field, and assign them an order correspondingly. We would then obtain a hierarchy of tangles as described in Chapter~\secFormalSetup.2. This would correspond to a hierarchy of more or less general or specific words, just as in reality.

In this last respect, our thesaurus might even learn from user interaction. Although a user may not know the meanings of the words on offer, they will be able to grade the questions put to them as more or less rele\-vant to their specific search: tangles, after all, reflect notions, not words, and users come with such notions in their minds. Our tangles will then have to be recomputed from time to time when enough user feedback has been collected, and re-checked editorially against the words on offer. At this point, editors might also comment on user-defined tangles that do {\it not\/} have a corresponding word to match: such tangles will exist, since notions exist that are not exactly matched by words. Compared against the notions in the minds of speakers of other languages, however, this seems even more likely and worth addressing in ongoing editorial work.

\subsection \secFormalExamples.6\ \ Economics: customer and product types

Let us rephrase our example from Chapter~\secInformalExamples.6 of an online shop, and the tangles of the dual feature systems this gives rise to, in terms of set partitions as introduced in  Chapter~\secFormalSetup.1.

In the example, $V\!$~was a set of customers of an online shop, and $S$ the set of items sold at this shop. We assumed that each customer~$v$ has made a single visit to the shop, specifying $s\in S$ as $v(s)=\vs$ if he or she included~$s$ in their purchase, and as $v(s)=\sv$ if not.%
   \Footnote{Note that, unlike in most other contexts, every $s\in S$ has a default specification here: we use~$\vs$ to indicate that $s$ is `present' rather than `absent', consistently for all $s\in S$.}
   In the dual setup, every item $s\in S$ specifies every~$v\in V$: as $s(v) = \vv$ if $v$ bought~$s$, and as~$\vvback$ otherwise. Thus, $s(v)=\vv$ if and only if $v(s)=\vs$.

In the terminology of Chapter~\secFormalSetup.1, every item $s\in S$ is a partition of~$V$: the partition $s = \{\vs,\sv\}$ of~$V\!$ into the complementary set
 $$\vs = \{\,v\in V\mid v \hbox{ bought } s\,\}$$%
   \COMMENT{}
 of customers that bought~$s$ and the set $\sv = V\sm\vs$ of those that did not. Similarly, every customer $v\in V\!$ is a partition of~$S$: the partition $v = \{\vv,\vvback\}$ of~$S$ into the complementary set
 $$\vv = \{\,s\in S\mid v \hbox{ bought } s\,\}$$
 of items that $v$ bought and the set $\vvback = S\sm\vv$ of items that $v$~did not buy.

We thought of the specification $v(S)$ of~$S$ as $v$'s `shopping basket', the division of~$S$ into the set $\vv$ of items that $v$ bought and the set~$\vvback$ of items that $v$ did not buy. Similarly, we thought of $s(V)$ as the `popularity score' of~$s$, the division of~$V\!$ into the set $\vs$ of fans of~$s$ and the set~$\sv$ of non-fans of~$s$.

The tangles of~$S$, then, were typical shopping baskets, or types of customer: ways of splitting the set $S$ of items into two parts~-- bought versus not bought~-- that were reflected, by and large, by the purchasing behaviour of many customers.%
   \COMMENT{}
   The tangles of~$V\!$ were typical popularity scores, or types of item: ways of splitting the set of customers into fans and non-fans which broadly reflects the popularity score of many items.

More formally, an $\F_n$-tangle $\tau$ of~$S$ is a specification of~$S$ such that for every three elements $s_1,s_2,s_3$ of~$S$ there are at least $n$ customers~$v$ that bought or failed to buy each~$s_i$ as prescribed by~$\tau$: such that $v(s_i) = \tau(s_i)$ for all~$i$, or more explicitly, such that $v$ bought~$s_i$ if $\tau(s_i) = \vsi$ and $v$ did not buy~$s_i$ if $\tau(s_i) = \svi$.
Similarly, an $\F_n$-tangle~$\rho$ of~$V\!$ is a specification of~$V\!$ such that for every three customers $v_1,v_2,v_3$ in~$V\!$ there are at least $n$ items~$s$ that were bought or not by each~$v_i$ as prescribed by~$\rho\,$: such that $s(v_i) = \rho(v_i)$ for all~$i$, or more explicitly, such that $v_i$ bought~$s$ if $\rho(v_i) = \vvi$ and $v_i$ did not buy~$s$ if $\rho(v_i) = \vviback$.

A~subset $U$ of~$V\!$ witnesses a tangle $\tau$ of~$S$ if, for every single $s\in S$, the tastes of most customers in~$U$ regarding~$s$ is as described by~$\tau$: if most customers in~$U$ bought~$s$ when $\tau(s)=\vs$, and most customers in~$U$ did not buy~$s$ when $\tau(s)=\sv$.\penalty-200\
   A subset $R$ of~$S$ witnesses a tangle $\rho$ of~$V\!$ if every $v\in V\!$ bought or avoided most of the items in~$R$ as described by~$\rho$: if $v$ bought most items in~$R$ when $\rho(v)=\vv$, and $v$ avoided most items in~$R$ if $\rho(v)=\vvback$.

With this formal setup in place, let us now revisit the examples discussed in Chapter~\secInformalExamples.6.

The clearest of these is the tangle~$\rho$ of~$V\!$ which specifies the rich customers~$v$ as~$\vv$, and all the others as~$\vvback$. This is a tangle, because there is a sizable group of items which the rich customers will tend to like {\it and\/} all the other customers will tend to dislike: the expensive high-quality items. Formally, we shall find for every triple $v_1,v_2,v_3$ of customers $n$ such items~$s$ to satisfy $s(v_i)=\rho(v_i)$ for all~$i=1,2,3\,$: so that, regardless of whether $v_i$ is rich or not, their purchasing behaviour with respect to our $n$ expensive high-quality items~$s$ will be as prescribed by~$\rho$.%
   \Footnote{Note, however, that the expensive high-quality items do not witness the tangle of rich people in the sense of Chapter~\secIdea.5, because not every rich person will have bought, in their one purchase, more than half of all the expensive items.}

If we assume that all the green items in our shop are more expensive than competing non-green items, there is a similar tangle of~$V\!$ which specifies a customer~$v$ as~$\vv$ whenever $v$ likes ecological items (even if they are expensive), and as~$\vvback$ if $v$ does not care enough about ecology to pay more. The green items in~$S$ then witness this tangle~$\rho$ of~$V\!$ in the same weak sense as above, i.e., in terms of the definition of an $\F_n$-tangle: for every triple $v_1,v_2,v_3$ in~$V\!$ there are at least~$n$ items, all among the green ones, which each $v_i$ buys if $\rho(v_i)=\vvi$ and does not buy if~$\rho(v_i) = \vviback$.%
   \COMMENT{}

Let us now change our assumptions slightly, and assume no longer that the green items are more expensive. Then $\rho$ is no longer a tangle of~$V\!$, since non-green customers might buy green items too, even if they do not prefer them: for a triple of non-green customers we may be unable to find $n$ items that they all fail to buy.%
   \COMMENT{}

Let us now look at tangles of~$S$. If we assume that every $s\in S$ is either green or positively dangerous to the environment, we have a tangle of~$S$ witnessed by the ecology-minded customers, again in the weak sense of the definition of $\F_n$-tangles. If we do not make this assumption, however, we may have no tangle of~$S$, regardless of whether the green goods are more expensive or not: the ecology-minded customers will no longer witness such a tangle, because they will disagree too much on the goods that are not ecologically critical. More fundamentally, we may be unable to extend the specification (as~$\vs$) of the green products~$s$ to a tangle of all of~$S$, because we may be unable to specify the environmentally neutral items in such a way that for every three of them we can find $n$ customers whose purchases agree on these three items with our specification. In this case, therefore, we would only have a tangle of the set $S'\sub S$ of environmentally critical items, those that are either good or bad for the environment.

Is there a way to capture the preference for green items by a tangle even in this more realistic case that not all items in our shop are environmentally critical? Instead of looking for tangles of all of~$S$ we might try to define an order function on~$S$ so that, for example, the set $S'$ above appears as one of the sets~$S_k$. Then there would be a `green' $k$-tangle witnessed by the ecology-minded customers. However it will be difficult or impossible to find an order function that singles out all the rele\-vant subsets of~$S$ in this way~-- not to mention the fact that we want to be open to unknown tangles, for which we cannot hope to encode a suitable subset of~$S$ in the form of~$S_k$.

However, the approach sketched in Chapter~\secFormalSetup.4 offers a way forward. Let us consider the set~$U$ of all partitions%
   \Footnote{Recall that {\it partitions\/} of a set, in our terminology, always split it into exactly two parts.}
   of~$S$ and define the {\it order\/} of ${\{A,B\}\in U}$ as
 $$|A,B| := \sum_{a\in A}\sum_{b\in B} \sigma (a,b)\,,$$\noskip
 where
 $$\sigma(a,b) = \big|\{\,v\in V\! : a(v)=b(v)=\vv\}\big|$$
 is the number of customers that bought both $a$ and~$b$. Thus, very broadly, a partition of~$S$ gets high order, and is therefore disregarded when $k$-tangles are computed for small~$k$, if it divides the contents~$\vv$ of many shopping baskets~$v(S)$~-- not the items missing from them~-- into two roughly equal parts.

Unlike our earlier idea, this unifies tangles of different subsets of~$S$ into tangles of the same set~$U$. For example, the set $S'$ of green items in~$S$ will witness a tangle of~$U$: partitions of low order cannot divide $S'$ into roughly equal parts, so most of $S'$ will lie on the same side of~-- and thereby consistently specify~-- all such partitions.

Similarly, the set of the comparatively inexpensive items in~$S$ witnesses a tangle of~$U$ (thought not a tangle of~$V$).%
   \Footnote{It witnesses a tangle of~$U$, because the pairs of inexpensive items will be bought together by many price-conscious customers, so no low-order partition of~$S$ will divide many pairs of inexpensive items, i.e.\ most of them lie on the same side of every partition. But as with the ecological items, the set of inexpensive items does not witness a tangle of~$V\!$ that specifies precisely the price-conscious customers~$v$ as~$\vv$, even in the weak sense of providing the $n$ elements needed for its triples, since customers that are not price-conscious have no motivation to avoid inexpensive items.}%
   \COMMENT{}
   By contrast, the set of expensive items in~$S$ does not witness a tangle of~$U$ (let alone one of~$V$): customers unable to afford them will not contribute to the value of~$\sigma$ for pairs of expensive items,%
   \COMMENT{}
   while customers that can afford them are not necessarily likely to buy them more often than cheaper items. Hence, $\sigma(a,b)$~will not be greater for pairs $a,b$ of expensive items than for arbitrary pairs of items.

This pair of examples shows neatly that tangles of~$U$ deliver what we hoped for: they describe, loosely and without any need for difficult judgments in single cases, groupings of items that provide some common motivation for purchase. In particular, they will do this for previously unknown kinds of motivation such as, perhaps, featuring a face on the packaging.

Conversely, tangles of~$U$ avoid some silly instances that can occur as tangles of $S$ or~$V\!$. For example, our online shop $S$~has a particularly silly tangle: the empty shopping basket $\{\,\sv\mid s\in S\,\}$. Indeed, if $S$ offers any reasonable amount of choice, this is a tangle witnessed by all of~$V\!$, even in the strong sense of Chapter~\secIdea.5: just because our shop has so many more items to offer than a single customer can buy, it is likely that for every item~$s\in S$ there are many more customers~$v$ that did not buy~$s$ (and hence specify it as $v(s):=\sv$) than there are customers that bought~$s$. Thus, $V\!$~is a witnessing set for the empty shopping basket as a tangle of~$S$.

However, one can use Theorem~\xxxTTD\ to show that there is no such tangle of~$U$; see~\cite{TanglesEmpirical,MonaLisa}.%
   \COMMENT{}

\beginsection References

\bibliographystyle{plain}
\bibliography{collective.bib}

\bigskip\vfill
\ninepoint\noindent Version 17 July, 2019

\eject\end